\documentclass[12pt,a4paper]{article}

\usepackage[utf8]{inputenc}
\usepackage[english]{babel}
\usepackage{amsmath,amssymb,amsfonts,amsthm,graphicx}
\usepackage{hyperref}
\usepackage{cases}

\begin{document}

\newcommand{\C}{\mathbb {C}}
\newcommand{\N}{\mathbb{N}}
\newcommand{\R}{\mathbb{R}}
\newcommand{\tf}{\mathcal{F}}

\renewcommand{\theenumi}{\roman{enumi}}
\renewcommand{\labelenumi}{\theenumi)}

\swapnumbers
\newtheorem{THM0.1}{Theorem}[section]
\newtheorem{THM0.2}[THM0.1]{Theorem}
\newtheorem{THM0.3}[THM0.1]{Theorem}
\newtheorem{THM0.4}[THM0.1]{Theorem}
\newtheorem{REM1}{Remark}[section]
\newtheorem{DEF1}[REM1]{Definition}
\newtheorem{THM1}[REM1]{Theorem}
\newtheorem{LEM6}[REM1]{Lemma}
\newtheorem{REM4}[REM1]{Remark}
\newtheorem{LEM7}[REM1]{Lemma}
\newtheorem{THM11}[REM1]{Theorem}
\newtheorem{DEF2}{Definition}[section]
\newtheorem{THM2}[DEF2]{Theorem}
\newtheorem{REM7}[DEF2]{Remark}
\newtheorem{THM3}[DEF2]{Theorem}
\newtheorem{THM4}{Theorem}[section]
\newtheorem{REM3}[THM4]{Remark}
\newtheorem{COR2}[THM4]{Corollary}
\newtheorem{REM6}[THM4]{Remark}
\newtheorem{THM5}[THM4]{Theorem}
\newtheorem{DEF3}{Definition}[section]
\newtheorem{REM5}[DEF3]{Remarks}
\newtheorem{THM6}[DEF3]{Theorem}
\newtheorem{REM0}[DEF3]{Remark}
\newtheorem{THM7}[DEF3]{Theorem}
\newtheorem{COR3}[DEF3]{Corollary}
\newtheorem{THM8}[DEF3]{Theorem}
\newtheorem{PROP1}[DEF3]{Proposition}
\newtheorem{THM9}[DEF3]{Theorem}
\newtheorem{THM10}[DEF3]{Theorem}
\newtheorem{LEM1}{Proposition}[section]
\newtheorem{LEM2}[LEM1]{Proposition}
\newtheorem{LEM3}[LEM1]{Lemma}
\newtheorem{LEM4}[LEM1]{Theorem}
\newtheorem{LEM5}[LEM1]{Theorem}
\newtheorem{COR1}[LEM1]{Corollary}
\newtheorem{REM2}[LEM1]{Remark}

\title{Lossless error estimates for the stationary phase method with applications to propagation features for the Schrödinger equation}

\author{Felix Ali Mehmeti\footnote{Université de Valenciennes et du Hainaut-Cambrésis, LAMAV, FR CNRS 2956, Le Mont Houy, 59313 Valenciennes Cedex 9, France. Email: felix.ali-mehmeti@univ-valenciennes.fr} , Florent Dewez\footnote{Université Lille 1, Laboratoire Paul Painlevé, CNRS U.M.R 8524, 59655 Villeneuve d'Ascq Cedex, France. Email: florent.dewez@math.univ-lille1.fr}}

\date{}

\maketitle

\begin{abstract}
	We consider a version of the stationary phase method in one dimension of A. Erdélyi, allowing the phase to have stationary points of non-integer order and the amplitude to have integrable singularities. After having completed the original proof and improved the error estimate in the case of regular amplitude, we consider a modification of the method by replacing the smooth cut-off function employed in the source by a characteristic function, leading to more precise remainder estimates. We exploit this refinement to study the time-asymptotic behaviour of the solution of the free Schrödinger equation on the line, where the Fourier transform of the initial data is compactly supported and has a singularity. We obtain asymptotic expansions with respect to time in certain space-time cones as well as uniform and optimal estimates in curved regions which are asymptotically larger than any space-time cone. These results show the influence of the frequency band and of the singularity on the propagation and on the decay of the wave packets.

\end{abstract}

\vspace{0.3cm}



\setcounter{section}{-1}
\section{Introduction}

\hspace{2.5ex} The asymptotic behaviour of oscillatory integrals with respect to a large parameter, sometimes used to study long-time asymptotics for solutions of dispersive equations, can often be described using the stationary phase method. A theorem of A. Erdélyi \cite[section 2.8]{erdelyi} permits to treat oscillatory integrals with singular amplitudes and furnishes asymptotic expansions with explicit remainder estimates. The approach is specific for one integration variable and the results are interesting for applications. Unfortunately the proof is only sketched in the source \cite{erdelyi}. In the present paper, we start by providing a complete proof and by improving the remainder estimates in the case without amplitude singularities. Then by applying the above method to a particular example, we exhibit an inherent blow-up of the expansion occurring when the endpoints of the integration interval tend to each other. In particular, we remark that the smooth cut-off function, employed in the original proof, prevents us from controlling explicitly the blow-up, restricting potentially the field of applications. This motivates an improvement of the above stationary phase method which consists in replacing the smooth cut-off function by a characteristic function, making the blow-up explicit in the applications. Finally we apply these abstract results to the solution of the free Schrödinger equation on the line for initial conditions having a singular Fourier transform with support in a compact interval. We calculate expansions to one term with respect to time in certain space-time cones, exhibiting the optimal time decay in these regions. Moreover by exploiting the above mentioned refinement of the stationary phase method, we provide uniform estimates of the solution in curved regions which are asymptotically larger than any space-time cone. We prove that the resulting decay rates are optimal by expanding the solution on the boundaries of these regions. These results lead to an asymptotic localization of the wave packets and highlight the influence of the singularity on the time decay rate.\\

Consider the free Schrödinger equation
\begin{equation*}
	(S) \qquad \left\{ \hspace{-3mm} \begin{array}{rl}
			& \big[ i \partial_t + \partial_x^2 \big] u(t,x) = 0 \\
			& \vspace{-0.3cm} \\
			& u(0,x) = u_0(x)
	\end{array} \right. \; ,
\end{equation*}
for $t > 0$ and $x \in \R$. If we suppose $u_0 \in L^{1}(\R)$ then
\begin{equation*}
	\left\| u(t,.) \right\|_{L^{\infty}(\R)} \leqslant \frac{\| u_0 \|_{L^1(\R)}}{2\sqrt{\pi}} \, t^{-\frac{1}{2}} \; ,
\end{equation*}
see for example \cite[p.60]{reed-simon}. If it is assumed that $u_0 \in L^2(\R)$ then we have by Strichartz' estimate (\cite{strichartz}, see also \cite{banica}),
\begin{equation*}
	\left\| u(t,.) \right\|_{L^{\infty}(\R)} \leqslant C \, \| u_0 \|_{L^2(\R)} \, t^{-\frac{1}{4}} \; ,
\end{equation*}
where $C>0$ is a certain constant.

Now consider the following class of initial data\\

\noindent \textbf{Condition ($\mathbf{C_{[p_1,p_2],\mu}}$).} \emph{Let $\mu \in (0,1)$ and $p_1 < p_2$ be two finite real numbers.\\
	A tempered distribution $u_0$ satisfies Condition (C$_{[p_1,p_2],\mu}$) if and only if $\tf u_0$ is a function which verifies $\tf u_0 \equiv 0$ on $\R \setminus [p_1,p_2)$ and
	\begin{equation*}
		\forall \, p \in (p_1, p_2) \qquad \tf u_0(p) = (p-p_1)^{\mu-1}	\, \tilde{u}(p) \; ,
	\end{equation*}
	where $\tilde{u} \in \mathcal{C}^1 \big( [p_1,p_2], \C \big)$ and $\tilde{u}(p_1) \neq 0$.}\\
	
\noindent Here $\tf u_0$ refers to the Fourier transform of $u_0$. Under this assumption, $u_0$ is a smooth function which never belongs to $L^1(\R)$ and belongs to $L^2(\R)$ if and only if $\mu \in \big(\frac{1}{2}, 1\big)$. The question of the $L^{\infty}$-time decay rate for the above problem when $\mu \in \big(0,\frac{1}{2}\big)$ has been answered in \cite[Theorem 2.5]{dewez}. There the author furnishes the following $L^{\infty}$-norm of the solution,
\begin{equation*}
	\left\| u(t,.) \right\|_{L^{\infty}(\R)} \leqslant c(u_0) \, t^{-\frac{\mu}{2}} \; ,
\end{equation*}
the constant $c(u_0) \geqslant 0$ is explicitly given in the statement of the theorem. In the present paper, we exploit a rewriting of the solution as an oscillatory integral in order to apply the abstract results established in the theorems \ref{THM1}, \ref{THM11}, \ref{THM2} and \ref{THM3}. Thanks to that, we prove that the above decay rate is optimal since it is attained on a well-chosen space-time direction. To do so, we use these abstract results to expand the solution on the critical direction $\frac{x}{t} = 2 \, p_1$, where the stationary point and the singularity coincide :

\begin{THM0.1} \label{THM0.1}
	Suppose that $u_0$ satisfies Condition \emph{(C$_{[p_1,p_2],\mu}$)}. For all $t > 0$, define $L_{\mu}(t,u_0)$ as follows :
	\begin{equation*}
		 \bullet \quad L_{\mu}(t,u_0) := \frac{1}{2} \, \Gamma \hspace{-1mm} \left( \frac{\mu}{2} \right) e^{-i \frac{\pi \mu}{4}} e^{i t p_1^{\, 2}} \, \tilde{u}(p_1) \; .
	\end{equation*}
	Then for all $(t,x) \in (0,+\infty) \times \R$ such that $\frac{x}{t} = 2 \, p_1$,	we have
	\begin{equation*}
		\left| u(t,x) - L_{\mu}(t,u_0) \, t^{-\frac{\mu}{2}} \right| \leqslant C(u_0) \, t^{-\frac{1}{2}} \; .
	\end{equation*}
	The constant $C(u_0) \geqslant 0$ is given in the proof.
\end{THM0.1}

\noindent See Theorem \ref{THM8} for a complete statement. Theorem \ref{THM0.1}, coupled with the paper \cite{dewez}, permits to construct in a natural way examples exhibiting the optimal decay rate $t^{-\frac{\mu}{2}}$ with $\mu \in (0,1)$. In particular for a given decay rate $t^{-\alpha}$ with $\alpha \in \big( \frac{1}{4} , \frac{1}{2} \big)$, we can furnish an initial data in $L^2(\R)$ such that the associated solution has an optimal $L^{\infty}$-time decay rate given by $t^{-\alpha}$.

In the present paper, we show also that the singular frequency affects the decay rate not only on the critical direction but also in other space-time regions. Especially, thanks to the preciseness coming from the refinement of the stationary phase method of Erdélyi, we establish a result which suggests that the influence of the singularity seems to be stronger in space-time regions along the direction $p_1$ :

\begin{THM0.2} \label{THM0.2}
	Suppose that $u_0$ satisfies Condition \emph{(C$_{[p_1,p_2],\mu}$)} with $\mu < \frac{1}{2}$, and fix $\varepsilon \in \big(0, \frac{1}{2}\big)$. Then for all $(t,x) \in (0,+\infty) \times \R$ satisfying $2 \, p_1 + 2 \, t^{-\varepsilon} \leqslant \frac{x}{t} < 2 \, p_2$ and $t > (p_2 - p_1)^{-\frac{1}{\varepsilon}}$, we have
	\begin{equation*}
		\big| u(t,x) \big| \leqslant C_0(u_0) \, t^{-\mu + \varepsilon \mu} + \sum_{k=1}^9 C_k(u_0) \, t^{- \sigma_k^{(1)} + \varepsilon \varsigma_k^{(1)}} \; ,
	\end{equation*}
	where $\displaystyle \max_{k \in \{1,..., 9\}} \left\{ -\sigma_k^{(1)} + \varepsilon \varsigma_k^{(1)} \right\} < -\mu + \varepsilon \mu$ and the decay rate $t^{-\mu + \varepsilon \mu}$ is optimal. The exponents $\sigma_k^{(1)}$, $\varsigma_k^{(1)}$ and the constants $C_k(u_0) \geqslant 0$ are given in the proof.
\end{THM0.2}

\noindent The resulting decay rate is optimal because it is attained on the left boundary of the region (see Theorem \ref{THM10}) and the case $\mu \geqslant \frac{1}{2}$ is also studied in the paper (see Theorem \ref{THM9} for a complete statement). This theorem is an application of our abstract results in a Schrödinger-like case (see Corollary \ref{COR2}) in which we exploit the explicit control of the blow-up of the remainder.

Furthermore the results of this paper permit to obtain propagation features for the wave packets. As for example, we furnish time asymptotic expansions of the solution in certain space-time cones. We provide leading terms exhibiting the optimal decay rates in these regions :

\begin{THM0.3} \label{THM0.3}
	Suppose that $u_0$ satisfies Condition \emph{(C$_{[p_1,p_2],\mu}$)} with $\mu > \frac{1}{2}$, and choose a real number $\varepsilon > 0$ such that $p_1 + \varepsilon < p_2$. For all $(t,x) \in (0,+\infty) \times \R$ satisfying $2 (p_1 + \varepsilon) < \frac{x}{t} < 2 \, p_2$, define $H(t,x,u_0) \in \C$ as follows,
	\begin{align*}
		& \bullet \hspace{5mm} H(t,x,u_0) := \frac{1}{2 \sqrt{\pi}} \, e^{-i \frac{\pi}{4}} \, e^{i \frac{x^2}{4t}} \, \tilde{u} \hspace{-1mm} \left( \frac{x}{2t} \right) \left( \frac{x}{2t} - p_1 \right)^{\mu-1} \; .
	\end{align*}
	Then for all $(t,x) \in (0,+\infty) \times \R$ satisfying $2 (p_1 + \varepsilon) < \frac{x}{t} < 2 \, p_2$, we have
	\begin{equation*}
		\left| u(t,x) - H(t,x,u_0) \, t^{-\frac{1}{2}} \right| \leqslant \sum_{k = 1}^9 \, C_k(u_0, \varepsilon) \, t^{-\sigma_k^{(2)}} \; ,
	\end{equation*}
	where $\displaystyle \max_{k \in \{1,...,9\}} \left\{ -\sigma_k^{(2)} \right\} < - \frac{1}{2}$.  The exponents $\sigma_k^{(2)}$ and the constants $C_k(u_0, \varepsilon) \geqslant 0$ are given in the proof.
\end{THM0.3}

\noindent We furnish also similar results in the case $\mu \leqslant \frac{1}{2}$ (see Theorem \ref{THM6} for a complete statement). Let us remark that the modulus of the coefficient $H(t,x,u_0)$ as well as the remainder tend to infinity when $\frac{x}{2t}$ tends to $p_1$; actually this is a consequence of the limitations of the stationary phase method and it is not specific to the free Schrödinger equation. Thus in order to describe the transition between the cone and the critical direction, the uniform estimate of the solution provided in Theorem \ref{THM0.2} seems to be more appropriate than the expansion calculated in Theorem \ref{THM0.3}.

Nevertheless the time asymptotic expansion in the cone permits to study precisely the energy flow of the solution. In the $L^2$-case, we prove that the energy\footnote[1]{For simplicity, we call in this paper \emph{energy} the $L^2$-norm of the solution on subsets of $\R$, in contrast to the usual quantum mechanics interpretation as probability of occurrence of the particle in a subset.} tends to be concentrated in the cone when the time tends to infinity :

\begin{THM0.4} \label{THM0.4}
	Suppose that the hypotheses of Theorem \ref{THM0.3} are satisfied with $\mu > \frac{1}{2}$ and define the interval $I_t$ as follows,
	\begin{equation*}
		I_t := \Big( \, 2 \, (p_1 + \varepsilon) \,  t, \, 2 \, p_2 \, t \, \Big) \; .
	\end{equation*}
	Then we have
	\begin{equation*}
		\left| \big\| u(t,.) \big\|_{L^2(I_t)} - \frac{1}{\sqrt{2 \pi}} \, \big\| \tf u_0 \big\|_{L^2(p_1+\varepsilon, p_2)} \right| \leqslant \sum_{k=1}^9 \tilde{C}_k(u_0,\varepsilon) \, t^{- \sigma_k^{(2)} + \frac{1}{2}} \; ,
	\end{equation*}
	where the constants $\tilde{C}_k(u_0, \varepsilon) \geqslant 0$ are given in the proof.
\end{THM0.4}

\noindent See Corollary \ref{COR3} for a complete statement. The physical principle of \textit{group velocity} applied to the free Schrödinger equation on the real line says roughly that the energy associated with a frequency $p$ propagates with the speed given by the stationary phase method, in this case $\frac{x}{t} = 2 p$. Theorems \ref{THM0.3} and \ref{THM0.4} furnish a precise meaning of this principle in our case.\\

In Section 1, we first recall Erdélyi's result concerning asymptotic expansions with remainder estimates of oscillatory integrals of the type
\begin{equation*}
	\int_{p_1}^{p_2} U(p) \, e^{i \omega \psi(p)} \, dp \; .
\end{equation*}
The amplitude $U$ can be singular at $p_1$ and $p_2$; factorizing the singularities, $U$ can be written as follows,
\begin{equation} \label{amplitude}
	U(p) = (p-p_1)^{\mu_1-1} \, (p_2 - p)^{\mu_2-1} \, \tilde{u}(p) \; ,
\end{equation}
where $\mu_1, \mu_2 \in (0,1)$ and $\tilde{u}$ is called the regular part of the amplitude.  The phase function $\psi$ is  allowed to have stationary points of real order at the endpoints; more precisely, we suppose the factorization
\begin{equation} \label{phase}
	\psi'(p) = (p-p_1)^{\rho_1-1} \, (p_2 - p)^{\rho_2-1} \, \tilde{\psi}(p) \; ,
\end{equation}
where $\rho_1, \rho_2 \in \R$ are larger than $1$ and $\tilde{\psi}$, which satisfies $\tilde{\psi} > 0$ on $[p_1,p_2]$, is called the non-degenerate part of the phase. We give the proof following the lines of the original demonstration: we start by splitting the integral employing a cut-off function which separates the endpoints of the interval. Then we use explicit substitutions to simplify the phases by exploiting the hypothesis \eqref{phase}. Afterwards  integrations by parts create the expansion of the integral and provide the remainder terms, that we estimate to conclude. The two last steps are carried out using complex analysis in one variable and using the factorization \eqref{amplitude}. Especially the application of Cauchy's theorem allows to shift the integration path of the integrals defining the functions created by the integrations by parts into a region of controllable oscillations of the integrands. This strategy, coupled with the explicit substitutions, leads to a precise estimate of the error.\\
Then we treat the case of the absence of amplitude singularities, which will be essential for certain applications. The previous error estimate furnishes here only the same decay rate for the highest term of the expansion and the remainder. The remedy proposed by Erdélyi \cite[p.55]{erdelyi} leads to complicated formulas when written down and does not seem possible in the case of stationary points of non integer order. To refine this analysis, we work on the integrals defining the functions created by the integrations by parts and involved in the remainder. Introducing a new parameter, we obtain estimates of these integrals permitting a balance between their singular behaviour with respect to the integration variable of the remainder and their decay with respect to $\omega$. Putting these new estimates into the remainder, we obtain the above mentioned refinement.

In Section 2, we modify the above stationary phase method by replacing the smooth cut-off function by a characteristic function. For this purpose, we consider a fixed cutting-point $q \in (p_1,p_2)$ and then we follow the lines of the original result: we carry out explicit substitutions and we integrate by parts. Here the characteristic function produces new terms related to the cutting-point and, due to the different substitutions employed in the two integrals, these terms are not the same. However by expanding them with respect to the parameter $\omega$, we observe  that the first terms cancel out but not the remainders. Finally, we estimate these new remainders related to the cutting-point as well as the remainders related to $p_1$ and $p_2$ to conclude.\\
As previously, we have to refine the estimate of the remainder related to $p_j$ if this point is a regular point of the amplitude. For this purpose, we follow simply the lines of the method described in the preceding section. Moreover, we note that the remainder related to the point $q$ is always negligible as compared with the expansion.

In Section 3, we apply the new method to an oscillatory integral where the phase is a polynomial function of degree 2 and where the amplitude has a unique singularity at the left endpoint of the interval, in view of applications to the free Schrödinger equation. Thanks to the explicitness of the phase function and the preciseness of the previous results, we furnish an asymptotic expansion together with remainder estimates which depend explicitly on the distance between the stationary point and the singularity. This explicitness permits to furnish uniform estimates of the remainder in parameter regions where the stationary point $p_0$ is far from the singularity $p_1$. Moreover we exploit this explicitness to establish uniform estimates of the integral in parameter regions approaching $p_1$ along curves parametrized by $\omega$, enlarging the regions in which we control the decay of the oscillatory integral. We expand the integral on the boundaries of these regions to show the optimality of the decay rates.

In Section 4, we consider the Fourier solution formula of the free Schrödinger equation on the line with initial conditions in a compact frequency band with a singular frequency at one of the endpoints. We are interested in the influence of the frequency band and of the singularity on the dispersion. Thanks to a rewriting of the solution as an oscillatory integral with respect to time, the results of Section 3 are applicable. Firstly we establish an explicit remainder estimate for an asymptotic expansion of the solution formula with respect to time in certain space-time cones related to the frequency band, providing the optimal decay rates in these cones (see Theorem \ref{THM0.3} for the case of weak singularities). The main ideas of the method have been used in \cite[Section 3]{fam}. In the $L^2$-case, we derive an estimate of the $L^2$-norm showing that the energy of waves packets in frequency bands is asymptotically localized in cones generated by the bands and behaves as a laminar flow. Then we exhibit the slow decay $t^{-\frac{\mu}{2}}$, coming from the interaction between the singularity and the stationary point, on a certain space-time direction. Finally we use the last results of the previous section in order to furnish optimal and uniform estimates of the solution in curved regions which are larger than any space-time cone (see Theorem \ref{THM0.2} for the case of strong singularities).

In Section 5, we furnish finally the proofs of intermediate results which have been only sketched in the book of Erdélyi \cite{erdelyi}. They play a key role in Erdélyi's stationary phase method.\\
We use essentially Taylor's theorem to exploit the factorization of the zeros of the derivative of the phase in view of defining a substitution which is a diffeomorphism. At the same time, this technique permits to extract the holomorphic part of the amplitude. Then integrations by parts create integrals of this holomorphic part and their integration path is shifted to a ray in the complex plane. Technically the completeness of the space of holomorphic functions is employed in this context. These aspects do not appear clearly in the original sketch of Erdélyi's proof.\\

This paper has been inspired by \cite{fam} and \cite{fam3}, where the authors consider the Klein-Gordon equation on a star shaped network with a potential which is constant but different on each semi-infinite branch. The authors are interested in the influence of the coefficients of the potential on the time asymptotic behaviour of the solution. To do so, they calculate asymptotic expansions to one term with respect to time of the solution with initial data in frequency bands in \cite{fam} and they exploit these expansions to describe the time asymptotic energy flow of wave packets in \cite{fam3}.\\
They notice that the asymptotic expansion degenerates when the frequency band approaches certain critical values coming from potential steps. These critical values play a similar role as the singularities of the Fourier transform of the initial conditions in the present paper. Hence our refined version of Erdélyi's expansion theorem could help to improve the comprehension of the problem of the blow-up of the expansion.\\
Moreover the paper \cite{fam} shows the way to obtain this type of results for the Schrödinger equation on domains with canonical geometry and canonical potential permitting sufficiently explicit spectral theoretic solution formulas.

Let us mention the paper \cite{dewez} in which $L^{\infty}$-norm estimates of solutions of dispersive equations defined by Fourier multipliers are provided. The initial conditions are supposed to have singular Fourier transforms which are compactly supported or having a sufficient decay at infinity, and the symbols of the Fourier multipliers satisfy a convexity-type hypothesis. In particular the results are applicable to the free Schrödinger equation.\\
To establish these results, the author establishes a generalization of the van der Corput lemma for oscillatory integrals (see for example \cite[p. 332]{stein}) allowing the phase to have a stationary point of non-integer order and the amplitude to have an integrable singularity, furnishing uniform estimates of the integral with respect to the position of the stationary point.\\
Combining \cite{dewez} with the present paper permits to solve the problem of the optimal $L^{\infty}$-time decay rate for the free Schrödinger equation with initial data in frequency bands and having singular frequencies. As perspective for future research, let us mention the same issue for the Klein-Gordon equation on the star shaped network \cite{fam}.

We can point out the usefulness of detailed informations on the motion of wave packets in frequency bands by citing \cite{fam0}. In this paper, the authors study the time-asymptotic behaviour of the solution of the Schrö\-din\-ger equation on star-shaped networks with a localized potential. They establish a perturbation inequality which shows that the evolution of high frequency signals is close to their evolution without potential.

One can also mention the methods and the results of the papers \cite{cazenave1998} and \cite{cazenave2010}, in which the time-decay rate of the free Schrö\-din\-ger equation is considered. In \cite{cazenave1998}, singular initial conditions are constructed to derive the exact $L^p$-time decay rates of the solution, which are slower than the classical results for regular initial conditions. In \cite{cazenave2010}, the authors construct initial conditions in Sobolev spaces (based on the Gaussian function), and they show that the related solutions has no definite $L^p$-time decay rates, nor coefficients, even though upper estimates for the decay rates are established.\\
The papers \cite{cazenave1998} and \cite{cazenave2010} use special formulas for functions and their Fourier transforms, which are themselves based on complex analysis. In our results, we furnish slower decay rates by considering initial conditions with singular Fourier transforms. Here, complex analysis is directly applied to the solution formula of the equation, which permits to obtain results for a whole class of functions. Our method seems therefore to be more flexible.

Further one could review other existing results on the Schrödinger equation with localized potential, e.g. \cite{goldbergschlag} and \cite{weder}, and non linear equations, e.g. \cite{strunk}, to check the possible use of refined estimates of the free equation. An interesting issue could be to find optimal decay conditions of the potential.

Then let us present for comparison two classical versions of the stationary phase method coming from the literature. In \cite[chapter 4]{evans}, the author assumes that the amplitude belongs to $\mathcal{C}_c^{\infty}(\R^d)$ (for $d \geqslant 1$) and supposes that the stationary points of the phase are non-degenerate. Firstly the author employs Morse's lemma to simplify the phase function. Then using Fubini's theorem, he obtains a product of functions, where each of them is a Fourier transform of a tempered distribution. Finally he computes these Fourier transforms and estimates them, leading to the result.
Nevertheless, the use of Morse's lemma implies a loss of precision regarding the estimate of the remainder. Indeed, Morse's lemma is based on the implicit function theorem and so the substitution is not explicit.

In \cite[chapter 7]{hormander}, the author provides a stationary phase method which is different from \cite{evans}. It is assumed that the amplitude $U$ and the phase have a certain regularity on $\R^d$ (for $d \geqslant 1$) and that $U$ has a compact support. Asymptotic expansions of the oscillatory integral are given by using Taylor's formula of the phase, where the stationary point is supposed non-degenerate. Morse's lemma is not needed in this situation. However, stronger hypothesis concerning the phase are required in order to bound uniformly the remainder by a constant, which is not explicit.

Finally we mention some papers exploiting results on oscillatory integrals. In \cite{strichartz}, the author furnishes the first Strichartz-type estimates for the Schrödinger equation and the Klein-Gordon equation. Using complex analysis, the author provides estimates of the $L^2(S)$-norm of Fourier transforms of functions belonging to $L^q(\R^d)$, for some $q \geqslant 1$, where $S$ is a quadratic surface. These considerations lead to the above mentioned estimates.

One-dimensional Schrödinger equations with singular coefficients are studied in \cite{banica}. Dispersion inequalities and Strichartz-type estimates are furnished and we observe that the singularities do not influence the results as compared with the regular case.

In \cite{fam94}, a global (probably not optimal) $L^{\infty}$-time decay estimate for the Klein-Gordon equation on $\R$ with potential steps has been proved using the van der Corput inequality, spectral theory and the methods of \cite{msw}.\\
Note that the idea of considering frequency bands to obtain qualitative informations on the solution is introduced in \cite{fam2}, in the same setting as in \cite{fam94}.

The article \cite{liess1} deals with the time decay rates for the system of crystal optics. The stationary phase method \cite{hormander} is employed to obtain the decay rates. Observe that a change of parameters is carried out to obtain a bounded phase in $\mathcal{C}^4$, which permits to apply this stationary phase method, like in \cite{fam}.

In \cite{liess2}, the authors are interested in the decay of Fourier transforms on singular surfaces. They consider phase functions having stationary points of higher integer order, motivated by Erdélyi's result \cite{erdelyi}.\\

\noindent \textbf{Acknowledgements:}\\
The authors thank S. De Bièvre, R. Haller-Dintelmann and V. Régnier for valuable discussions, and E. Creusé for valuable support. The second named author has been supported by a research grant from the excellence laboratory in mathematics and physics CEMPI and the region Nord-Pas-de-Calais.

\section{Erdélyi's expansion formula: complete proofs and slight improvements}

\hspace{2.5ex} Before formulating the results of this section, let us introduce the assumptions related to the amplitude and to the phase.\\

Let $p_1,p_2$ be two finite real numbers such that $p_1 < p_2$.\\

\noindent \textbf{Assumption ($\mathbf{P_{\rho_1,\rho_2,N}}$).} Let $\rho_1, \rho_2 \geqslant 1$ and $N \in \N \backslash \{0\}$.\\
A function $\psi : [p_1,p_2] \longrightarrow \R$ satifies Assumption (P$_{\rho_1, \rho_2, N}$) if and only if $\psi \in \mathcal{C}^1\big([p_1,p_2], \R\big)$ and there exists a function $\tilde{\psi} : [p_1,p_2] \longrightarrow \R$ such that
	\begin{equation*}
		\forall \, p \in [p_1,p_2] \qquad \psi'(p) = (p-p_1)^{\rho_1 -1} (p_2-p)^{\rho_2 -1} \tilde{\psi}(p) \; ,
	\end{equation*}
	where $\tilde{\psi} \in \mathcal{C}^N\big([p_1,p_2], \R\big)$ is assumed positive.\\ The points $p_j$ $(j=1,2)$ are called \emph{stationary points} of $\psi$ of order $\rho_j -1$, and $\tilde{\psi}$ the \emph{non-degenerate part} of $\psi$.\\
	
\noindent \textbf{Assumption ($\mathbf{A_{\mu_1,\mu_2,N}}$).} Let $\mu_1 , \mu_2 \in (0,1]$ and $N \in \N \backslash \{0\}$.\\
A function $U : (p_1, p_2) \longrightarrow \C$ satisfies Assumption (A$_{\mu_1,\mu_2,N}$) if and only if there exists a function $\tilde{u} : [p_1,p_2] \longrightarrow \C$ such that
	\begin{equation*}
		\forall \, p \in (p_1,p_2) \qquad U(p )= (p - p_1)^{\mu_1 -1} (p_2 - p)^{\mu_2 -1} \tilde{u}(p) \; ,
	\end{equation*}
	where $\tilde{u} \in \mathcal{C}^N\big([p_1,p_2], \C\big)$, and $\tilde{u}(p_j) \neq 0$ if $\mu_j \neq 1$ $(j=1,2)$.\\
	The points $p_j$ are called \emph{singularities} of $U$, and $\tilde{u}$ the \emph{regular part} of $U$.

\begin{REM1}
	\em The hypothesis $\tilde{u}(p_j) \neq 0$ if $\mu_j \neq 1$ prevents the function $\tilde{u}$ from affecting the behaviour of the singularity $p_j$.
\end{REM1}

\subsection*{Non-vanishing singularities: Erdélyi's theorem}

\hspace{2.5ex} The aim of this subsection is to state Erdélyi's result \cite[section 2.8]{erdelyi} and to provide a complete proof.\\

Let us define some objects that will be used throughout this section.

\begin{DEF1} \label{DEF1}
	\begin{enumerate}
	\item Let $\eta \in \big(0, \frac{p_2 - p_1}{2}\big)$. For $j=1,2$, let $\varphi_j : I_j \longrightarrow \R$ be the functions defined by
	\begin{equation*}
		\varphi_1(p) := \big( \psi(p) - \psi(p_1) \big)^{\frac{1}{\rho_1}} \qquad , \qquad \varphi_2(p) := \big( \psi(p_2) - \psi(p) \big)^{\frac{1}{\rho_2}} \; ,
	\end{equation*}
	with $I_1 := [p_1, p_2 - \eta]$, $I_2 := [p_1 + \eta, p_2]$ and $s_1 := \varphi_1(p_2 - \eta)$, $s_2 := \varphi_2(p_1 + \eta)$.
	\item For $j=1,2$, let $k_j : (0,s_j] \longrightarrow \C$ be the functions defined by
	\begin{equation*}
			k_j(s) := U\big( \varphi_j^{-1}(s) \big) \, s^{1- \mu_j} \, \big( \varphi_j^{-1} \big)'(s) \; ,
	\end{equation*}
	which can be extended to $[0,s_j]$ (see Proposition \ref{LEM2}).
	\item Let $\nu : [p_1,p_2] \longrightarrow \R$ be a smooth function such that
	\begin{equation*}
		\left\{ \begin{array}{rl}
					& \nu = 1 \qquad \text{on} \quad [p_1, p_1 + \eta] \; , \\
					& \nu = 0 \qquad \text{on} \quad [p_2 - \eta, p_2] \; ,\\
					& 0 \leqslant \nu \leqslant 1 \; ,
				\end{array} \right.
	\end{equation*}
	where $\eta$ is defined above.\\ For $j=1,2$, let $\nu_j : [0,s_j] \longrightarrow \R$ be the functions defined by
	\begin{equation*}
		\nu_1(s) := \nu \circ \varphi_1^{-1}(s) \qquad , \qquad \nu_2(s) := (1 - \nu) \circ \varphi_2^{-1}(s) \; .
	\end{equation*}
	\item For $s > 0$, let $\Lambda^{(j)}(s)$ be the complex curve defined by
	\begin{equation*}
		\Lambda^{(j)}(s) := \left\{ s + t e^{(-1)^{j+1} i \frac{\pi}{2\rho_j}} \, \Big| \, t \geqslant 0 \right\} \; .
	\end{equation*}
	\end{enumerate}
\end{DEF1}

\begin{THM1} \label{THM1}
	Let $N \in \N \backslash \{0\}$, $\rho_1, \rho_2 \geqslant 1$ and $\mu_1, \mu_2 \in (0,1)$. Suppose that the functions $\psi : [p_1,p_2] \longrightarrow \R$ and $U : (p_1, p_2) \longrightarrow \C$ satisfy Assumption \emph{(P$_{\rho_1,\rho_2,N}$)} and Assumption \emph{(A$_{\mu_1,\mu_2,N}$)}, respectively.\\ Then for $j=1,2$ , there are functions $A_N^{(j)}$, $R_N^{(j)} : (0,+\infty) \longrightarrow \C$ such that
	\begin{equation*}
		\left\{ \begin{array}{rl}
			& \displaystyle \int_{p_1}^{p_2} U(p) \, e^{i \omega \psi(p)} \, dp = \sum_{j=1,2} \left( A_N^{(j)}(\omega) + R_N^{(j)}(\omega) \right) \; , \\
			& \displaystyle \left| R_N^{(j)}(\omega) \right| \leqslant \frac{1}{(N-1)!} \, \frac{1}{\rho_j} \, \Gamma\bigg(\frac{N}{\rho_j}\bigg) \int_0^{s_j} s^{\mu_j-1} \left| \frac{d^N}{ds^N}\big[\nu_j k_j\big](s) \right| \, ds \; \omega^{-\frac{N}{\rho_j}} \; ,
		\end{array} \right.
	\end{equation*}
	for all $\omega > 0$, where
	\begin{equation*}
		\begin{aligned}
			\bullet \quad & A_N^{(j)}(\omega) := e^{i \omega \psi(p_j)} \, \sum_{n=0}^{N-1} \Theta_{n+1}^{(j)}(\rho_j,\mu_j) \, \frac{d^n}{ds^n}\big[ k_j \big](0) \, \omega^{-\frac{n+\mu_j}{\rho_j}} \; , \\
			\bullet \quad & R_N^{(j)}(\omega) := (-1)^{N+1+j} \, e^{i\omega \psi(p_j)} \int_0^{s_j} \phi_{N}^{(j)}(s,\omega, \rho_j, \mu_j) \, \frac{d^N}{ds^N}\big[ \nu_j k_j \big] (s) \, ds \; ,
		\end{aligned}
	\end{equation*}
	and for $n=0,...,N-1$,
	\begin{equation*}
		\begin{aligned}
			\bullet \quad & \Theta_{n+1}^{(j)}(\rho_j, \mu_j) := \frac{(-1)^{j+1}}{n! \, \rho_j} \, \Gamma \Big(\frac{n+\mu_j}{\rho_j}\Big) \, e^{(-1)^{j+1} i \frac{\pi}{2} \, \frac{n+\mu_j}{\rho_j}} \; , \\
			\bullet \quad & \phi_{n+1}^{(j)}(s,\omega,\rho_j, \mu_j) := \frac{(-1)^{n+1}}{n!} \int_{\Lambda^{(j)}(s)} (z-s)^n z^{\mu_j -1} e^{(-1)^{j+1} i\omega z^{\rho_j}} \, dz \; .
		\end{aligned}
	\end{equation*}
\end{THM1}

\begin{proof}
	For fixed $\omega > 0$, $\rho_j \geqslant 1$ and $\mu_j \in (0,1)$, we shall note $\phi_{n}^{(j)}(s,\omega)$ instead of $\phi_{n}^{(j)}(s,\omega,\rho_j,\mu_j)$. Now let us divide the proof in five steps.\\
	
	\noindent \textit{First step: Splitting of the integral.} Using the cut-off function $\nu$, we can write the integral as follows,
	\begin{equation*}
		\int_{p_1}^{p_2} U(p) \, e^{i \omega \psi(p)} \, dp = \tilde{I}^{(1)}(\omega) + \tilde{I}^{(2)}(\omega) \; ,
	\end{equation*}
	where
	\begin{equation*}
		\tilde{I}^{(1)}(\omega) := \int_{p_1}^{p_2-\eta} \nu(p) \, U(p) \, e^{i \omega \psi(p)} \, dp \quad , \quad \tilde{I}^{(2)}(\omega) := \int_{p_1+\eta}^{p_2} \big( 1-\nu(p) \big) \, U(p) \, e^{i \omega \psi(p)} \, dp \; .
	\end{equation*}
	\textit{Second step: Substitution.} Proposition \ref{LEM1} affirms that $\varphi_j : I_j \longrightarrow [0,s_j]$ is a $\mathcal{C}^{N+1}$-diffeomorphism with $q_1 = p_2 - \eta$ here. Using the substitution $s = \varphi_1(p)$, we obtain
	\begin{equation*}
		\begin{aligned}
			\tilde{I}^{(1)}(\omega)	& = \int_{p_1}^{p_2-\eta} \nu(p) \, U(p) \, e^{i \omega \psi(p)} \, dp \\
								& = e^{i \omega \psi(p_1)} \int_0^{s_1} \nu \big( \varphi_1^{-1}(s) \big) \, U \big( \varphi_1^{-1}(s) \big) \, e^{i \omega s^{\rho_1}} (\varphi_1^{-1})'(s) \, ds \\
								& = e^{i \omega \psi(p_1)} \int_0^{s_1} \nu \big( \varphi_1^{-1}(s) \big) \, U \big( \varphi_1^{-1}(s) \big) s^{1-\mu_1} (\varphi_1^{-1})'(s) \, s^{\mu_1-1}  e^{i \omega s^{\rho_1}} \, ds \\
								& = e^{i \omega \psi(p_1)} \int_0^{s_1} \nu_1(s) \, k_1(s) \, s^{\mu_1-1} e^{i \omega s^{\rho_1}} \, ds \; ,
		\end{aligned}
	\end{equation*}
	where $k_1$ and $\nu_1$ are introduced in Definition \ref{DEF1}. In a similar way, we obtain
	\begin{equation*}
		\tilde{I}^{(2)}(\omega) = - e^{i \omega \psi(p_2)} \int_0^{s_2} \nu_2(s) \, k_2(s) \, s^{\mu_2-1} e^{-i \omega s^{\rho_2}} \, ds \; .
	\end{equation*}
	 Note that the minus sign comes from the decrease of $\varphi_2$.\\

	\noindent \textit{Third step: Integrations by parts.} Corollary \ref{COR1} provides successive primitives of the function $\displaystyle s \longmapsto s^{\mu_j - 1} e^{(-1)^{j+1} i\omega s^{\rho_j}}$. Moreover Proposition \ref{LEM2} ensures that $k_j \in \mathcal{C}^N\big([0,s_j], \C \big)$. Thus by $N$ integrations by parts, we obtain
	\begin{equation*}
		\begin{aligned}
			e^{-i \omega \psi(p_1)} \, \tilde{I}^{(1)}(\omega)	& = \int_0^{s_1} \nu_1(s) \, k_1(s) \, s^{\mu_1-1} e^{i \omega s^{\rho_1}} \, ds \\
														& = \Big[ \phi_1^{(1)}(s,\omega) \big( \nu_1 k_1 \big)(s)\Big]_0^{s_1} - \int_0^{s_1} \phi_1^{(1)}(s,\omega) \, \frac{d}{ds}\big[ \nu_1 k_1 \big](s) \, ds \\
														& = \dots \\
														& = \sum_{n=0}^{N-1} (-1)^n \left[ \phi_{n+1}^{(1)}(s,\omega) \, \frac{d^n}{ds^n}\big[ \nu_1 k_1 \big](s)\right]_0^{s_1} \\
														& \qquad \qquad + (-1)^N \int_0^{s_1} \phi_{N}^{(1)}(s,\omega) \, \frac{d^N}{ds^N}\big[ \nu_1 k_1 \big](s) \, ds \; .
		\end{aligned}
	\end{equation*}
	One can simplify the last expression by using the properties of the function $\nu_1$: by hypothesis, $\nu(p_1)=1$, $\nu(p_2-\eta)=0$ and $\displaystyle \frac{d^n}{dp^n}[\nu](p_1)=\frac{d^n}{dp^n}[\nu](p_2-\eta)=0$, for $n \geqslant 1$. So the definition of $\nu_1$ implies
	\begin{equation*}
		\nu_1(0) = \nu(p_1) = 1 \qquad , \qquad \nu_1(s_1) = \nu(p_2-\eta) = 0 \; ;
	\end{equation*}
	and by the product rule applied to $\nu_1 k_1$, it follows
	\begin{equation*}
		\frac{d^n}{ds^n}\big[ \nu_1 k_1 \big](0) = \frac{d^n}{ds^n}[k_1](0) \qquad , \qquad \frac{d^n}{ds^n}\big[ \nu_1 k_1 \big](s_1) = 0 \; .
	\end{equation*}
	This leads to
	\begin{equation*}
		\begin{aligned}
			\tilde{I}^{(1)}(\omega)	& = \sum_{n=0}^{N-1} (-1)^{n+1} \phi_{n+1}^{(1)}(0,\omega) \, \frac{d^n}{ds^n} \big[ k_1 \big](0) \, e^{i \omega \psi(p_1)}\\
							& \qquad \qquad + (-1)^{N} \, e^{i \omega \psi(p_1)} \int_0^{s_1} \phi_{N}^{(1)}(s,\omega) \, \frac{d^N}{ds^N}\big[ \nu_1 k_1 \big](s) \, ds \; .
		\end{aligned}
	\end{equation*}
	By similar computations, we obtain
	\begin{equation*}
		\begin{aligned}
			\tilde{I}^{(2)}(\omega)	& = \sum_{n=0}^{N-1} (-1)^n \phi_{n+1}^{(2)}(0,\omega) \, \frac{d^n}{ds^n} \big[ k_2 \big] (0) \, e^{i \omega \psi(p_2)}\\
							& \qquad \qquad + (-1)^{N+1} \, e^{i \omega \psi(p_2)} \int_0^{s_2} \phi_{N}^{(2)}(s,\omega) \, \frac{d^N}{ds^N}\big[ \nu_2 k_2 \big](s) \, ds \; . \\
		\end{aligned}
	\end{equation*}
	
	\vspace{0.4cm}
	
	\noindent \textit{Fourth step: Calculation of the main terms.} Let us compute the coefficient $\phi_{n+1}^{(j)}(0,\omega)$. We recall the expression of $\phi_{n+1}^{(j)}(s,\omega)$ from Corollary \ref{COR1},
	\begin{equation*}
		\phi_{n+1}^{(j)}(s,\omega) = \frac{(-1)^{n+1}}{n!} \int_{\Lambda^{(j)}(s)} (z-s)^n \, z^{\mu_j - 1} \, e^{(-1)^{j+1} i \omega z^{\rho_j}} \, dz \; ,
	\end{equation*}
	for all $s \in [0,s_j]$ and $n=0,...,N-1$. Choosing $j=1$, putting $s=0$ and parametrizing the curve $\Lambda^{(1)}(0)$ with $\displaystyle z = t \, e^{i \frac{\pi}{2 \rho_1}}$ lead to
	\begin{equation*}
		\phi_{n+1}^{(1)}(0,\omega) = \frac{(-1)^{n+1}}{n!} \, e^{i \frac{\pi}{2} \, \frac{n+\mu_1}{\rho_1}} \int_0^{+\infty} t^{n+\mu_1 - 1} \, e^{-\omega t^{\rho_1}} \, dt \; .
	\end{equation*}
	Setting $y= \omega \, t^{\rho_1}$ in the previous integral gives
	\begin{equation*}
		\begin{aligned}
			\phi_{n+1}^{(1)}(0,\omega)	& = \frac{(-1)^{n+1}}{n!} \, e^{i \frac{\pi}{2} \, \frac{n+\mu_1}{\rho_1}} \, (\rho_1 \, \omega)^{-1} \int_0^{+\infty} \left( \frac{y}{\omega} \right)^{\frac{n+\mu_1}{\rho_1} - 1} \, e^{-y} \, dy \\
										& = \frac{(-1)^{n+1}}{n!} \, e^{i \frac{\pi}{2} \, \frac{n+\mu_1}{\rho_1}} \, \frac{1}{\rho_1} \, \Gamma \hspace{-1mm} \left(\frac{n + \mu_1}{\rho_1} \right) \, \omega^{-\frac{n+\mu_1}{\rho_1}} \; ,
		\end{aligned}
	\end{equation*}
	where $\Gamma$ is the Gamma function defined by
	\begin{equation*}
		\Gamma \, : \, z \in \big\{ z \in \C \, \big| \, \Re(z) > 0 \big\} \, \longmapsto \, \int_0^{+\infty} t^{z-1} e^{-t} \, dt \in \C \; .
	\end{equation*}
	A similar work provides
	\begin{equation*}
		\phi_{n+1}^{(2)}(0,\omega) = \frac{(-1)^{n+1}}{n!} \, e^{-i \frac{\pi}{2} \, \frac{n+\mu_2}{\rho_2}} \, \frac{1}{\rho_2} \, \Gamma \hspace{-1mm} \left(\frac{n + \mu_2}{\rho_2} \right) \, \omega^{-\frac{n+\mu_2}{\rho_2}} \; .
	\end{equation*}
	Then we obtain
	\begin{equation*}
		e^{i \omega \psi(p_j)} \sum_{n=0}^{N-1} (-1)^{n+j} \, \phi_{n+1}^{(j)}(0,\omega) \, \frac{d^n}{ds^n}\big[ k_j \big](0) = e^{i \omega \psi(p_j)} \sum_{n=0}^{N-1} \Theta_{n+1}^{(j)}(\rho_j,\mu_j) \, \frac{d^n}{ds^n} \big[ k_j \big] (0) \, \omega^{- \frac{n+\mu_j}{\rho_j} } \; ,
	\end{equation*}
	where $\displaystyle \Theta_{n+1}^{(j)}(\rho_j,\mu_j) := \frac{(-1)^{j+1}}{n! \, \rho_j} \, \Gamma \hspace{-1mm} \left( \frac{n + \mu_j}{\rho_j} \right) \, e^{(-1)^{j+1} \frac{\pi}{2} \, \frac{n + \mu_j}{\rho_j}}$.
	
	\vspace{0.4cm}
	
	\noindent \textit{Fifth step: Remainder estimates.} The last step consists in estimating the remainders $R_N^{(j)}(\omega)$. For $j=1,2$, we have for all $s \in (0,s_j]$ and for all $t \geqslant 0$,
	\begin{equation}\label{est3}
		s \leqslant \left| s + t e^{(-1)^{j+1} i \frac{\pi}{2 \rho_j}} \right| \qquad \Longrightarrow \qquad s^{\mu_j-1} \geqslant \left| s + t e^{(-1)^{j+1} i \frac{\pi}{2 \rho_j}} \right|^{\mu_j-1} \; ,
	\end{equation}
	since $\mu_j \in (0,1)$. Now parametrize the integral defining $\phi_{N}^{(j)}(s,\omega)$ with $\displaystyle z = s + t e^{(-1)^{j+1} i \frac{\pi}{2 \rho_j}}$ and employ the previous inequality \eqref{est3} to obtain
	\begin{align}
		\Big| \phi_{N}^{(j)}(s,\omega) \Big|	& \leqslant \frac{1}{(N-1)!} \int_0^{+\infty} t^{N-1} \left| \, s + t e^{(-1)^{j+1} i\frac{\pi}{2 \rho_j}} \right|^{\mu_j -1} \left| \, e^{(-1)^{j+1} i \omega \big( s + t e^{(-1)^{j+1} i\frac{\pi}{2 \rho_j}} \big)^{\rho_j}} \right| \, dt \nonumber \\
													& \leqslant \frac{1}{(N-1)!} \, s^{\mu_j-1} \int_0^{+\infty} t^{N-1} e^{-\omega t^{\rho_j}} \, dt \nonumber \\
													& = \label{est4} \frac{1}{(N-1)!} \, s^{\mu_j-1} \, \frac{1}{\rho_j} \, \Gamma \hspace{-1mm} \left(\frac{N}{\rho_j}\right) \omega^{-\frac{N}{\rho_j}} \; ,
	\end{align}
	where the last equality has been obtained by using the substitution $y = \omega \, t^{\rho_j}$. Employing the definition of $R_N^{(j)}(\omega)$ and inequality \eqref{est4} leads to
	\begin{align*}
		\Big| R_{N}^{(j)}(\omega) \Big|	& \leqslant \int_0^{s_j} \left| \, \phi_{N}^{(j)}(s,\omega) \right| \, \left| \frac{d^N}{ds^N} \big[ \nu_j k_j \big](s) \right| \, ds \\
												& \leqslant \frac{1}{(N-1)!} \, \frac{1}{\rho_j} \, \Gamma \hspace{-1mm} \left(\frac{N}{\rho_j} \right) \int_0^{s_j} s^{\mu_j-1} \left| \frac{d^N}{ds^N} \big[\nu_j k_j \big](s) \right| \, ds \; \omega^{-\frac{N}{\rho_j}} \; .
	\end{align*}
	We note that the last integral is well-defined because $\displaystyle \frac{d^N}{ds^N}\big[\nu_j k_j \big] : [0,s_j] \longrightarrow \R$ is continuous and $s \longmapsto s^{\mu_j - 1}$ is locally integrable on $[0,s_j]$.\\
	Finally, we remark that the highest term of the expansion $A_N^{(j)}(\omega)$ behaves like $\omega^{-\frac{N-1+\mu_j}{\rho_j}}$ when $\omega$ tends to infinity. Moreover $R_N^{(j)}(\omega)$ is estimated by $\omega^{-\frac{N}{\rho_j}}$. This implies that the decay rate of the remainder with respect to $\omega$ is higher than the one of the highest term of the expansion. This ends the proof.
\end{proof}

\subsection*{Amplitudes without singularities: refinement of the error estimate}

\hspace{2.5ex} The preceding theorem remains true if we suppose $\mu_j = 1$, that is to say if the amplitude $U$ is regular at the point $p_j$. But in this case, we observe that the decay rates of the  highest term of the expansion related to $p_j$ and of the remainder related to $p_j$ are the same, namely $\omega^{-\frac{N}{\rho_j}}$. Hence the aim of this subsection is to refine the estimate of the remainder in this specific case.\\

For this purpose, we establish the two following lemmas. In the first one, we provide two estimates of the function $s \longmapsto \phi_{N}^{(j)}(s,\omega,\rho_j,1)$: the first estimate is uniform with respect to $s$ but the decay with respect to $\omega$ is not sufficiently fast; on the other hand, the second one provides a better decay with respect to $\omega$ but is singular with respect to $s$. The first estimate is actually established in the proof of Theorem \ref{THM1} and we integrate by parts to establish the second one.

\begin{LEM6} \label{LEM6}
	Let $j \in \{1,2\}$, $\rho_j \geqslant 1$ and $N \in \N \backslash \{0\}$. Then for all $s, \omega > 0$, we have
	\begin{equation*}
		\left\{ \begin{array}{rl}
			& \displaystyle \left| \phi_{N}^{(j)}(s,\omega,\rho_j,1) \right| \leqslant a_{N,\rho_j} \, \omega^{-\frac{N}{\rho_j}} \; , \\
			& \displaystyle \left| \phi_{N}^{(j)}(s,\omega,\rho_j,1) \right| \leqslant b_{N,\rho_j} \, \omega^{-\left(1+\frac{N-1}{\rho_j}\right)} \, s^{1-\rho_j} \, + c_{N,\rho_j} \, \omega^{-\left(1+\frac{N}{\rho}\right)} \, s^{-\rho_j} \; ,
		\end{array} \right.
	\end{equation*}
	where the constants $a_{N,\rho_j}, b_{N,\rho_j}, c_{N,\rho_j} > 0$ are given in the proof.	
\end{LEM6}

\begin{REM4}
	\em Note that we can extend $\phi_{N}^{(j)}(.,\omega,\rho_j,1) \, : \, [0,s_j] \longrightarrow \R$ to $[0,+\infty)$, see Remark \ref{REM2}.
\end{REM4}

\begin{proof}[Proof of Lemma \ref{LEM6}]
	Let us fix $s > 0$, $\omega > 0$ and for example $j=1$. We recall the expression of $\phi_{N}^{(1)}(s,\omega,\rho_1,1)$ with the parametrization of the path $\Lambda^{(1)}(s)$ given in Definition \ref{DEF1}:
	\begin{equation*}
		\phi_{N}^{(1)}(s,\omega,\rho_1,1) = \frac{(-1)^N}{(N-1)!} \int_0^{+\infty} t^{N-1} \, e^{i \frac{\pi (N-1)}{2 \rho_1}} \, e^{i \omega \left(s+te^{i \frac{\pi}{2 \rho_1}}\right)^{\rho_1}} \, dt \, e^{i \frac{\pi}{2 \rho_1}} \; .
	\end{equation*}
	On the one hand, estimate \eqref{est4} is still valid for $\mu_1 = 1$, namely,
	\begin{equation*}
		\left| \phi_{N}^{(1)}(s,\omega,\rho_1,1) \right| \leqslant \frac{1}{(N-1)!} \, \frac{1}{\rho_1} \, \Gamma\bigg(\frac{N}{\rho_1}\bigg) \, \omega^{-\frac{N}{\rho_1}} =: a_{N,\rho_1} \, \omega^{-\frac{N}{\rho_1}} \; ,
	\end{equation*}
	furnishing the first estimate of the lemma.\\
	On the other hand, we establish the second inequality by using integrations by parts. To do so, we remark that for all $s>0$ the first derivative of the function $t \in (0,+\infty) \longmapsto i \omega \big(s+te^{i \frac{\pi}{2 \rho_1}}\big)^{\rho_1}$ does not vanish; therefore we can write
	\begin{equation*}
		e^{i \omega \left(s+te^{i \frac{\pi}{2 \rho_1}}\right)^{\rho_1}} = (i \omega \rho_1)^{-1} \, e^{-i \frac{\pi}{2 \rho_1}} \left(s+te^{i \frac{\pi}{2 \rho_1}}\right)^{1-\rho_1} \, \frac{d}{dt} \hspace{-1mm} \left[e^{i \omega \left(s+te^{i \frac{\pi}{2 \rho_1}}\right)^{\rho_1}}\right] \; .
	\end{equation*}
	Moreover Lemma \ref{LEM3} implies
	\begin{equation} \label{estinf}
		\forall \, s > 0 \qquad \left| e^{i \omega \left(s+te^{i \frac{\pi}{2 \rho_1}}\right)^{\rho_1}} \right| \leqslant e^{- \omega t^{\rho_1}} \longrightarrow  0 \quad , \quad t \longrightarrow +\infty \; .
	\end{equation}
	Now we distinguish the two following cases:
	\begin{itemize}
		\item \textit{Case $N=1$.} Thanks to the two previous observations, we can integrate by parts, providing
		\begin{equation*}
			\begin{aligned}
				\phi_1^{(1)}(s,\omega,\rho_1,1)	& = - (i \omega \rho_1)^{-1} \int_0^{+\infty} \left(s+te^{i \frac{\pi}{2 \rho_1}}\right)^{1-\rho_1} \frac{d}{dt} \hspace{-1mm} \left[e^{i \omega \left(s+te^{i \frac{\pi}{2 \rho_1}}\right)^{\rho_1}}\right] \, dt \\
													& = (i \omega \rho_1)^{-1} s^{1-\rho_1} e^{i \omega s^{\rho_1}} \\
													& \qquad + \frac{1-\rho_1}{i \omega \rho_1} \, e^{i \frac{\pi}{2 \rho_1}} \int_0^{+\infty} \left(s+te^{i \frac{\pi}{2 \rho_1}}\right)^{-\rho_1} e^{i \omega \left(s+te^{i \frac{\pi}{2 \rho_1}}\right)^{\rho_1}} \, dt \; ,
			\end{aligned}
		\end{equation*}
		where the boundary term at infinity is equal to zero according to \eqref{estinf}. It follows
		\begin{align}
			\Big| \phi_1^{(1)}(s,\omega,\rho_1,1) \Big|	& \leqslant (\omega \rho_1)^{-1} \, s^{1-\rho_1} \nonumber \\
															& \qquad + \frac{\rho_1-1}{\omega \rho_1} \int_0^{+\infty} \left| s+te^{i \frac{\pi}{2 \rho_1}}\right|^{-\rho_1} \left|e^{i \omega \left(s+te^{i \frac{\pi}{2 \rho_1}}\right)^{\rho_1}}\right| \, ds \nonumber \\
															& \label{ineq11}\leqslant (\omega \rho_1)^{-1} \, s^{1-\rho_1} + \frac{\rho_1-1}{\omega \rho_1} \, s^{-\rho_1} \int_0^{+\infty} e^{-\omega t^{\rho_1}} \, dt \\
															& = \label{eq1} \frac{1}{\rho_1} \, \omega^{-1} \, s^{1-\rho_1} + \frac{\rho_1 -1}{\rho_1^2} \, \Gamma\bigg(\frac{1}{\rho_1}\bigg) \, \omega^{-\left(1+\frac{1}{\rho_1}\right)} \, s^{-\rho_1} \; \\
															& =: b_{1,\rho_1} \, \omega^{-1} \, s^{1-\rho_1} + c_{1,\rho_1} \, \omega^{-\left(1+\frac{1}{\rho_1}\right)} \, s^{-\rho_1} \; . \nonumber
		\end{align}
		Lemma \ref{LEM3} permits to obtain \eqref{ineq11} by estimating the exponential and we use the substitution $\displaystyle y= \omega \, t^{\rho_1}$ to get \eqref{eq1}.
		\item \textit{Case $N \geqslant 2$.} We proceed as above by using an integration by parts :
		\begin{align}
			\phi_{N}^{(1)}(s,\omega,\rho_1,1)	& = \frac{(-1)^N}{(N-1)!} \, e^{i \frac{\pi (N-1)}{2 \rho_1}} \, (i \omega \rho_1)^{-1} \nonumber \\
												& \qquad \times \int_0^{+\infty} t^{N-1} \left(s+te^{i \frac{\pi}{2 \rho_1}}\right)^{1-\rho_1} \frac{d}{dt} \hspace{-1mm} \left[ e^{i \omega \left(s+te^{i \frac{\pi}{2 \rho_1}}\right)^{\rho_1}} \right] \, dt \nonumber \\
												& = \frac{(-1)^{N+1}}{(N-1)!} \, e^{i \frac{\pi (N-1)}{2 \rho_1}} \, (i \omega \rho_1)^{-1} \nonumber \\
												& \label{ipp} \qquad \times \int_0^{+\infty} \frac{d}{dt} \hspace{-1mm} \left[ t^{N-1} \left(s+te^{i \frac{\pi}{2 \rho_1}}\right)^{1-\rho_1} \right] e^{i \omega \left(s+te^{i \frac{\pi}{2 \rho_1}}\right)^{\rho_1}} \, dt \\
												& = \frac{(-1)^{N+1}}{(N-1)!} \, e^{i \frac{\pi (N-1)}{2 \rho_1}} \, (i \omega \rho_1)^{-1} \nonumber \\
												& \qquad \times \bigg( (N-1) \int_0^{+\infty} t^{N-2} \left(s+te^{i \frac{\pi}{2 \rho_1}}\right)^{1-\rho_1} e^{i \omega \left(s+te^{i \frac{\pi}{2 \rho_1}}\right)^{\rho_1}} \, dt \nonumber \\												& \qquad \qquad + (1-\rho_1) \, e^{i \frac{\pi}{2 \rho_1}} \int_0^{+\infty} t^{N-1} \left(s+te^{i \frac{\pi}{2 \rho_1}}\right)^{-\rho_1} e^{i \omega \left(s+te^{i \frac{\pi}{2 \rho_1}}\right)^{\rho_1}} \, dt \bigg) \, . \nonumber
		\end{align}
		The boundary terms in \eqref{ipp} are zero; indeed we can remark that the term at $0$ vanishes and we can use \eqref{estinf} once again to show that the term at infinity is equal to $0$. Then by similar arguments to those of the preceding case, we obtain
		\begin{align}
			\Big| \phi_{N}^{(1)}(s,\omega,\rho_1,1) \Big|	& \leqslant \frac{\omega^{-1}}{\rho_1 (N-2)!} \int_0^{+\infty} t^{N-2} \left| \left(s+te^{i \frac{\pi}{2 \rho_1}}\right)^{1-\rho_1} e^{i \omega \left(s+te^{i \frac{\pi}{2 \rho_1}}\right)^{\rho_1}} \right| \, dt \nonumber \\
															& \qquad + \frac{(\rho_1-1) \omega^{-1}}{\rho_1 (N-1)!} \int_0^{+\infty} t^{N-1} \left| \left(s+te^{i \frac{\pi}{2 \rho_1}}\right)^{-\rho_1} e^{i \omega \left(s+te^{i \frac{\pi}{2 \rho_1}}\right)^{\rho_1}} \right| \, dt \nonumber \\
															& \leqslant \frac{1}{\rho_1 (N-2)!} \, \omega^{-1} \, s^{1 - \rho_1} \int_0^{+\infty} t^{N-2} \, e^{- \omega t^{\rho_1}} \, dt \nonumber \\
															& \qquad + \frac{(\rho_1-1)}{\rho_1 (N-1)!} \, \omega^{-1} \, s^{-\rho_1} \int_0^{+\infty} t^{N-1} \, e^{- \omega t^{\rho_1}} \, dt \nonumber \\
															& = \frac{1}{\rho_1^2 (N-2)!} \, \Gamma\bigg( \frac{N-1}{\rho_1} \bigg) \, \omega^{-\left(1+ \frac{N-1}{\rho_1}\right)} \, s^{1-\rho_1} \nonumber \\
															& \qquad + \frac{(\rho_1 -1)}{\rho_1^2 (N-1)!} \, \Gamma \bigg( \frac{N}{\rho_1} \bigg) \, \omega^{-\left(1+\frac{N}{\rho_1}\right)} \, s^{-\rho_1} \nonumber \\
															& =: b_{N,\rho_1} \, \omega^{-\left(1+ \frac{N-1}{\rho_1}\right)} \, s^{1-\rho_1} + c_{N,\rho_1} \,  \omega^{-\left(1+\frac{N}{\rho_1}\right)} \, s^{-\rho_1} \; , \nonumber
		\end{align}
		concluding this point.
	\end{itemize}
	A very similar work for $j=2$ provides the conclusion; you replace $\rho_1$ by $\rho_2$ in the expressions of the constants $a_{N,\rho_1}$, $b_{N,\rho_1}$ and $c_{N,\rho_1}$ to obtain $a_{N,\rho_2}$, $b_{N,\rho_2}$ and $c_{N,\rho_2}$.
\end{proof}

Given a function satisfying a system of inequalities similar to the one given in Lemma \ref{LEM6}, we furnish a new estimate which has a sufficiently fast decay with respect to $\omega$ and which is integrable with respect to $s$. We exploit the balance between blow-up and the decay to establish this result. Note that a technical argument requires $\rho \geqslant 2$.

\begin{LEM7}\label{LEM7}
	Let $N \in \N \backslash \{0\}$, $\rho \geqslant 2$ and $f : (0,+\infty) \times (0,+\infty) \longrightarrow \R$ be a function which satisfies the following inequalities:
	\begin{equation*}
		\forall \, s, \omega > 0 \qquad \left\{ \begin{array}{rl} \hspace{-0.3cm}
			& \displaystyle \big| f(s,\omega) \big| \leqslant a \, \omega^{-\frac{N}{\rho}} \; , \\[0.2cm]
			& \displaystyle \big| f(s,\omega) \big| \leqslant b \, \omega^{-\left(1+\frac{N-1}{\rho}\right)} \, s^{1-\rho} \, + c \, \omega^{-\left(1+\frac{N}{\rho}\right)} \, s^{-\rho} \; , \\
		\end{array} \right.
	\end{equation*}
	where $a,b,c > 0$ are three constants.\\
	Fix $\gamma \in (0,1)$ and define $\displaystyle \delta := \frac{\gamma+N}{\rho} \in \left( \frac{N}{\rho},\frac{1+N}{\rho} \right)$. Then we have
	\begin{equation*}
		\forall \, s,\omega > 0  \qquad \big| f(s,\omega) \big| \leqslant L_{\gamma,\rho} \, s^{-\gamma} \, \omega^{-\delta} \; ,
	\end{equation*}
	where $L_{\gamma,\rho} := a \, K_{\rho}^{\; \gamma} > 0$, with $K_{\rho}$ the unique positive solution of
	\begin{equation*}
		a K^{\rho} - b K - c = 0 \; .
	\end{equation*}	
\end{LEM7}

\begin{proof}
	Let $g_1 , g_2 \, : \, (0,+\infty) \, \times \, (0,+\infty) \, \longrightarrow \R $ defined by
	\begin{equation*}
		g_1(s,\omega) := a \, \omega^{-\frac{N}{\rho}} \qquad , \qquad g_2(s,\omega):=  b \, \omega^{-\left(1+\frac{N-1}{\rho}\right)} \, s^{1-\rho} \, + c \, \omega^{-\left(1+\frac{N}{\rho}\right)} \, s^{-\rho} \; .
	\end{equation*}
	Now we fix $\omega > 0$ and we define the function $h_{\omega} : (0,+\infty) \longrightarrow \R$ by
	\begin{equation*}
		h_{\omega}(s) := s^{\rho} \, \big(g_1(s,\omega) - g_2(s,\omega) \big) = a \, \omega^{-\frac{N}{\rho}} \, s^{\rho} - b \, \omega^{-\left(1+\frac{N-1}{\rho}\right)} \, s - c \, \omega^{-\left(1+\frac{N}{\rho}\right)} \; .
	\end{equation*}
	Then the point $s_{\omega} := K_{\rho} \, \omega^{-\frac{1}{\rho}}$, where $K_{\rho}$ is the unique positive solution of the equation $a K^{\rho} - b K - c = 0$, is the unique positive solution of the equation $h_{\omega}(s)=0$. So $g_1(.,\omega)$ and $g_2(.,\omega)$ intersect each other at the point $s_{\omega}$ and we have $g_1(.,\omega) \leqslant g_2(.,\omega)$ for $s \in (0,s_{\omega}]$ and $g_1(.,\omega) \geqslant g_2(.,\omega)$ otherwise. Hence we obtain more precise estimates :
	\begin{equation*}
	\left\{ \begin{array}{rl} \hspace{-0.3cm}
			& \displaystyle \forall \, s \in (0,s_{\omega}] \qquad \big| f(s,\omega) \big| \leqslant a \, \omega^{-\frac{N}{\rho}} = g_1(s,\omega) \; , \\[0.2cm]
			& \displaystyle \forall \, s \in [s_{\omega},+\infty) \qquad \big| f(s,\omega) \big| \leqslant b \, \omega^{-\left(1+\frac{N-1}{\rho}\right)} \, s^{1-\rho} \, + c \, \omega^{-\left(1+\frac{N}{\rho}\right)} \, s^{-\rho} = g_2(s,\omega) \; . \\
		\end{array} \right.
	\end{equation*}
	Now we seek a function $g : (0,+\infty) \times (0,+\infty) \, \longrightarrow \R$ which is locally integrable with respect to the variable $s$ and which satisfies the following inequalities for any $\omega > 0$ :	\begin{equation} \label{int_cond}
		\left\{ \begin{array}{rl} \hspace{-0.3cm}
			& \displaystyle \forall \, s \in (0,s_{\omega}] \qquad \big| f(s,\omega) \big| \leqslant g_1(s,\omega) \leqslant g(s,\omega) \; , \\[0.2cm]
			& \displaystyle \forall \, s \in [s_{\omega},+\infty) \qquad \big| f(s,\omega) \big| \leqslant g_2(s,\omega) \leqslant g(s,\omega) \; , \\
		\end{array} \right.
	\end{equation}
	Here we choose $g_{\gamma,\delta}(s,\omega) := L_{\gamma,\rho} \, s^{-\gamma} \omega^{-\delta}$, where $L_{\gamma,\rho}, \, \delta, \, \gamma > 0$ must be clarified. To this end, we require the following condition :
	\begin{equation*}
		\forall \, \omega > 0 \qquad g_{\gamma,\delta}(s_{\omega},\omega) = g_1(s_{\omega},\omega) = g_2(s_{\omega},\omega) \; ,
	\end{equation*}
	leading to
	\begin{equation*}
		g_{\gamma,\delta}\left(K_{\rho} \, \omega^{-\frac{1}{\rho}},\omega \right) = L_{\gamma,\rho} \, K_{\rho}^{\; -\gamma} \, \omega^{\frac{\gamma}{\rho}-\delta} = a \, \omega^{-\frac{N}{\rho}} \; .
	\end{equation*}
	Since this equality holds for all $\omega > 0$, we obtain
	\begin{equation*}
		\left\{ \begin{array}{rl} \hspace{-0.3cm}
			& \displaystyle L_{\gamma,\rho} = a \, K_{\rho}^{\; \gamma} \\[0.2cm]
			& \displaystyle \frac{\gamma}{\rho} - \delta = -\frac{N}{\rho}
		\end{array} \right. \qquad \Longrightarrow \qquad
		\left\{ \begin{array}{rl} \hspace{-0.3cm}
			& \displaystyle L_{\gamma,\rho} = a \, K_{\rho}^{\; \gamma} \\[0.2cm]
			& \displaystyle \delta = \rho^{-1}(\gamma+N)
		\end{array} \right. \; .
	\end{equation*}
	Here we choose $\gamma \in (0,1)$ so that $g_{\gamma, \delta}(.,\omega) : (0,+\infty) \longrightarrow \R$ is locally integrable with respect to $s$; it follows $\displaystyle \delta = \rho^{-1} (\gamma+N) \in \left(\frac{N}{\rho},\frac{1+N}{\rho} \right)$. To conclude, we have to check the system of inequalities \eqref{int_cond} :
	\begin{itemize}
		\item \textit{Case $s \leqslant s_{\omega}$.} We have
		\begin{equation*}
			g_{\gamma,\delta}(s,\omega) = a \, K_{\rho}^{\; \gamma} \, \omega^{-\delta} \, s^{-\gamma} \geqslant a \, \omega^{\frac{\gamma}{\rho} - \delta} = a \, \omega^{-\frac{N}{\rho}} = g_1(s,\omega) \; ,
		\end{equation*}
		since $\displaystyle \frac{\gamma}{\rho} - \delta = -\frac{N}{\rho}$.
		\item \textit{Case $s \geqslant s_{\omega}$.} Here, we want to show that $g_2(s,\omega) \leqslant g_{\gamma,\delta}(s,\omega)$, that is to say
		\begin{equation} \label{toproove1}
			s^{\rho} \big( g_{\gamma,\delta}(s,\omega) - g_2(s,\omega) \big) = a \, K_{\rho}^{\; \gamma} \, \omega^{-\delta} \, s^{\rho-\gamma} - b \, \omega^{-\left(1+\frac{N-1}{\rho}\right)} \, s - c \, \omega^{-\left(1+\frac{N}{\rho}\right)} \geqslant 0 \; .
		\end{equation}
		We define the function $k_{\omega} : (0,+\infty) \longrightarrow \R$ by $k_{\omega}(s) := s^{\rho}\big(g_{\gamma,\delta}(s,\omega) - g_2(s,\omega)\big)$, and we differentiate it,
		\begin{equation*}
			(k_{\omega})'(s) = a \, K_{\rho}^{\; \gamma} (\rho-\gamma) \, \omega^{-\delta} \, s^{\rho-\gamma -1} - b \, \omega^{-\left(1+\frac{N-1}{\rho}\right)} \; .
		\end{equation*}
		Since $s > 0$ and $\rho \geqslant 2$, $(k_{\omega})'$ is an increasing function and vanishes at the point
		\begin{equation*}
			s_{\omega}' = \left( \frac{b}{a \, K_{\rho}^{\; \gamma} (\rho - \gamma)} \right)^{\frac{1}{\rho - \gamma -1}} \omega^{-\frac{1}{\rho}} \; .
		\end{equation*}
		Now we want to show the inequality : $s_{\omega}' \leqslant s_{\omega}$. To this end, we observe that
		\begin{equation*}
			0 \leqslant b \, K_{\rho} \, (\rho - \gamma - 1) + (\rho-\gamma) \, c \; ,
		\end{equation*}
		because $\rho \geqslant 2$. Furthermore since $K_{\rho}$ satisfies $a \, K_{\rho}^{\; \rho} - b K_{\rho} - c = 0$, we obtain
		\begin{equation*}
			\frac{b K_{\rho}}{\rho-\gamma} \leqslant b K_{\rho} + c = a \, K_{\rho}^{\, \rho} \qquad \Longleftrightarrow \qquad \frac{b}{a \, K_{\rho}^{\; \gamma} (\rho - \gamma)} \leqslant K_{\rho}^{\; \rho - \gamma - 1} \; ,
		\end{equation*}
		and so
		\begin{equation*}
			s_{\omega}' = \bigg(\frac{b}{a \, K_{\rho}^{\; \gamma} (\rho - \gamma)}\bigg)^{\frac{1}{\rho-1-\gamma}} \omega^{-\frac{1}{\rho}} \leqslant K_{\rho} \, \omega^{-\frac{1}{\rho}} = s_{\omega} \; .
		\end{equation*}
		Hence for all $s \geqslant s_{\omega} \geqslant s_{\omega}'$, $k_{\omega}$ is an increasing function and
		\begin{equation*}
			k_{\omega}(s) \geqslant k_{\omega}(s_{\omega}) = s_{\omega}^{\; \rho}\big(g_{\gamma,\delta}(s_{\omega},\omega) - g_2(s_{\omega},\omega)\big) = 0 \; .
		\end{equation*}
		Hence inequality \eqref{toproove1} is satisfied and $g_{\delta,\gamma}(s,\omega) - g_2(s,\omega) \geqslant 0$.
	\end{itemize}
\end{proof}

Supposing that the amplitude is regular at $p_j$, we derive a more precise estimate of the remainder term related to $p_j$ from the two preceding results.

\begin{THM11} \label{THM11}
	Let $N \in \N \backslash \{0\}$ and assume $\mu_j = 1$ and $\rho_j \geqslant 2$ for a certain $j \in \{1,2\}$. Suppose that the functions $\psi : [p_1,p_2] \longrightarrow \R$ and $U : (p_1, p_2) \longrightarrow \C$ satisfy Assumption \emph{(P$_{\rho_1,\rho_2,N}$)} and Assumption \emph{(A$_{\mu_1,\mu_2,N}$)}, respectively. Then the statement of Theorem \ref{THM1} is still true and, for $\gamma \in (0,1)$ and
	\begin{equation*}
		\delta := \frac{\gamma + N}{\rho} \in \left( \frac{N}{\rho_j}, \frac{N+1}{\rho_j} \right) \; ,
	\end{equation*}
	we have a more precise estimate for the remainder term $R_N^{(j)}(\omega)$, namely,
	\begin{equation*}
		\left| R_N^{(j)}(\omega) \right| \leqslant L_{\gamma, \rho_j, N} \int_0^{s_j} s^{-\gamma} \left| \frac{d^N}{ds^N}\left[ \nu_j k_j\right](s) \right| \, ds \; \omega^{-\delta} \; ,
	\end{equation*}
	for all $\omega > 0$, where $L_{\gamma, \rho_j, N} > 0$ is the constant given by Lemma \ref{LEM6} and Lemma \ref{LEM7}.
\end{THM11}

\begin{proof}
	We have to check only the error estimate since the first four steps of the proof of Theorem \ref{THM1} remain valid with $\mu_j = 1$.\\
	Since Lemma \ref{LEM6} ensures that $ \phi_{N}^{(j)}(.,\omega,\rho_j,1)$ satisfies the assumptions of Lemma \ref{LEM7}, we obtain
	\begin{equation*}
		\forall \, s \in (0,s_j] \quad  \forall \,\omega > 0 \qquad \left| \phi_{N}^{(j)}(s,\omega,\rho_j,1) \right| \leqslant L_{\gamma, \rho_j, N} \, s^{-\gamma} \, \omega^{-\delta} \; ,
	\end{equation*}
	where $\gamma, \delta >0$ are defined as above and $L_{\gamma, \rho_j, N} > 0$ is given in Lemma \ref{LEM7}. Using the expression of the remainder term from Theorem \ref{THM1} and the preceding estimate leads to the conclusion, namely,
	\begin{equation*}
		\begin{aligned}
			\left| R_{N}^{(j)}(\omega) \right|	& \leqslant \int_0^{s_j} \left| \, \phi_{N}^{(j)}(s,\omega,\rho_j, 1) \right| \, \left| \frac{d^N}{ds^N}[\nu_j k_j](s) \right| ds \\
												& \leqslant L_{\gamma, \rho_j, N} \int_0^{s_j} s^{-\gamma} \left| \frac{d^N}{ds^N}[\nu_j k_j](s) \right| ds \; \omega^{-\delta} \; .
		\end{aligned}
	\end{equation*}
	And we observe that the decay rate of the remainder term $R_N^{(j)}(\omega)$ with respect to $\omega$ is higher than the one of the highest term of the expansion $A_N^{(j)}(\omega)$.
\end{proof}

\section{Lossless error estimates}

\hspace{2.5ex} In order to motivate the results of this section, let us consider the following oscillatory integral

\begin{equation*}
	I(\omega, p_2) = \int_{p_1}^{p_2} (p-p_1)^{-\frac{1}{4}} \, e^{-i\omega (p-p_2)^2} \, dp \; ,
\end{equation*}
with $\omega > 0$ and $p_2 > p_1$; here $p_1$ is a singularity of the amplitude of order $\mu_1 = \frac{3}{4}$ and $p_0$ a stationary point of order $\rho_2 = 2$. By applying the results of the preceding section, we obtain an asymptotic expansion of the above oscillatory integral, namely,
\begin{align*}
	& \int_{p_1}^{p_2} (p-p_1)^{-\frac{1}{4}} \, e^{-i\omega (p-p_2)^2} \, dp = \frac{\sqrt{\pi}}{2} \, e^{-i\frac{\pi}{4}} \, (p_2-p_1)^{-\frac{1}{4}} \, \omega^{-\frac{1}{2}} \\
	& \hspace{1cm} + \frac{\Gamma\big(\frac{3}{4}\big)}{2^{\frac{1}{4}}} \, e^{i \frac{3 \pi}{8}} \, e^{-i \omega (p_2-p_1)^2} \, (p_2-p_1)^{-\frac{1}{4}} \, \omega^{-\frac{3}{4}} + R_1^{(1)}(\omega,p_2) + R_1^{(2)}(\omega,p_2) \; ,
\end{align*}
and for $\delta \in \big( \frac{3}{4},1 \big)$, we have the following estimates for the remainders :
\begin{align*}
	& \bullet \quad \Big| R_1^{(1)}(\omega,p_2) \Big| \leqslant \int_0^{\frac{8}{9} (p_2-p_1)^2} s^{-\frac{1}{4}} \left| \big(\nu_1 k_1 \big)'(s) \right| ds \, \omega^{-1} ; \\
	& \bullet \quad \Big| R_1^{(2)}(\omega,p_2) \Big| \leqslant L_{\gamma,2,1} \int_0^{\frac{2}{3} (p_2-p_1)} s^{-\gamma} \left| \big(\nu_2 k_2 \big)'(s) \right| ds \, \omega^{-\delta} .
\end{align*}
On one hand, one can remark that $I(\omega, p_2)$ is bounded by $\frac{3}{4} \, (p_2 - p_1)^{\frac{3}{4}}$. On the other hand, if $p_2$ tends to $p_1$ then the expansion of the integral blows up. This implies that this expansion does not furnish a good approximation of $I(\omega, p_2)$ for fixed $\omega > 0$ when $p_2$ is too close to $p_1$.\\
Moreover since the integral is bounded and the expansion blows up, the remainder tends also to infinity when $p_2$ tends to $p_1$. Now let us note that the graphs of the smooth cut-off functions $\nu_j$ compress when $p_2$ tends to $p_1$ implying the fact that the $L^{\infty}$-norm of $(\nu_j)'$ tends to infinity. This observation leads to the idea that the smooth cut-off function contributes artificially to the blow-up of the remainder.\\
Hence the aim of this section is to provide lossless error estimate for the stationary phase method by replacing the smooth cut-off function by a characteristic function.\\

We start by modifying slightly the functions $\varphi_j$ and $k_j$ defined in the previous section by changing only their domains of definition. We shall use the notations of the first section and these new definitions will be used throughout this section.\\

Let $p_1,p_2$ be two finite real numbers such that $p_1 < p_2$, and choose $q \in (p_1 , p_2)$.

\begin{DEF2} \label{DEF2}
	Let $\psi : [p_1,p_2] \longrightarrow \R$ and $U : (p_1, p_2) \longrightarrow \C$ be two functions satisfying Assumption \emph{(P$_{\rho_1,\rho_2,1}$)} and Assumption \emph{(A$_{\mu_1,\mu_2,1}$)}, respectively.
	\begin{enumerate}
	\item For $j=1,2$, let $\varphi_j : I_j \longrightarrow \R$ be the functions defined by
	\begin{equation*}
		\varphi_1(p) := \big( \psi(p) - \psi(p_1) \big)^{\frac{1}{\rho_1}} \qquad , \qquad \varphi_2(p) := \big( \psi(p_2) - \psi(p) \big)^{\frac{1}{\rho_2}} \; ,
	\end{equation*}
	with $I_1 := [p_1, q]$, $I_2 := [q, p_2]$ and $s_1 := \varphi_1(q)$, $s_2 := \varphi_2(q)$.
	\item For $j=1,2$, let $k_j : (0,s_j] \longrightarrow \C$ be the functions defined by
	\begin{equation*}
			k_j(s) := U\big( \varphi_j^{-1}(s) \big) \, s^{1- \mu_j} \, \big( \varphi_j^{-1} \big)'(s) \; ,
	\end{equation*}
	which can be extended to $[0,s_j]$ (see Proposition \ref{LEM2}).
	\end{enumerate}
\end{DEF2}

Now we state and prove a refinement of the version of the stationary phase method of Erdélyi which consists in replacing the smooth cut-off function by a characteristic function. The hypotheses on the regularity of the phase and of the amplitude are weakened as compared with Theorem \ref{THM1}, because we establish an expansion to one term only which requires a single integration by parts.

\begin{THM2} \label{THM2}
	Let $\rho_1, \rho_2 \geqslant 1$ and $\mu_1, \mu_2 \in (0,1)$. Suppose that the functions $\psi : [p_1,p_2] \longrightarrow \R$ and $U : (p_1, p_2) \longrightarrow \C$ satisfy Assumption \emph{(P$_{\rho_1,\rho_2,1}$)} and Assumption \emph{(A$_{\mu_1,\mu_2,1}$)}, respectively.\\ Then for $j=1,2$, there are functions $A^{(j)}$, $R_1^{(j)}(.,q)$, $R_2^{(j)}(.,q) : (0,+\infty) \longrightarrow \C$ such that
	\begin{equation*}
		\left\{ \begin{array}{rl}
			& \displaystyle \int_{p_1}^{p_2} U(p) \, e^{i \omega \psi(p)} \, dp = \sum_{j=1,2} \left( A^{(j)}(\omega) + R_1^{(j)}(\omega,q) + R_2^{(j)}(\omega,q) \right) \; , \\[2mm]
			& \displaystyle \left| R_1^{(j)}(\omega,q) \right| \leqslant \frac{1}{\rho_j} \, \Gamma \bigg(\frac{1}{\rho_j}\bigg) \int_0^{s_j} s^{\mu_j-1} \big|(k_j)'(s) \big| \, ds \; \omega^{-\frac{1}{\rho_j}} \; , \\
			& \displaystyle \left| R_2^{(j)}(\omega,q) \right| \leqslant \frac{\rho_j-\mu_j}{\rho_j} \, \Gamma \hspace{-1mm} \left(\frac{1}{\rho_j}\right) \, \Big| U(q) \, (\varphi_j)'(q)^{-1} \Big| \, \varphi_j(q)^{-\rho_j} \, \omega^{-\left(1+\frac{1}{\rho_j}\right)} \; ,
		\end{array} \right.
	\end{equation*}
	for all $\omega > 0$ and for a fixed $q \in (p_1,p_2)$, where
	\begin{equation*}
		\begin{aligned}
			\bullet \quad & A^{(j)}(\omega) :=  e^{i \omega \psi(p_j)} \, k_j(0) \, \Theta^{(j)}(\rho_j,\mu_j) \, \omega^{-\frac{\mu_j}{\rho_j}} \; , \\
			\bullet \quad & R_1^{(j)}(\omega,q) := (-1)^{j} \, e^{i\omega \psi(p_j)} \int_0^{s_j} \phi^{(j)}(s,\omega, \rho_j, \mu_j) \, (k_j)'(s) \, ds \; , \\
			\bullet \quad & R_2^{(j)}(\omega,q) := (-1)^j \, i \, \frac{\mu_j-\rho_j}{\rho_j} \, e^{i \omega \psi(p_j)} \, k_j(s_j) \int_{\Lambda^{(j)}(s_j)} z^{\mu_j-\rho_j-1} e^{(-1)^{j+1} i \omega z^{\rho_j}} dz \; \omega^{-1} \; , \\
			\bullet \quad & \Theta^{(j)}(\rho_j,\mu_j) := \frac{(-1)^{j+1}}{\rho_j} \, \Gamma \bigg(\frac{\mu_j}{\rho_j}\bigg) \, e^{(-1)^{j+1} i \frac{\pi}{2} \, \frac{\mu_j}{\rho_j}} \; ,  \\
			\bullet \quad & \phi^{(j)}(s,\omega,\rho_j, \mu_j) := - \int_{\Lambda^{(j)}(s)} z^{\mu_j -1} e^{(-1)^{j+1} i\omega z^{\rho_j}} \, dz \; .
		\end{aligned}
	\end{equation*}
\end{THM2}

\begin{REM7}
	\emph{Note that the quantities $\Theta^{(j)}(\rho_j,\mu_j)$ and $\phi^{(j)}(s,\omega,\rho_j, \mu_j)$ correspond to $\Theta_1^{(j)}(\rho_j,\mu_j)$ and $\phi_1^{(j)}(s,\omega,\rho_j, \mu_j)$, respectively, which are defined in Theorem \ref{THM1}. We have removed the subscripts for simplicity.}
\end{REM7}

\begin{proof}[Proof of Theorem \ref{THM2}]
	This proof follows the steps of the proof of Theorem \ref{THM1}. Hence the steps which are identical will be only sketched and we shall focus on the new arguments coming from the cutting-point $q$.\\
	As previously, we shall note $\phi^{(j)}(s,\omega)$ instead of $\phi^{(j)}(s,\omega,\rho_j,\mu_j)$ in the proof.\\
	
	\noindent \textit{First step: Splitting of the integral.} We fix a point $q \in (p_1,p_2)$ and we split the integral at this point,
	\begin{equation*}
		\int_{p_1}^{p_2} U(p) \, e^{i \omega \psi(p)} \, dp = \tilde{I}^{(1)}(\omega,q) + \tilde{I}^{(2)}(\omega,q) \; ,
	\end{equation*}
	where
	\begin{equation*}
		\tilde{I}^{(1)}(\omega,q) := \int_{p_1}^{q} U(p) \, e^{i \omega \psi(p)} \, dp \quad , \quad \tilde{I}^{(2)}(\omega,q) := \int_{q}^{p_2} U(p) \, e^{i \omega \psi(p)} \, dp \; . \\ \\
	\end{equation*}
	\textit{Second step: Substitution.} Since $\varphi_1 : [p_1,q] \longrightarrow [0,s_1]$ is a $\mathcal{C}^{2}$-diffeomorphism, we obtain by setting $s = \varphi_1(p)$,
	\begin{equation*}
		\tilde{I}^{(1)}(\omega,q) = e^{i \omega \psi(p_1)} \int_0^{s_1} k_1(s) \, s^{\mu_1-1} e^{i \omega s^{\rho_1}} \, ds \; ,
	\end{equation*}
	where $k_1$ is given in Definition \ref{DEF2}. In a similar way, we obtain
	\begin{equation*}
		\tilde{I}^{(2)}(\omega,q) = - e^{i \omega \psi(p_2)} \int_0^{s_2} k_2(s) \, s^{\mu_2-1} e^{-i \omega s^{\rho_2}} \, ds \; ,
	\end{equation*}
	with $k_2$ defined above.\\
	
	\noindent \textit{Third step: Integration by parts.} An integration by parts leads to
	\begin{equation*}
		\begin{aligned}
			\tilde{I}^{(1)}(\omega,q)	& = \phi^{(1)}(s_1,\omega) \, k_1(s_1) \, e^{i \omega \psi(p_1)} - \phi^{(1)}(0,\omega) \, k_1(0) \, e^{i \omega \psi(p_1)}\\
							& \qquad \qquad - e^{i \omega \psi(p_1)} \int_0^{s_1} \phi^{(1)}(s,\omega) \, (k_1)'(s) \, ds \; ,
		\end{aligned}
	\end{equation*}
	and similarly,
	\begin{align}
			\tilde{I}^{(2)}(\omega,q)	& = \phi^{(2)}(0,\omega) \, k_2(0) \, e^{i \omega \psi(p_2)} - \phi^{(2)}(s_2,\omega) \, k_2(s_2) \, e^{i \omega \psi(p_2)} \nonumber \\
								& \label{I_2} \qquad \qquad + e^{i \omega \psi(p_2)} \int_0^{s_2} \phi^{(2)}(s,\omega) (k_2)'(s) \, ds \; .
	\end{align}
	
	\noindent \textit{Fourth step: Cancellation.} The aim of this step is to simplify the difference:
	\begin{equation} \label{difference}
		\phi^{(1)}(s_1,\omega) \, k_1 (s_1) \, e^{i \omega \psi(p_1)} - \phi^{(2)}(s_2,\omega) \, k_2 (s_2) \, e^{i \omega \psi(p_2)} \; .
	\end{equation}
	Since the two functions $s \longmapsto \phi^{(j)}(s,\omega)$ are given by oscillatory integrals, we can expand them with respect to $\omega$ and we shall show that the first terms cancel out.\\ Since $s_1 > 0$, we note that the derivative of the function $t \longmapsto \left(s_1+t e^{i \frac{\pi}{2\rho_1}} \right)^{\rho_1}$ does not vanish for all $t > 0$  and one can write
	\begin{equation*}
		e^{i \omega \left(s_1+t e^{i \frac{\pi}{2\rho_1}}\right)^{\rho_1}} = e^{-i \frac{\pi}{2\rho_1}} \, (i \omega \rho_1)^{-1} \left(s_1+t e^{i \frac{\pi}{2\rho_1}}\right)^{1-\rho_1} \frac{d}{dt} \hspace{-1mm} \left[e^{i \omega \left(s_1+t e^{i \frac{\pi}{2\rho_1}}\right)^{\rho_1}}\right] \; .
	\end{equation*}
	Putting this equality in the definition of $\phi^{(1)}(s_1,\omega)$ and carrying out an integration by parts lead to
	\begin{align}
			\phi^{(1)}(s_1,\omega)	& = - (i \omega \rho_1)^{-1} \int_0^{+\infty} \left(s_1+te^{i \frac{\pi}{2 \rho_1}}\right)^{\mu_1-\rho_1} \frac{d}{dt}\left[e^{i \omega \left(s_1+te^{i \frac{\pi}{2 \rho_1}}\right)^{\rho_1}}\right] \, dt \nonumber \\
														& = (i \omega \rho_1)^{-1} s_1^{\; \mu_1-\rho_1} \, e^{i \omega s_1^{\; \rho_1}} \nonumber \\
														& \label{expphi1} \qquad + \frac{\mu_1-\rho_1}{i \omega \rho_1} \, e^{i \frac{\pi}{2 \rho_1}} \int_0^{+\infty} \left(s_1+te^{i \frac{\pi}{2 \rho_1}}\right)^{\mu_1-\rho_1-1} e^{i \omega \left(s_1+te^{i \frac{\pi}{2 \rho_1}}\right)^{\rho_1}} \, dt \; .
	\end{align}
	We remark that the boundary term at infinity is $0$; indeed, we observe that
	\begin{equation*}
		s_1 \leqslant \left| s_1+te^{i \frac{\pi}{2 \rho_1}} \right| \qquad \Longrightarrow \qquad s_1^{\; \mu_1-\rho_1} \geqslant \left| s_1+te^{i \frac{\pi}{2 \rho_1}} \right|^{\mu_1-\rho_1} \; ,
	\end{equation*}
	because $\mu_1 \leqslant \rho_1$, and by using Lemma \ref{LEM3}, we obtain
	\begin{equation*}
		\forall \, t > 0 \qquad \left| \left(s_1+te^{i \frac{\pi}{2 \rho_1}}\right)^{\mu_1-\rho_1} e^{i \omega \left(s_1+te^{i \frac{\pi}{2 \rho_1}}\right)^{\rho_1}} \right| \leqslant s_1^{\; \mu_1-\rho_1} e^{-\omega t^{\rho_1}} \; \longrightarrow \; 0 \quad , \quad t \rightarrow + \infty \; .
	\end{equation*}
	In a similar way, we have
	\begin{align*}
			\phi^{(2)}(s_2,\omega)	& = -(i \omega \rho_2)^{-1} s_2^{\; \mu_2-\rho_2} \, e^{-i \omega s_2^{\rho_2}} \\
														& \qquad - \frac{\mu_2-\rho_2}{i \omega \rho_2} \, e^{-i \frac{\pi}{2 \rho_2}} \int_0^{+\infty} \left(s_2+te^{-i \frac{\pi}{2 \rho_2}}\right)^{\mu_2-\rho_2-1} e^{-i \omega \left(s_2+te^{-i \frac{\pi}{2 \rho_2}}\right)^{\rho_2}} \, dt \; .
	\end{align*}
	Furthermore, by the definitions of $k_j$ and $s_j := \varphi_j(q)$, we obtain
	\begin{equation*}
		k_j(s_j) = U\big( \varphi_j^{-1}(s_j) \big) \, s_j^{1-\mu_j} \, (\varphi_j^{-1})'(s_j) = U(q) \, \varphi_j(q)^{1-\mu_j} \, (\varphi_j)'(q)^{-1} \; .
	\end{equation*}
	Now we multiply the last expression in the case $j=1$ by the expansion of $\phi^{(1)}(s_1,\omega)$ given by \eqref{expphi1},
	\begin{equation} \label{simpphi1}
		\begin{aligned}
			\phi^{(1)} \, (s_1,\omega) & \, k_1(s_1) \, e^{i \omega \psi(p_1)} = (i \omega \rho_1)^{-1} \, e^{i \omega \left(s_1^{\; \rho_1} + \psi(p_1)\right)} \, U(q) \, \varphi_1(q)^{1-\rho_1} \, (\varphi_1)'(q)^{-1} \\
									& - i \, \frac{\mu_1 - \rho_1}{\rho_1} \, e^{i \omega \psi(p_1)} \, e^{i \frac{\pi}{2 \rho_1}} \, k_1(s_1) \int_0^{+\infty} \left(s_1+te^{i \frac{\pi}{2 \rho_1}}\right)^{\mu_1-\rho_1-1} e^{i \omega \left(s_1+te^{i \frac{\pi}{2 \rho_1}}\right)^{\rho_1}} \, dt \; \omega^{-1} \; .
		\end{aligned}
	\end{equation}
	The definition of $\varphi_1(q)$ gives $\varphi_1(q)^{\rho_1} = \psi(q) - \psi(p_1)$ and by the regularity of $\varphi_1$, one has	\begin{equation*}
		\rho_1 \, (\varphi_1)'(q) \, \varphi_1(q)^{\rho_1-1} = \frac{d}{dp}\Big[ (\varphi_1)^{\rho_1} \Big](q) = \psi'(q) \; ,
	\end{equation*}
	simplifying the first term in \eqref{simpphi1}; moreover the integral in \eqref{simpphi1} can be written as an integral on the curve $\Lambda^{(1)}(s_1)$ in the complex plane. These considerations lead to
	\begin{align*}
			\phi^{(1)}(s_1,\omega) \, & k_1(s_1)  \, e^{i \omega \psi(p_1)} = -i \omega^{-1} \, e^{i \omega \psi(q)} \, \frac{U(q)}{\psi'(q)} \\
										& \qquad - i \, \frac{\mu_1 - \rho_1}{\rho_1}  \, e^{i \omega \psi(p_1)} \, k_1(s_1)  \int_{\Lambda^{(1)}(s_1)} z^{\mu_1-\rho_1-1} e^{i \omega z^{\rho_1}} \, dz \; \omega^{-1} \; .
	\end{align*}
	By similar calculations, we obtain
	\begin{align*}
			\phi^{(2)} \, (s_2,\omega) & \, k_2(s_2) \, e^{i \omega \psi(p_2)}	= -i \omega^{-1} \, e^{i \omega \psi(q)} \, \frac{U(q)}{\psi'(q)} \\
										& \qquad - i \, \frac{\mu_2 - \rho_2}{\rho_2} \, e^{i \omega \psi(p_2)} \, k_2(s_2) \int_{\Lambda^{(2)}(s_2)} z^{\mu_2-\rho_2-1} e^{-i \omega z^{\rho_2}} \, dz \; \omega^{-1} \; .
	\end{align*}
	Hence we remark that the difference \eqref{difference} is  equal to
	\begin{equation*}
		\sum_{j = 1}^2 \, (-1)^{j} \, i \, \frac{\mu_j - \rho_j}{\rho_j} \, e^{i \omega \psi(p_j)} \,  k_j(s_j) \int_{\Lambda^{(j)}(s_j)} z^{\mu_j-\rho_j-1} \, e^{(-1)^{j+1} i \omega z^{\rho_j}} \, dz \; \omega^{-1} \; .	
	\end{equation*}
	Consequently, we are able to write the initial integral as follows,
	\begin{equation*}
		\begin{aligned}
			\int_{p_1}^{p_2} U(p) e^{i \omega \psi(p)} dp	& = - \phi^{(1)}(0,\omega) \, k_1(0) \, e^{i \omega \psi(p_1)} - e^{i \omega \psi(p_1)} \int_0^{s_1} \phi^{(1)}(s,\omega) (k_1)'(s) \, ds \\
							& \qquad \qquad - i \, \frac{\mu_1-\rho_1}{\rho_1}  \, e^{i \omega \psi(p_1)} \, k_1(s_1)  \int_{\Lambda^{(1)}(s_1)} z^{\mu_1-\rho_1-1} e^{i \omega z^{\rho_1}} \, dz \; \omega^{-1} \\
							& \quad + \phi^{(2)}(0,\omega) \, k_2(0) \, e^{i \omega \psi(p_2)} + e^{i \omega \psi(p_2)} \int_0^{s_2} \phi^{(2)}(s,\omega) (k_2)'(s) \, ds \\
							& \qquad \qquad + i \, \frac{\mu_2-\rho_2}{\rho_2} \, e^{i \omega \psi(p_2)} \, k_2(s_2) \int_{\Lambda^{(2)}(s_2)} z^{\mu_2-\rho_2-1} e^{-i \omega z^{\rho_2}} \, dz \; \omega^{-1} \\
							& =: \sum_{j=1,2} \left( A^{(j)}(\omega) + R_1^{(j)}(\omega,q) + R_2^{(j)}(\omega,q) \right) \; .
		\end{aligned}
	\end{equation*}
	According to the fourth step of the proof of Theorem \ref{THM1}, we have
	\begin{equation*}
		\phi^{(j)}(0,\omega) = -\frac{1}{\rho_j} \, \Gamma \hspace{-1mm} \left(\frac{\mu_j}{\rho_j}\right) \, e^{(-1)^{j+1} i \frac{\pi}{2} \, \frac{\mu_j}{\rho_j}} \, \omega^{-\frac{\mu_j}{\rho_j}} =: (-1)^j \, \Theta^{(j)}(\rho_j,\mu_j) \, \omega^{-\frac{\mu_j}{\rho_j}}  \; .
	\end{equation*}
	leading to the definition of $A^{(j)}(\omega)$,
	\begin{equation*}
		A^{(j)}(\omega) := (-1)^{j} \, \phi^{(j)}(0,\omega) \, k_j(0) \, e^{i \omega \psi(p_j)} = 	e^{i \omega \psi(p_j)} \, k_j(0) \, \Theta^{(j)}(\rho_j,\mu_j) \, \omega^{-\frac{\mu_j}{\rho_j}} \; .
	\end{equation*}
	
	\noindent \textit{Fifth step: Remainder estimates.} Using the last step of the proof of Theorem \ref{THM1}, we obtain
	\begin{equation*}
		\left| e^{i \omega \psi(p_j)} \int_0^{s_j} \phi^{(j)}(s,\omega, \rho_j, \mu_j) \, (k_j)'(s) \, ds \right| \leqslant \frac{1}{\rho_j} \, \Gamma \hspace{-1mm} \left(\frac{1}{\rho_j}\right) \int_0^{s_j} s^{\mu_j-1} \, \big| (k_1)'(s) \big| \, ds \, \omega^{-\frac{1}{\rho_j}} \; .
	\end{equation*}
	Now let us estimate $R_2^{(j)}(\omega,q)$. We have
	\begin{align}
		\Bigg| i \, \frac{\mu_j-\rho_j}{\omega \rho_j}	& \, e^{i \omega \psi(p_j)} \, k_j(s_j) \, e^{(-1)^{j+1} i \frac{\pi}{2 \rho_j}} \nonumber \\
												& \qquad \times \int_0^{+\infty} \left(s_j+te^{(-1)^{j+1} i \frac{\pi}{2\rho_j}}\right)^{\mu_j-\rho_j-1} e^{(-1)^{j+1} i \omega \left(s_j+te^{(-1)^{j+1}i \frac{\pi}{2 \rho_j}}\right)^{\rho_j}} \, dt \, \Bigg| \nonumber \\
												& \label{ineq1} \leqslant \frac{\rho_j-\mu_j}{\rho_j} \, \big| k_j(s_j) \big| \, \omega^{-1} \,  s_j^{\; \mu_j-\rho_j-1} \int_0^{+\infty} \left| e^{(-1)^{j+1} i \omega \left(s_j+te^{(-1)^{j+1}i \frac{\pi}{2 \rho_j}}\right)^{\rho_j}} \right| \, dt \\
												& \label{ineq2} \leqslant \frac{\rho_j-\mu_j}{\rho_j} \, \Big| U \big( \varphi_j^{-1}(s_j)\big) \, (\varphi_j^{-1})'(s_j) \Big| \, s_j^{\; -\rho_j} \, \omega^{-1} \int_0^{+\infty} e^{-\omega t^{\rho_j}} \, dt \\
												& \label{ineq3} = \frac{\rho_j-\mu_j}{\rho_j} \, \Gamma \hspace{-1mm} \left(\frac{1}{\rho_j}\right) \, \big| U(q) \, (\varphi_j)'(q)^{-1} \big| \, \varphi_j(q)^{-\rho_j} \,  \omega^{-\left(1+\frac{1}{\rho_j}\right)} \; ;
	\end{align}
	\begin{itemize}
		\item \eqref{ineq1}: use the fact $\displaystyle s_j \leqslant \left| s_j+te^{(-1)^{j+1}i \frac{\pi}{2 \rho_j}} \right|$ ;
		\item \eqref{ineq2}: employ the definition of the function $k_j$ and Lemma \ref{LEM3} ;
		\item \eqref{ineq3}: use the equalities $\displaystyle \int_0^{+\infty} e^{-\omega t^{\rho_j}} dt = \Gamma \hspace{-1mm} \left(\frac{1}{\rho_j} \right) \omega^{-\frac{1}{\rho_j}}$ and $q = \varphi_j^{-1}(s_j)$.
	\end{itemize}
	We remark finally that the decay rates of $A^{(j)}(\omega)$, $R_1^{(j)}(\omega,q)$ and $R_2^{(j)}(\omega,q)$ are $\omega^{-\frac{\mu_j}{\rho_j}}$, $\omega^{-\frac{1}{\rho_j}}$ and $\omega^{-\big(1+\frac{1}{\rho_j}\big)}$, respectively. Thus the decay rates of the remainder related to $p_j$ and of the remainder related to $q$ are higher than the one of the first term related to $p_j$. This ends the proof.
\end{proof}

For fixed $q \in (p_1,p_2)$, we observe that $R_2^{(j)}(\omega,q)$ (for $j=1,2$) is always negligible as compared with $A^{(j)}(\omega)$ when $\omega$ tends to infinity, even if $\mu_j=1$. Neverthless if $\mu_j=1$  then the decay rates of $R_1^{(j)}(\omega,q)$ and $A^{(j)}(\omega)$ with respect to $\omega$ are the same. So we shall use the ideas of the proof of Theorem \ref{THM11} to obtain a better decay rate for $R_1^{(j)}(\omega,q)$.

\begin{THM3} \label{THM3}
	Assume $\mu_j = 1$ and $\rho_j \geqslant 2$ for a certain $j \in \{1,2\}$. Suppose that the functions $\psi : [p_1,p_2] \longrightarrow \R$ and $U : (p_1, p_2) \longrightarrow \C$ satisfy Assumption \emph{(P$_{\rho_1,\rho_2,1}$)} and Assumption \emph{(A$_{\mu_1,\mu_2,1}$)}, respectively. Then the statement of Theorem \ref{THM2} is still true and, for $\gamma \in (0,1)$ and
	\begin{equation*}
		\delta := \frac{\gamma + 1}{\rho} \in \left( \frac{1}{\rho_j}, \frac{2}{\rho_j} \right) \; ,
	\end{equation*}
	we have a more precise estimate for the remainder term $R_1^{(j)}(\omega,q)$, namely,
	\begin{equation*}
		\forall \, \omega > 0 \qquad \left| R_1^{(j)}(\omega,q) \right| \leqslant L_{\gamma, \rho_j} \int_0^{s_j} s^{-\gamma} \, \big| (k_j)'(s) \big| \, ds \; \omega^{-\delta} \; ,
	\end{equation*}
	where $L_{\gamma, \rho_j} := L_{\gamma, \rho_j, 1} > 0$ with $L_{\gamma, \rho_j, 1}$ is given in Theorem \ref{THM11}.
\end{THM3}

\begin{proof}
	We obtain the above estimate by following the lines of the proof of Theorem \ref{THM11}.
\end{proof}

\section{Approaching stationary points and amplitude singularities: the first and error term between blow-up and decay}

\hspace{2.5ex} In this section, we consider a family of oscillatory integrals with respect to a large parameter $\omega$. In view of applications to the solution formula of the free Schrödinger equation, we suppose that the phase function is a polynomial of degree $2$ and has its stationary point $p_0$ inside $(p_1,p_2)$, which contains the support of the amplitude. We suppose in addition that the amplitude has a singularity at $p_1$, the left endpoint of the  integration interval.\\
The aim of this section is to furnish remainder estimates with explicit blow-up and to exploit this in order to find curves in the parameter space on which blow-up and decay balance out.\\

In the first result, we split the integral at the stationary point and then we expand the two new integrals by using the results of the previous section.

\begin{THM4} \label{THM4}
	Let $p_1 < p_2$ be two finite real numbers. Let $p_0 \in (p_1,p_2)$ and $c \in \R$ be two parameters and define $\psi : [p_1,p_2] \longrightarrow \R$ by
	\begin{equation*}
		\psi(p) := - (p-p_0 )^2 + c \; .
	\end{equation*}
	Define the following integrals for all $\omega > 0$,
	\begin{equation*}
		I^{(1)}(\omega,p_0) := \int_{p_1}^{p_0} U(p) \, e^{i \omega \psi(p)} \, dp \qquad , \qquad I^{(2)}(\omega,p_0) := \int_{p_0}^{p_2} U(p) \, e^{i \omega \psi(p)} \, dp \; ,
	\end{equation*}
	where $U$ satisfies Assumption \emph{(A$_{\mu,1,1}$)} on $[p_1,p_2]$ with $\mu \in (0,1)$, and $U(p_2)=0$. Let us define $\tilde{H}(\omega,\psi,U)$ and $\tilde{K}_{\mu}(\omega,\psi,U)$ as follows,
	\begin{equation*}
		\tilde{H}(\omega,\psi,U) := \sqrt{\pi} \, e^{-i \frac{\pi}{4}} \, e^{i \omega c} \, \tilde{u}(p_0) \qquad , \qquad \tilde{K}_{\mu}(\omega, \psi, U) := \frac{\Gamma(\mu)}{2^{\mu}} \, e^{i \frac{\pi \mu}{2}} \, e^{i \omega \psi(p_1)} \, \tilde{u}(p_1) \; .
	\end{equation*}
	Then
	\begin{itemize}
		\item we have
		\begin{equation*}
			\begin{aligned}
				\Big| I^{(1)}(\omega,p_0) \, - \, \tilde{K}_{\mu}(\omega, \psi,U) \, (p_0 - p_1)^{-\mu} \, \omega^{-\mu}  \, -	& \, \frac{1}{2} \, \tilde{H}(\omega,\psi,U) \, (p_0 - p_1)^{\mu - 1} \, \omega^{-\frac{1}{2}} \Big| \\
											& \leqslant \sum_{k = 1}^6 \, R_k^{(1)}(U) \, (p_0 - p_1)^{-\alpha_k^{(1)}} \, \omega^{-\beta_k^{(1)}} \; ,
			\end{aligned}
		\end{equation*}
		where the constants $R_k^{(1)}(U) \geqslant 0$ and the exponents $\alpha_k^{(1)} \in \R$, $\beta_k^{(1)} > 0$ are given in the proof ;
		\item we have
		\begin{equation*}
			\Big| I^{(2)}(\omega,p_0) \, -	\, \frac{1}{2} \, \tilde{H}(\omega,\psi,U) \, (p_0 - p_1)^{\mu - 1} \, \omega^{-\frac{1}{2}} \Big| \leqslant \sum_{k = 1}^2 \, R_k^{(2)}(U) \, (p_0 - p_1)^{-\alpha_k^{(2)}} \, \omega^{-\beta_k^{(2)}} \; ,
		\end{equation*}
		where the constants $R_k^{(2)}(U) \geqslant 0$ and the exponents $\alpha_k^{(2)} \in \R$, $\beta_k^{(2)} > 0$ are given in the proof.
	\end{itemize}
\end{THM4}

\begin{proof}
	\noindent \textit{Study of $I^{(1)}(\omega,p_0)$}. For all $p \in [p_1,p_0]$, we have
	\begin{equation*}
		\psi'(p) = 2 (p_0-p) \; .
	\end{equation*}
	By setting $\tilde{\psi} := 2$, we observe that $\psi$ verifies Assumption \emph{(P$_{1,2,N}$)} on $[p_1,p_0]$, for all $N \geqslant 1$. Then Theorem \ref{THM2} is applicable. Here we choose $q := q(p_0) = p_1 + \frac{p_0-p_1}{2} = p_0 - \frac{p_0-p_1}{2}$ for simplicity. Then we obtain
	\begin{equation*}
		I^{(1)}(\omega,p_0) = \sum_{j=1,2} \left( A^{(j)}(\omega,p_0) + R_1^{(j)}(\omega,p_0) + R_2^{(j)}(\omega,p_0) \right) \; ,	\end{equation*}
	with
	\begin{equation*}
		\begin{aligned}
			& \bullet \quad A^{(1)}(\omega,p_0) = e^{i \omega \psi(p_1)} \, k_1(0) \, \Theta^{(1)}(1,\mu) \, \omega^{-\mu} = \Gamma(\mu) \, e^{i \frac{\pi \mu}{2}} \, e^{i \omega \psi(p_1)} \, k_1(0) \, \omega^{-\mu} \; , \\
			& \bullet \quad A^{(2)}(\omega,p_0) = e^{i \omega \psi(p_0)} \, k_2(0) \, \Theta^{(2)}(2,1) \, \omega^{-\frac{1}{2}} = - \frac{\sqrt{\pi}}{2} \, e^{-i \frac{\pi}{4}} \, e^{i \omega \psi(p_0)} \, k_2(0) \, \omega^{-\frac{1}{2}} \; .
		\end{aligned}
	\end{equation*}
	To compute the values of $k_1(0)$ and $k_2(0)$, let us study the functions $(\varphi_1^{-1})'$ and $(\varphi_2^{-1})'$. On the one hand, we obtain by the simple definition of $\varphi_1$,
	\begin{equation} \label{phi122}
		\varphi_1(p) = \psi(p) - \psi(p_1) \qquad \Longrightarrow \qquad (\varphi_1)'(p) = \psi'(p) = 2 (p_0-p) \; ,
	\end{equation}
	for all $p \in [p_1,q]$. On the other hand, by the definition of $\varphi_2$ and the expression of $\psi$, one has
	\begin{equation*}
		\varphi_2(p) = \big( \psi(p_0) - \psi(p) \big)^{\frac{1}{2}} = (p_0-p) \; ,
	\end{equation*}
	for every $p \in [q,p_0]$. So $(\varphi_2^{-1})'(s) = -1$ and $(\varphi_2^{-1})''(s) = 0$ for all $s \in [0,s_2]$. Then it is possible to compute $k_1(0)$ by using the representation of $k_1$ given in the proof of Proposition \ref{LEM2},
	\begin{align}
		\label{k_1} k_1(s)	& = \left( \int_0^1 (\varphi_1^{-1})'(sy) \, dy \right)^{\mu-1} \tilde{u}\big(\varphi_1^{-1}(s) \big) (\varphi_1^{-1})'(s) \\
				& \underset{s \longrightarrow 0^+}{\longrightarrow} \; (\varphi_1^{-1})'(0)^{\mu} \, \tilde{u}\big( \varphi_1^{-1}(0) \big) = \psi'(p_1)^{-\mu} \, \tilde{u}(p_1) = \big( 2 (p_0 - p_1) \big)^{-\mu} \, \tilde{u}(p_1) \; . \nonumber
	\end{align}
	The calculation of $k_2(0)$ is easier than above since $p_0$ is a regular point of the amplitude :
	\begin{equation*}
		k_2(s) = U\big( \varphi_2^{-1}(s) \big) (\varphi_2^{-1})'(s) = -U\big( \varphi_2^{-1}(s) \big) \; \underset{s \longrightarrow 0^+}{\longrightarrow} \; -U(p_0) = -(p_0 - p_1)^{\mu-1} \, \tilde{u}(p_0) \; .
	\end{equation*}
	Therefore we obtain
	\begin{align*}
		& \bullet \quad A^{(1)}(\omega,p_0) = \frac{\Gamma(\mu)}{2^{\mu}} \, e^{i \frac{\pi \mu}{2}} \, e^{i \omega \psi(p_1)} \, \tilde{u}(p_1) \, (p_0 - p_1)^{-\mu} \, \omega^{-\mu} \; , \\
		& \bullet \quad A^{(2)}(\omega,p_0) = \frac{\sqrt{\pi}}{2} \, e^{-i \frac{\pi}{4}} \, e^{i \omega c} \, \tilde{u}(p_0) \, (p_0 - p_1)^{\mu-1} \, \omega^{-\frac{1}{2}} \; .
	\end{align*}
	
	Now let us control precisely the remainder terms. To do so, we have to study the functions $k_j$ and $\varphi_j$. Firstly, employing the equalities $\frac{1}{2}(p_0-p_1) = p_0-q = q-p_1$ and \eqref{phi122} leads to
	\begin{equation} \label{24}
		(p_0-p_1) = 2(p_0-q) \leqslant (\varphi_1)'(p) \leqslant 2 (p_0-p_1) \; ,
	\end{equation}
	for all $p \in [p_1,q]$. It follows
	\begin{equation*}
		\forall \,  s \in [0,s_1] \qquad \big(2 (p_0-p_1) \big)^{-1} \leqslant (\varphi_1^{-1})'(s) \leqslant (p_0 - p_1)^{-1} \; .
	\end{equation*}
	Moreover by the equality $(\varphi_1^{-1})''(s) = - (\varphi_1)''\big(\varphi_1^{-1}(s)\big) \, (\varphi_1^{-1})'(s)^3$, we have
	\begin{equation*}
		\forall \, s \in [0,s_1] \qquad (\varphi_1^{-1})''(s) = 2 \, (\varphi_1^{-1})'(s)^3 \leqslant 2 \, (p_0-p_1)^{-3} \; .
	\end{equation*}
	Then it is possible compute the value of $s_1$ by using its definition,
	\begin{equation} \label{s1}
		s_1 = \varphi_1(q) = \psi(q) - \psi(p_0) = \frac{3}{4} \, (p_0 - p_1)^2 \leqslant (p_0 - p_1)^2 \; .
	\end{equation}
	Concerning $s_2$, we have
	\begin{equation} \label{s2}
		s_2 = \varphi_2(q) = (p_0 - q) = \frac{1}{2} \, (p_0-p_1) \leqslant p_0 - p_1 \; .
	\end{equation}
	Now we study the functions $k_j$. For this purpose, we shall use the expression of $k_1$ given in \eqref{k_1}. Since $\psi$ satisfies Assumption (P$_{1,2,N}$) on $[p_1,p_0]$ for all $N \geqslant 1$, it follows that $\varphi_1$ is a $\mathcal{C}^{N+1}$-diffeomorphism by Proposition \ref{LEM1}. Thus one has the ability to differentiate under the integral sign the function $\displaystyle s \longmapsto \int_0^1 (\varphi_1^{-1})'(sy) dy$. Hence for all $s \in [0,s_1]$,
	\begin{equation*}
		\begin{aligned}
			(k_1)'(s) =	& (\mu-1) \left( \int_0^1 y \, (\varphi_1^{-1})''(sy) \, dy \right) \left( \int_0^1 (\varphi_1^{-1})'(sy) \, dy \right)^{\mu-2} \tilde{u}\big(\varphi_1^{-1}(s) \big) \, (\varphi_1^{-1})'(s) \\
						& \qquad + \left( \int_0^1 (\varphi_1^{-1})'(sy) \, dy \right)^{\mu-1} \tilde{u}'\big(\varphi_1^{-1}(s) \big) \, (\varphi_1^{-1})'(s)^2 \\
						&  \qquad \qquad + \left( \int_0^1 (\varphi_1^{-1})'(sy) \, dy \right)^{\mu-1} \tilde{u}\big(\varphi_1^{-1}(s) \big) \, (\varphi_1^{-1})''(s) \; .
		\end{aligned}
	\end{equation*}
	In consequence, we obtain the following estimate:
	\begin{align}
		\big\| (k_1)' \big\|_{L^{\infty}(0,s_1)} & \leqslant \frac{1-\mu}{2} \, 2 \, (p_0-p_1)^{-3} \, \big(2 (p_0-p_1) \big)^{2-\mu} \, \| \tilde{u} \|_{L^{\infty}(p_1,p_2)} \, (p_0 - p_1)^{-1} \nonumber \\
												& \qquad + \big(2 (p_0-p_1) \big)^{1-\mu} \, \| \tilde{u}' \|_{L^{\infty}(p_1,p_2)} \, (p_0 - p_1)^{-2} \nonumber \\												& \qquad \qquad + \big(2 (p_0-p_1) \big)^{1-\mu} \| \tilde{u} \|_{L^{\infty}(p_1,p_2)} \, 2 (p_0-p_1)^{-3} \nonumber \\
												& \label{derk1} \leqslant 2^{1-\mu} \, \| \tilde{u} \|_{W^{1,\infty}(p_1,p_2)} \, \Big( 2(2-\mu)(p_0-p_1)^{-(2+\mu)} + (p_0 - p_1)^{-(1+\mu)} \Big) \; .
	\end{align}
	
	As above, we study the function $k_2$ by employing its definition. We differentiate $k_2$ by using the fact that $U(p) = (p-p_1)^{\mu-1} \tilde{u}(p)$,
	\begin{equation*}
		(k_2)'(s) = \Big( (\mu-1) \big( \varphi_2^{-1}(s) - p_1\big)^{\mu-2} \, \tilde{u}\big( \varphi_2^{-1}(s) \big) + \big( \varphi_2^{-1}(s) - p_1 \big)^{\mu-1} \, \tilde{u}'\big( \varphi_2^{-1}(s) \big) \Big) (\varphi_2^{-1})'(s)^2 \; ,
	\end{equation*}
	for all $s \in [0,s_2]$. We employ the fact that $\varphi_2^{-1}(s) \in [q,p_2]$ for any $s \in [0,s_2]$ and the equality $(\varphi_2^{-1})'(s) = -1$ to obtain
	\begin{align}
		\big\| (k_2)' \big\|_{L^{\infty}(0,s_2)}	& \leqslant \Big((1-\mu) \, 2^{2-\mu} (p_0-p_1)^{\mu-2} \, \| \tilde{u} \|_{L^{\infty}(p_1,p_2)} + 2^{1-\mu} (p_0-p_1)^{\mu-1} \| \tilde{u}' \|_{L^{\infty}(p_1,p_2)} \Big) \nonumber \\
													& \label{k2} \leqslant 2^{1-\mu} \, \| \tilde{u} \|_{W^{1,\infty}(p_1,p_2)} \Big( 2(1-\mu) (p_0-p_1)^{\mu-2} + (p_0-p_1)^{\mu-1} \Big) \; .
	\end{align}
	These considerations permit to estimate the four remainders.
	\begin{itemize}
		\item \emph{Estimate of $R_1^{(1)}(\omega,p_0)$.} Theorem \ref{THM2} furnishes an estimate of $R^{(1)}(\omega,p_0)$. We combine it with the estimates of $(k_1)'$ given in \eqref{derk1} and $s_1$ given in \eqref{s1}:
		\begin{align}
				\left| R_1^{(1)}(\omega,p_0) \right|	& \leqslant \int_0^{s_1} s^{\mu-1} \big|(k_1)'(s) \big| \, ds \; \omega^{-1} \nonumber \\
														& \leqslant \frac{1}{\mu} \, s_1^{\; \mu} \, \big\| (k_1)' \big\|_{L^{\infty}(0,s_1)} \, \omega^{-1} \nonumber \\
														& \leqslant \frac{2^{1-\mu}}{\mu} \, \| \tilde{u} \|_{W^{1,\infty}(p_1,p_2)} \Big( 2(2-\mu)(p_0-p_1)^{\mu-2} + (p_0 - p_1)^{\mu-1} \Big) \, \omega^{-1} \nonumber \\
														& \label{R111} =: R_1^{(1)}(U) \, (p_0 - p_1)^{-\alpha_1^{(1)}} \, \omega^{-\beta_1^{(1)}} + R_2^{(1)}(U) \, (p_0-p_1)^{-\alpha_2^{(1)}} \, \omega^{-\beta_2^{(1)}} \; ,
		\end{align}
		where
		\begin{align*}
			& \bullet \hspace{2mm} R_1^{(1)}(U) := \frac{2^{2-\mu}}{\mu} \, (2-\mu) \, \| \tilde{u} \|_{W^{1,\infty}(p_1,p_2)} \quad , \quad R_2^{(1)}(U) := \frac{2^{1-\mu}}{\mu} \, \| \tilde{u} \|_{W^{1,\infty}(p_1,p_2)} \; , \\
			& \bullet \hspace{2mm} \alpha_1^{(1)} := 2 - \mu \quad , \quad \alpha_2^{(1)} := 1 - \mu \quad , \quad \beta_1^{(1)} = \beta_2^{(1)} := 1 \; .
		\end{align*}
		\item \emph{Estimate of $R_2^{(1)}(\omega,p_0)$.} The estimate of $R_2^{(1)}(\omega,p_0)$ from Theorem \ref{THM2} provides
		\begin{align}
			\left| R_2^{(1)}(\omega,p_0) \right|	& \leqslant (1-\mu) \, \left| U(q) \, (\varphi_1)'(q)^{-1} \right| \, \varphi_1(q)^{-1} \, \omega^{-2} \nonumber \\
													& \leqslant \frac{1-\mu}{2^{\mu-1}} \, \| \tilde{u} \|_{L^{\infty}(p_1,p_2)} (p_0 - p_1)^{\mu-1} \, (p_0-p_1)^{-1} \, \left( \frac{3}{4} \, (p_0-p_1)^2 \right)^{-1} \, \omega^{-2} \nonumber \\
													& = (1-\mu) \, \frac{2^{3-\mu}}{3} \, \| \tilde{u} \|_{L^{\infty}(p_1,p_2)} \, (p_0 - p_1)^{\mu-4} \, \omega^{-2} \nonumber \\
													& \label{R211} =: R_3^{(1)}(U) \,  (p_0 - p_1)^{-\alpha_3^{(1)}} \, \omega^{-\beta_3^{(1)}} \; ,
		\end{align}
		where the definition of $U$, inequality \eqref{24} and the value of $s_1$ given in \eqref{s1} were used. We have defined
		\begin{align*}
			& \bullet \hspace{2mm} R_3^{(1)}(U) := (1-\mu) \, \frac{2^{3-\mu}}{3} \, \| \tilde{u} \|_{L^{\infty}(p_1,p_2)} \; , \\
			& \bullet \hspace{2mm} \alpha_3^{(1)} := 4 - \mu \quad , \quad \beta_3^{(1)} := 2 \; .
		\end{align*}		
		\item \emph{Estimate of $R_1^{(2)}(\omega,p_0)$.} Here $\mu_2 = 1$, so we have to employ the estimate of $R_1^{(2)}(\omega,p_0)$ provided by Theorem \ref{THM3},
		\begin{equation} \label{R121}
			\begin{aligned}
				\left| R_1^{(2)}(\omega,p_0) \right|	& \leqslant L_{\gamma, 2} \int_0^{s_2} s^{-\gamma} \big| (k_2)'(s) \big| \, ds \; \omega^{-\delta} \\
														& \leqslant \frac{L_{\gamma, 2}}{1-\gamma} \, s_2^{\; 1-\gamma} \big\| (k_2)' \big\|_{L^{\infty}(0,s_2)} \, \omega^{-\delta} \\
														& \leqslant \frac{L_{\gamma, 2}}{1-\gamma} \, 2^{1-\mu} \, \| \tilde{u} \|_{W^{1,\infty}(p_1,p_2)} \Big( 2(1-\mu) (p_0-p_1)^{\mu-1-\gamma} + (p_0-p_1)^{\mu-\gamma} \Big) \, \omega^{-\delta} \\
														& =: R_4^{(1)}(U) \, (p_0 - p_1)^{-\alpha_4^{(1)}} \, \omega^{-\beta_4^{(1)}} + R_5^{(1)}(U) \, (p_0 - p_1)^{-\alpha_5^{(1)}} \, \omega^{-\beta_5^{(1)}} \; ,
			\end{aligned}
		\end{equation}
		where the last inequality was obtained by using \eqref{s2} and \eqref{k2}. Here the parameter $\delta$ is arbitrarily chosen in $\big( \frac{1}{2} , 1 \big)$, and the parameter $\gamma \in (0,1)$ is given by $\gamma = 2 \delta - 1$. We have defined
		\begin{align*}
			& \bullet \hspace{2mm} R_4^{(1)}(U) := \frac{L_{\gamma, 2}}{1-\gamma} \, 2^{2-\mu} \, (1-\mu) \, \| \tilde{u} \|_{W^{1,\infty}(p_1,p_2)} \quad , \quad R_5^{(1)}(U) := \frac{L_{\gamma, 2}}{1-\gamma} \, 2^{1-\mu} \, \| \tilde{u} \|_{W^{1,\infty}(p_1,p_2)} \; , \\
			& \bullet \hspace{2mm} \alpha_4^{(1)} := - \mu + 1 + \gamma \quad , \quad \alpha_5^{(1)} := \gamma - \mu \quad , \quad \beta_4^{(1)} = \beta_5^{(1)} := \delta \; .
		\end{align*}
		\item \emph{Estimate of $R_2^{(2)}(\omega,p_0)$.} We employ Theorem \ref{THM2} once again to control $R_2^{(2)}(\omega,p_0)$. The definition of $U$, the relation $(\varphi_2^{-1})' = -1$ and the value of $s_2$ given in \eqref{s2} lead to
		\begin{align}
			\left| R_2^{(2)}(\omega,p_0) \right|	& \leqslant \frac{1}{2} \, \Gamma \bigg(\frac{1}{2}\bigg) \, \left| U(q) \, (\varphi_2)'(q)^{-1} \right| \, \varphi_2(q)^{-2} \, \omega^{-\frac{3}{2}} \nonumber \\
												& \leqslant \frac{\sqrt{\pi}}{2^{\mu}} \, \| \tilde{u} \|_{L^{\infty}(p_1,p_2)} \, (p_0 - p_1)^{\mu-1} \, \big(2^{-1} (p_0-p_1) \big)^{-2} \, \omega^{-\frac{3}{2}} \nonumber \\
												& = \frac{\sqrt{\pi}}{2^{\mu-2}} \, \| \tilde{u} \|_{L^{\infty}(p_1,p_2)} (p_0-p_1)^{\mu-3} \, \omega^{-\frac{3}{2}} \nonumber \\
												& \label{R221} =: R_6^{(1)}(U) \, (p_0 - p_1)^{-\alpha_6^{(1)}} \, \omega^{-\beta_6^{(1)}} \; ,
		\end{align}
		where
		\begin{align*}
			& \bullet \hspace{2mm} R_6^{(1)}(U) := \frac{\sqrt{\pi}}{2^{\mu-2}} \, \| \tilde{u} \|_{L^{\infty}(p_1,p_2)} \; , \\
			& \bullet \hspace{2mm} \alpha_6^{(1)} := 3 - \mu \quad , \quad \beta_6^{(1)} := \frac{3}{2} \; .
		\end{align*}
	\end{itemize}
	
	\noindent \textit{Study of $I^{(2)}(\omega,p_0)$}. Firstly we remark that $\psi'$ is negative for all $p \in [p_0,p_2]$. To apply Theorem \ref{THM2}, we make the change of variables $p \mapsto -p$ in order to have an increasing phase. We obtain
	\begin{equation*}
		I^{(2)}(\omega,p_0) = \int_{p_0}^{p_2} U(p) e^{i \omega \psi(p)} \, dp = \int_{\check{p}_2}^{\check{p}_0} \check{U}(p) e^{i \omega \breve{\psi}(p)} \, dp \; ,
	\end{equation*}
	where we put $\check{U}(p) := U(-p)$, $\check{\psi}(p) := \psi(-p)$, $\check{p}_0 := -p_0$ and $\check{p}_2 := -p_2$.\\
	Thanks to this substitution, $\check{\psi}$ is now an increasing function that satisfies Assumption \emph{(P$_{1,2,N}$)} on $[\check{p}_2, \check{p}_0]$, for all $N \geqslant 1$, and by hypothesis $\check{U}$ verifies \emph{(A$_{1,1,1}$)} on $[\check{p}_2, \check{p}_0]$. Furthermore we remark that $\check{p_2}$ is not a singularity of the amplitude and not a stationary point. This observation indicates the non-necessity of a cutting-point. Hence, in the notations of Theorem \ref{THM2}, we employ only the expansion of the integral $\tilde{I}^{(2)}(\omega,p_0)$ with $p_2 = \check{p}_0$ and $q = \check{p}_2$. So we obtain from \eqref{I_2},
	\begin{equation*}
		\begin{aligned}
			I^{(2)}(\omega,p_0)	& = \phi^{(2)}(0,\omega,2,1) \, k_2(0) \, e^{i \omega \check{\psi}(\check{p}_0)} - \phi^{(2)}(s_2,\omega,2,1) \, k_2(s_2) \, e^{i \omega \check{\psi}(\check{p}_0)} \\
								& \qquad \qquad + e^{i \omega \check{\psi}(\check{p}_0)} \int_0^{s_2} \phi^{(2)}(s,\omega,2,1) (k_2)'(s) \, ds \; .
		\end{aligned}
	\end{equation*}
	Let us clarify the first terms by studying the function $\varphi_2$. The definition of $\varphi_2$ and the expression of $\psi$ yield
	\begin{equation*}
		\varphi_2(p) = \big( \check{\psi}(\check{p}_0) - \check{\psi}(p) \big)^{\frac{1}{2}} = \check{p}_0 - p \; ,
	\end{equation*}
	for all $p \in [\check{p}_2, \check{p}_0]$. It follows that $(\varphi_2)'(p) = (\varphi_2^{-1})'(s) = - 1$, and by the definition of $k_2$, we obtain
	\begin{equation*}
		k_2(s) = \check{U}\big(\varphi_2^{-1}(s) \big) (\varphi_2^{-1})'(s) = -\check{U}\big(\varphi_2^{-1}(s) \big) \; .
	\end{equation*}
	Since $\varphi_2(\check{p_0}) = 0$ and $\varphi_2(\check{p_2}) = s_2$, we have $k_2(0) = -\check{U}(\check{p_0}) = - U(p_0)$ and $k_2(s_2) = -U(p_2) = 0$, by the hypothesis on $U$.\\ Combining this with the expression of $\Theta^{(2)}(2,1)$ coming from Theorem \ref{THM2}, we obtain
	\begin{align*}
			I^{(2)}(\omega,p_0)	& = \frac{\sqrt{\pi}}{2} \, e^{-i \frac{\pi}{4}} \, e^{i \omega \psi(p_0)} \, \tilde{u}(p_0) \, (p_0 - p_1)^{\mu-1} \, \omega^{-\frac{1}{2}} \\
								& \qquad \qquad + e^{i \omega \psi(p_0)} \int_0^{s_2} \phi^{(2)}(s,\omega,2,1) (k_2)'(s) \, ds \; .
	\end{align*}
	As in the preceding step, let us estimate the remainder term. Firstly, we bound the number $s_2$ as follows:
	\begin{equation*}
		s_2 = \check{p}_0 - \check{p}_2 = p_2 - p_0 \leqslant p_2 - p_1 \; .
	\end{equation*}
	Now we establish an estimate of the first derivative of $k_2$. By the definition of this function and by the fact that $(\varphi_2^{-1})'=-1$, we have
	\begin{equation*}
		(k_2)'(s) = (1-\mu) \left(\check{p}_1 - \varphi_2^{-1}(s) \right)^{\mu-2} \check{\tilde{u}}\left(\varphi_2^{-1}(s) \right) + \left(\check{p}_1 - \varphi_2^{-1}(s) \right)^{\mu-1} \Big(\check{\tilde{u}}\Big)' \left(\varphi_2^{-1}(s)\right) \; ,
	\end{equation*}
	where $\check{p}_1 := -p_1$. Since $\varphi_2^{-1}(s) \in [\check{p}_2,\check{p}_0]$, it follows
	\begin{equation*}
		\big\| (k_2)' \big\|_{L^{\infty}(0,s_2)} \leqslant \left( (1-\mu) (p_0 - p_1)^{\mu-2} + (p_0 - p_1)^{\mu-1} \right) \left\| \tilde{u} \right\|_{W^{1,\infty}(p_1,p_2)} \; .
	\end{equation*}
	Combine this inequality with the estimate of the remainder given in Theorem \ref{THM2} to obtain
	\begin{align}
		\left| 	R_1^{(2)}(\omega,p_0) \right|	& \leqslant \frac{L_{\gamma,2}}{1-\gamma} \, (p_2-p_1)^{1-\gamma} \left\| \tilde{u} \right\|_{W^{1,\infty}(p_1,p_2)} \nonumber \\
												& \qquad \qquad \times \left( (1-\mu)(p_0-p_1)^{\mu-2} + (p_0-p_1)^{\mu-1} \right) \omega^{-\delta} \nonumber \\
												& \label{Re12} =: R_1^{(2)}(U) \, (p_0 - p_1)^{-\alpha_1^{(2)}} \, \omega^{-\beta_1^{(2)}} + R_2^{(2)}(U) \, (p_0 - p_1)^{-\alpha_2^{(2)}} \, \omega^{-\beta_2^{(2)}} \; ,
	\end{align}
	where $\gamma$ and $\delta$ are defined above. We have defined
	\begin{align*}
			& \bullet \hspace{2mm} R_1^{(2)}(U) := \frac{L_{\gamma,2}}{1-\gamma} \, (1-\mu) \, (p_2-p_1)^{1-\gamma} \, \left\| \tilde{u} \right\|_{W^{1,\infty}(p_1,p_2)} \quad , \\
			& \hspace{1cm} \quad R_2^{(2)}(U) := \frac{L_{\gamma,2}}{1-\gamma} \, (p_2-p_1)^{1-\gamma} \left\| \tilde{u} \right\|_{W^{1,\infty}(p_1,p_2)} \; , \\
			& \bullet \hspace{2mm} \alpha_1^{(2)} := 2 - \mu \quad , \quad \alpha_2^{(2)} := 1 - \mu \quad , \quad \beta_1^{(2)} = \beta_2^{(2)} := \delta \; .
		\end{align*}
\end{proof}

Employing the preceding result, we derive asymptotic expansions of the oscillatory integral with explicit error estimates. We distinguish three cases depending on the strenght of the singularity for readability.

\begin{COR2} \label{COR2}
	Under the assumptions of Theorem \ref{THM4}, let us define $\tilde{H}(\omega,\psi,U)$ and $\tilde{K}(\omega,\psi,U)$ as follows,
	\begin{equation*}
		\tilde{H}(\omega,\psi,U) := \sqrt{\pi} \, e^{-i \frac{\pi}{4}} \, e^{i \omega c} \, \tilde{u}(p_0) \qquad , \qquad \tilde{K}_{\mu}(\omega, \psi, U) := \frac{\Gamma(\mu)}{2^{\mu}} \, e^{i \frac{\pi \mu}{2}} \, e^{i \omega \psi(p_1)} \, \tilde{u}(p_1) \; .
	\end{equation*}
	Then we have
	\begin{itemize}
		\item \textit{Case $\displaystyle \mu > \frac{1}{2}$ :}
	\begin{equation*}
		\left| \int_{p_1}^{p_2} U(p) \, e^{i \omega \psi(p)} \, dp - \tilde{H}(\omega,\psi,U) \, (p_0-p_1)^{\mu-1} \, \omega^{-\frac{1}{2}} \right| \leqslant \sum_{k = 1}^9 \, \tilde{R}_k^{(1)}(U) \, (p_0-p_1)^{-\tilde{\alpha}_k^{(1)}} \, \omega^{-\tilde{\beta}_k^{(1)}} \; ,
	\end{equation*}
	where the constants $\tilde{R}_k^{(1)}(U) \geqslant 0$ and the exponents $\tilde{\alpha}_k^{(1)} \in \R$, $\tilde{\beta}_k^{(1)} > \frac{1}{2}$ are given in the proof ;
		\item \textit{Case $\displaystyle \mu = \frac{1}{2}$ :}
	\begin{align*}
		& \left| \int_{p_1}^{p_2} U(p) \, e^{i \omega \psi(p)} \, dp - \Big( \tilde{H}(\omega,\psi,U) + \tilde{K}_{\mu}(\omega,\psi,U) \Big) \, (p_0-p_1)^{-\frac{1}{2}} \, \omega^{-\frac{1}{2}} \right| \\
		& \hspace{4cm} \leqslant \sum_{k = 1}^8 \, \tilde{R}_k^{(2)}(U)  \, (p_0-p_1)^{-\tilde{\alpha}_k^{(2)}} \, \omega^{-\tilde{\beta}_k^{(2)}} \; ,
	\end{align*}
	where the constants $\tilde{R}_k^{(2)}(U) \geqslant 0$ and the exponents $\tilde{\alpha}_k^{(2)} \in \R$, $\tilde{\beta}_k^{(2)} > \frac{1}{2}$ are given in the proof ;
		\item \textit{Case $\displaystyle \mu < \frac{1}{2}$ :}
	\begin{equation*}
		\left| \int_{p_1}^{p_2} U(p) \, e^{i \omega \psi(p)} \, dp - \tilde{K}_{\mu}(\omega,\psi,U) \, (p_0-p_1)^{-\mu} \, \omega^{-\mu} \right| \leqslant \sum_{k = 1}^9 \, \tilde{R}_k^{(3)}(U) \, (p_0-p_1)^{-\tilde{\alpha}_k^{(3)}} \, \omega^{-\tilde{\beta}_k^{(3)}} \; ,
	\end{equation*}
	where the constants $\tilde{R}_k^{(3)}(U) \geqslant 0$ and the exponents $\tilde{\alpha}_k^{(3)} \in \R$, $\tilde{\beta}_k^{(3)} > \mu$ are given in the proof.
	\end{itemize}
\end{COR2}

\begin{proof}
	This result is a direct consequence of Theorem \ref{THM4}. We start by splitting the integral as follows,
	\begin{equation*}
		\int_{p_1}^{p_2} U(p) \, e^{i \omega \psi(p)} \, dp = \int_{p_1}^{p_0} \dots \; + \int_{p_0}^{p_2} \dots \; =: I^{(1)}(\omega,p_0) + I^{(2)}(\omega,p_0) \; .
	\end{equation*}
	Then we apply Theorem \ref{THM4} by distinguishing the three following cases :
	\begin{itemize}
		\item \textit{Case $\displaystyle \mu > \frac{1}{2}$ :}
	\begin{align*}
		& \left| \int_{p_1}^{p_2} U(p) \, e^{i \omega \psi(p)} \, dp - \tilde{H}(\omega, \psi, U) \, (p_0 - p_1)^{\mu-1} \, \omega^{-\frac{1}{2}} \right| \\
		& \hspace{1cm} = \bigg| \int_{p_1}^{p_0} U(p) \, e^{i \omega \psi(p)} \, dp -  \frac{1}{2} \, \tilde{H}(\omega,\psi,U) \, (p_0 - p_1)^{\mu-1} \, \omega^{-\frac{1}{2}}  \\
		& \hspace{2cm} + \int_{p_0}^{p_2} U(p) \, e^{i \omega \psi(p)} \, dp -  \frac{1}{2} \, \tilde{H}(\omega,\psi,U) \, (p_0 - p_1)^{\mu-1} \, \omega^{-\frac{1}{2}} \bigg| \\
		& \hspace{1cm} \leqslant \left| \tilde{K}_{\mu}(\omega,\psi,U) \, (p_0-p_1)^{-\mu} \, \omega^{-\mu} \right| + \sum_{k = 1}^6 \, R_k^{(1)}(U) \, (p_0 - p_1)^{-\alpha_k^{(1)}} \, \omega^{-\beta_k^{(1)}} \\
		& \hspace{2cm} + \sum_{k = 1}^2 \, R_k^{(2)}(U) \, (p_0 - p_1)^{-\alpha_k^{(2)}} \, \omega^{-\beta_k^{(2)}} \\
		& \hspace{1cm} \leqslant \frac{\Gamma(\mu)}{2^{\mu}} \, \left\| \tilde{u} \right\|_{L^{\infty}(p_1,p_2)} \, (p_0-p_1)^{-\mu} \, \omega^{-\mu} + \sum_{k = 1}^6 \, R_k^{(1)}(U) \, (p_0 - p_1)^{-\alpha_k^{(1)}} \, \omega^{-\beta_k^{(1)}} \\
		& \hspace{2cm} + \sum_{k = 1}^2 \, R_k^{(2)}(U) \, (p_0 - p_1)^{-\alpha_k^{(2)}} \, \omega^{-\beta_k^{(2)}} \\
		& \hspace{1cm} =: \sum_{k = 1}^9 \, \tilde{R}_k^{(1)}(U) \, (p_0 - p_1)^{-\tilde{\alpha}_k^{(1)}} \, \omega^{-\tilde{\beta}_k^{(1)}} \; ,
	\end{align*}
	where
	\begin{align*}
		& \bullet \hspace{2mm} \tilde{R}_1^{(1)}(U) := \frac{\Gamma(\mu)}{2^{\mu}} \, \left\| \tilde{u} \right\|_{L^{\infty}(p_1,p_2)} \quad , \quad \tilde{\alpha}_1^{(1)} := \mu \quad , \quad \tilde{\beta}_1^{(1)} := \mu \; ; \\[3mm]
		& \bullet \hspace{2mm} \tilde{R}_{k+1}^{(1)}(U) := R_k^{(1)}(U) \quad , \quad \tilde{\alpha}_{k+1}^{(1)} := \alpha_k^{(1)} \quad , \quad \tilde{\beta}_{k+1}^{(1)} := \beta_k^{(1)} \qquad k = 1, \dots, 6 \; ;\\[3mm]
		& \bullet \hspace{2mm} \tilde{R}_{k+7}^{(1)}(U) := R_k^{(2)}(U) \quad , \quad \tilde{\alpha}_{k+7}^{(1)} := \alpha_k^{(2)} \quad , \quad \tilde{\beta}_{k+7}^{(1)} := \beta_k^{(2)} \qquad k = 1, 2 \; .
	\end{align*}
	One can check that each $\tilde{\beta}_k^{(1)}$ is strictly larger than $\frac{1}{2}$, ensuring that the decay rate of each remainder term is higher than $\omega^{-\frac{1}{2}}$.
		\item \textit{Case $\displaystyle \mu = \frac{1}{2}$ :}
	\begin{align*}
		& \left| \int_{p_1}^{p_2} U(p) \, e^{i \omega \psi(p)} \, dp - \Big( \tilde{H}(\omega,\psi,U) + \tilde{K}_{\frac{1}{2}}(\omega, \psi,U) \Big) \, (p_0 - p_1)^{\mu-1} \, \omega^{-\frac{1}{2}} \right|	\\
		& \hspace{1cm} = \bigg| \int_{p_1}^{p_0} U(p) \, e^{i \omega \psi(p)} \, dp - \left( \frac{1}{2} \, \tilde{H}(\omega,\psi,U) + \tilde{K}_{\frac{1}{2}}(\omega, \psi,U) \right) \, (p_0 - p_1)^{\mu-1} \, \omega^{-\frac{1}{2}} \\
		& \hspace{2cm} + \int_{p_0}^{p_2} U(p) \, e^{i \omega \psi(p)} \, dp - \frac{1}{2} \, \tilde{H}(\omega,\psi,U) \, (p_0 - p_1)^{\mu-1} \, \omega^{-\frac{1}{2}} \bigg|	\\
		& \hspace{1cm} \leqslant \sum_{k = 1}^6 \, R_k^{(1)}(U) \, (p_0 - p_1)^{-\alpha_k^{(1)}} \, \omega^{-\beta_k^{(1)}} + \sum_{k = 1}^2 \, R_k^{(2)}(U) \, (p_0 - p_1)^{-\alpha_k^{(2)}} \, \omega^{-\beta_k^{(2)}} \\
		& \hspace{1cm} =: \sum_{k = 1}^8 \, \tilde{R}_k^{(2)}(U) \, (p_0 - p_1)^{-\tilde{\alpha}_k^{(2)}} \, \omega^{-\tilde{\beta}_k^{(2)}} \; ,
	\end{align*}
	where
	\begin{align*}
		& \bullet \hspace{2mm} \tilde{R}_{k}^{(2)}(U) := R_k^{(1)}(U) \quad , \quad \tilde{\alpha}_{k}^{(2)} := \alpha_k^{(1)} \quad , \quad \tilde{\beta}_{k}^{(2)} := \beta_k^{(1)} \qquad k = 1, \dots, 6 \; ;\\[3mm]
		& \bullet \hspace{2mm} \tilde{R}_{k+6}^{(2)}(U) := R_k^{(2)}(U) \quad , \quad \tilde{\alpha}_{k+6}^{(2)} := \alpha_k^{(2)} \quad , \quad \tilde{\beta}_{k+6}^{(2)} := \beta_k^{(2)} \qquad k = 1, 2 \; .
	\end{align*}
	As above, we have $\tilde{\beta}_k^{(2)} > \frac{1}{2}$.
		\item \textit{Case $\displaystyle \mu < \frac{1}{2}$ :}
	\begin{align*}
		& \left| \int_{p_1}^{p_2} U(p) \, e^{i \omega \psi(p)} \, dp - \tilde{K}_{\mu}(\omega, \psi ,U) \, (p_0-p_1)^{-\mu} \, \omega^{-\mu} \right|	\\
		& \hspace{1cm} = \left| \int_{p_1}^{p_0} U(p) \, e^{i \omega \psi(p)} \, dp - \tilde{K}_{\mu}(\omega,\psi ,U) \, (p_0-p_1)^{-\mu} \, \omega^{-\mu} + \int_{p_0}^{p_2} U(p) \, e^{i \omega \psi(p)} \, dp \right| \\
		& \hspace{1cm} \leqslant \left| \frac{1}{2} \, \tilde{H}(\omega, \psi,U) (p_0 - p_1)^{\mu-1} \, \omega^{-\frac{1}{2}} \right| + \sum_{k = 1}^6 \, R_k^{(1)}(U) \, (p_0 - p_1)^{-\alpha_k^{(1)}} \, \omega^{-\beta_k^{(1)}}  \\
		& \hspace{2cm} + \left| \frac{1}{2} \, \tilde{H}(\omega, \psi ,U) \, (p_0 - p_1)^{\mu-1} \, \omega^{-\frac{1}{2}} \right| + \sum_{k = 1}^2 \, R_k^{(2)}(U) \, (p_0 - p_1)^{-\alpha_k^{(2)}} \, \omega^{-\beta_k^{(2)}} \\
		& \hspace{1cm} \leqslant \sqrt{\pi} \, \left\| \tilde{u} \right\|_{L^{\infty}(p_1,p_2)} \, (p_0-p_1)^{\mu - 1} \, \omega^{-\frac{1}{2}} \\
		& \hspace{2cm} + \sum_{k = 1}^6 \, R_k^{(1)}(U) \, (p_0 - p_1)^{-\alpha_k^{(1)}} \, \omega^{-\beta_k^{(1)}} + \sum_{k = 1}^2 \, R_k^{(2)}(U) \, (p_0 - p_1)^{-\alpha_k^{(2)}} \, \omega^{-\beta_k^{(2)}} \\
		& \hspace{1cm} =: \sum_{k = 1}^9 \, \tilde{R}_k^{(3)}(U) \, (p_0 - p_1)^{-\tilde{\alpha}_k^{(3)}} \, \omega^{-\tilde{\beta}_k^{(3)}} \; ,
	\end{align*}
	where
	\begin{align*}
		& \bullet \hspace{2mm} \tilde{R}_1^{(3)}(U) := \sqrt{\pi} \, \left\| \tilde{u} \right\|_{L^{\infty}(p_1,p_2)} \quad , \quad \tilde{\alpha}_1^{(3)} := 1-\mu \quad , \quad \tilde{\beta}_1^{(3)} := \frac{1}{2} \; ; \\[3mm]
		& \bullet \hspace{2mm} \tilde{R}_{k+1}^{(3)}(U) := R_k^{(1)}(U) \quad , \quad \tilde{\alpha}_{k+1}^{(3)} := \alpha_k^{(1)} \quad , \quad \tilde{\beta}_{k+1}^{(3)} := \beta_k^{(1)} \qquad k = 1, \dots, 6 \; ;\\[3mm]
		& \bullet \hspace{2mm} \tilde{R}_{k+7}^{(3)}(U) := R_k^{(2)}(U) \quad , \quad \tilde{\alpha}_{k+7}^{(3)} := \alpha_k^{(2)} \quad , \quad \tilde{\beta}_{k+7}^{(3)} := \beta_k^{(2)} \qquad k = 1, 2 \; .
	\end{align*}
	Here we note that $\tilde{\beta}_5^{(3)} = \tilde{\beta}_6^{(3)} = \delta$. So we can choose $\delta \in \big( \frac{1}{2} , 1 \big)$ so that $\delta > \mu$ and thanks to that, every $\tilde{\beta}_k^{(3)}$ are strictly larger than $\mu$.
\end{itemize}
\end{proof}

\begin{REM6} \label{REM6}
	\emph{Looking carefully to the value of each $\tilde{\alpha}_k^{(j)}$, we note that  only $\tilde{\alpha}_6^{(1)} = \tilde{\alpha}_5^{(2)} = \tilde{\alpha}_6^{(3)} = \alpha_5^{(1)} = \gamma - \mu$ can be negative. To prevent this situation, we can choose $\delta \in \big( \frac{1}{2} , 1 \big)$ so that $\gamma - \mu \geqslant 0$, namely $\delta \geqslant \frac{\mu + 1}{2}$. This will be useful to simplify slightly the proof in Theorems \ref{THM6} and \ref{THM9}, and Corollary \ref{COR3}.}
\end{REM6}

From the previous corollary, we observe that the blow-up of the asymptotic expansion comes from the terms $(p_0-p_1)^{-\tilde{\alpha}_k^{(j)}}$. This motivates the idea of considering $p_0$ approaching $p_1$ with a certain convergence speed, described by the parameter $\varepsilon > 0$, when the large parameter $\omega$ tends to infinity, as for example $p_0 - p_1 = \omega^{-\varepsilon}$. This procedure modifies the decay rates of the asymptotic expansion of the integral. We shall exploit the idea that below a certain threshold, the convergence speed is sufficiently slow so that the decay with respect to $\omega$ compensates the blow-up. Thus the proof consists in finding the values of $\varepsilon > 0$ for which the decay rate of the remainder is strictly larger than the one of the first term. This leads to asymptotic expansions on curves in the space of the parameters.

\begin{THM5} \label{THM5}
	Let $\varepsilon \in \big(0, \,  \frac{1}{2} \big)$ and suppose that $\displaystyle p_0 := p_1 + \omega^{-\varepsilon}$ for $\displaystyle \omega > (p_2 - p_1)^{-\frac{1}{\varepsilon}}$. Then under the assumptions of Theorem \ref{THM4} and with the notations of Corollary \ref{COR2}, we have
	\begin{itemize}
		\item \textit{Case $\displaystyle \mu > \frac{1}{2}$ :}
	\begin{equation*}
			\left| \int_{p_1}^{p_2} U(p) \, e^{i \omega \psi(p)} \, dp - \tilde{H}(\omega, \psi,U) \, \omega^{-\frac{1}{2} + \varepsilon(1-\mu)} \right| \leqslant \sum_{k = 1}^9 \, \tilde{R}_k^{(1)}(U) \, \omega^{-\tilde{\beta}_k^{(1)} + \varepsilon \, \tilde{\alpha}_k^{(1)}} \; ,
	\end{equation*}
	where $\displaystyle \max_{k \in \{1,...,9\}} \left\{ -\tilde{\beta}_k^{(1)} + \varepsilon \, \tilde{\alpha}_k^{(1)} \right\} < -\frac{1}{2} + \varepsilon(1-\mu)$ ;
		\item \textit{Case $\displaystyle \mu = \frac{1}{2}$ :}
	\begin{equation*}
		\left| \int_{p_1}^{p_2} U(p) \, e^{i \omega \psi(p)} \, dp - \Big( \tilde{H}(\omega,\psi,U) + \tilde{K}_{\mu}(\omega,\psi,U) \Big) \, \omega^{-\frac{1}{2} + \frac{1}{2} \, \varepsilon} \right| \leqslant \sum_{k = 1}^8 \, \tilde{R}_k^{(2)}(U)  \, \omega^{-\tilde{\beta}_k^{(2)} + \varepsilon \, \tilde{\alpha}_k^{(2)}} \; ,
	\end{equation*}
	where $\displaystyle \max_{k \in \{1,...,8\}} \left\{ -\tilde{\beta}_k^{(2)} + \varepsilon \, \tilde{\alpha}_k^{(2)} \right\} < -\frac{1}{2} + \frac{1}{2} \, \varepsilon$ ;
		\item \textit{Case $\displaystyle \mu < \frac{1}{2}$ :}
	\begin{equation*}
		\left| \int_{p_1}^{p_2} U(p) \, e^{i \omega \psi(p)} \, dp - \tilde{K}_{\mu}(\omega,\psi,U) \, \omega^{-\mu + \varepsilon \mu} \right| \leqslant \sum_{k = 1}^9 \, \tilde{R}_k^{(3)}(U) \, \omega^{-\tilde{\beta}_k^{(3)} + \varepsilon \, \tilde{\alpha}_k^{(3)}} \; ,
	\end{equation*}
	where $\displaystyle \max_{k \in \{1,...,9\}} \left\{ -\tilde{\beta}_k^{(3)} + \varepsilon \, \tilde{\alpha}_k^{(3)} \right\} < -\mu + \varepsilon \mu$ .
	\end{itemize}
\end{THM5}

\begin{proof}
	First of all, note that $I^{(1)}(\omega,p_0)$ and $I^{(2)}(\omega,p_0)$ are well-defined since the hypothesis $\omega > (p_2 - p_1)^{-\frac{1}{\varepsilon}}$ implies $p_0 \in (p_1,p_2)$. Now we replace $p_0 - p_1$ by $\omega^{-\varepsilon}$ in the estimates of Corollary \ref{COR2} and we compare the decay rates of the expansion with those of the remainder. In the following, we choose the parameter $\delta \in \big( \frac{1}{2} , 1 \big)$ in such way that we have $ \delta > \frac{1}{2} + \varepsilon$ and $\delta > \mu$.	
	\begin{itemize}
		\item \textit{Case $\displaystyle \mu > \frac{1}{2}$ :}	here we have
		\begin{equation*}
			\tilde{H}(\omega,\psi,U) \, (p_0 - p_1)^{\mu-1} \, \omega^{-\frac{1}{2}} = \tilde{H}(\omega,\psi,U) \, \omega^{-\frac{1}{2} + \varepsilon(1-\mu)} \; .
		\end{equation*}
		We note that the modulus of the coefficient $\tilde{H}(\omega,p_0,U)$ can be bounded from above and below by a non-zero constant when $\omega$ is sufficiently large, due to the hypothesis $\tilde{u}(p_1) \neq 0$. It follows that this coefficient does not influence the decay and so the expansion behaves like $\omega^{-\frac{1}{2} + \varepsilon(1-\mu)}$ when $\omega$ tends to infinity.\\
		To compare the decay rates of the expansion and of the remainder, it is sufficient to compare the exponents of $\omega$. Putting $\omega^{-\varepsilon} = p_0 - p_1$ in \eqref{R111}, \eqref{R211}, \eqref{R121}, \eqref{R221} and \eqref{Re12}, we obtain new decay rates for the remainder terms and the exponents have to be less than $\displaystyle -\frac{1}{2} + \varepsilon(1-\mu)$. This leads to the following system of inequalities :
		\begin{subnumcases}{}
		\displaystyle -\frac{1}{2} + \varepsilon(1-\mu) > -\mu + \varepsilon \mu \label{systineq1}\\
		\displaystyle -\frac{1}{2} + \varepsilon(1-\mu) > -1 + \varepsilon(2-\mu) \label{systineq2} \\
		\displaystyle -\frac{1}{2} + \varepsilon(1-\mu) > -1 + \varepsilon(1-\mu) \label{systineq3} \\
		\displaystyle -\frac{1}{2} + \varepsilon(1-\mu) > -2 + \varepsilon(4-\mu) \label{systineq4} \\
		\displaystyle -\frac{1}{2} + \varepsilon(1-\mu) > -\delta + \varepsilon(1+\gamma-\mu) \label{systineq5} \\
		\displaystyle -\frac{1}{2} + \varepsilon(1-\mu) > -\delta + \varepsilon(\gamma - \mu) \label{systineq6} \\
		\displaystyle -\frac{1}{2} + \varepsilon(1-\mu) > - \frac{3}{2} + \varepsilon(3-\mu) \label{systineq7}\\
		\displaystyle -\frac{1}{2} + \varepsilon(1-\mu) > -\delta + \varepsilon(1-\mu) \label{systineq8} \\
		\displaystyle -\frac{1}{2} + \varepsilon(1-\mu) > -\delta + \varepsilon(2-\mu) \label{systineq9}
	\end{subnumcases}
	\begin{itemize}
		\item \emph{Inequalities \eqref{systineq1}, \eqref{systineq2}, \eqref{systineq4}, \eqref{systineq5}, \eqref{systineq7}.} These inequalities are satisfied because they are equivalent to $\varepsilon < \frac{1}{2}$, which is true by hypothesis. Note that we used the relation $\gamma = 2 \delta - 1$ to study \eqref{systineq5}.
		\item \emph{Inequality \eqref{systineq3}.} This inequality is equivalent to $\frac{1}{2} < 1$ which is clearly true.
		\item \emph{Inequality \eqref{systineq6}.} This inequality is equivalent to $ \delta - \frac{1}{2} > \varepsilon(\gamma-1)$. By hypothesis, we have $\delta > \frac{1}{2}$, $\varepsilon > 0$ and $\gamma < 1$. So \eqref{systineq6} is satisfied.
		\item \emph{Inequality \eqref{systineq8}.} This inequality is true because $\delta > \frac{1}{2}$ by hypothesis.
		\item \emph{Inequality \eqref{systineq9}.} We have supposed that $\varepsilon < \delta - \frac{1}{2}$, which is equivalent to \eqref{systineq9}. The latter is then satisfied.
	\end{itemize}
	\item \textit{Case $\displaystyle \mu = \frac{1}{2}$ :}	in this situation,
		\begin{equation*}
			\Big( \tilde{H}(\omega,\psi,U) + \tilde{K}_{\frac{1}{2}}(\omega, \psi, U) \Big) \, (p_0 - p_1)^{-\frac{1}{2}} \, \omega^{-\frac{1}{2}} = \Big( \tilde{H}(\omega,p_0,U) + \tilde{K}_{\frac{1}{2}}(\omega,p_0,U) \Big) \, \omega^{-\frac{1}{2} + \frac{1}{2} \varepsilon} \; .
		\end{equation*}
		As explained above, the expansion behaves like $\omega^{-\frac{1}{2} + \frac{1}{2} \varepsilon}$ when $\omega$ tends to infinity.\\
		Here the decay rates of the remainder have to be faster than $\omega^{-\frac{1}{2} + \frac{1}{2} \varepsilon}$. So we have to check the new system of inequalities,
		\begin{subnumcases}{}
		\displaystyle -\frac{1}{2} + \frac{1}{2} \varepsilon > -1 + \frac{3}{2} \varepsilon \label{systineq10} \\
		\displaystyle -\frac{1}{2} + \frac{1}{2} \varepsilon > -1 + \frac{1}{2} \varepsilon \label{systineq11} \\
		\displaystyle -\frac{1}{2} + \frac{1}{2} \varepsilon > -2 + \frac{7}{2} \varepsilon \label{systineq12} \\
		\displaystyle -\frac{1}{2} + \frac{1}{2} \varepsilon > -\delta + \left(\frac{1}{2}+\gamma \right) \varepsilon \label{systineq13} \\
		\displaystyle -\frac{1}{2} + \frac{1}{2} \varepsilon > -\delta + \left(\gamma - \frac{1}{2} \right) \varepsilon \label{systineq14} \\
		\displaystyle -\frac{1}{2} + \frac{1}{2} \varepsilon > - \frac{3}{2} + \frac{5}{2} \varepsilon \label{systineq15}\\
		\displaystyle -\frac{1}{2} + \frac{1}{2} \varepsilon > -\delta + \frac{1}{2} \varepsilon \label{systineq16} \\
		\displaystyle -\frac{1}{2} + \frac{1}{2} \varepsilon > -\delta + \frac{3}{2} \varepsilon \label{systineq17}
	\end{subnumcases}
	\begin{itemize}
		\item \emph{Inequalities \eqref{systineq10}, \eqref{systineq12}, \eqref{systineq13}, \eqref{systineq15}.} These inequalities are satisfied because they are equivalent to $\varepsilon < \frac{1}{2}$. The relation $\gamma = 2 \delta - 1$ is used to study \eqref{systineq13}.
		\item \emph{Inequality \eqref{systineq11}.} This inequality is equivalent to $\frac{1}{2} < 1$ which is clearly true.
		\item \emph{Inequality \eqref{systineq14}.} Similar to inequality \eqref{systineq6}.
		\item \emph{Inequality \eqref{systineq16}.} True because $\delta > \frac{1}{2}$ by hypothesis.
		\item \emph{Inequality \eqref{systineq17}.} Equivalent to $\varepsilon < \delta - \frac{1}{2}$, which is true.
	\end{itemize}
	\item \textit{Case $\displaystyle \mu < \frac{1}{2}$ :}	following the steps of the two preceding cases, we obtain
		\begin{equation*}
			\tilde{K}_{\mu}(\omega,\psi,U) \, (p_0 - p_1)^{-\mu} \, \omega^{-\mu}	 = K_{\mu}(\omega,\psi,U) \, \omega^{-\mu + \varepsilon \mu} \; ,
		\end{equation*}
		and in this situation, the system of inequalities is given by
		\begin{subnumcases}{}
		\displaystyle -\mu + \varepsilon \mu > -\frac{1}{2} + \varepsilon(1-\mu) \label{systineq18}\\
		\displaystyle -\mu + \varepsilon \mu > -1 + \varepsilon(2-\mu) \label{systineq19} \\
		\displaystyle -\mu + \varepsilon \mu > -1 + \varepsilon(1-\mu) \label{systineq20} \\
		\displaystyle -\mu + \varepsilon \mu > -2 + \varepsilon(4-\mu) \label{systineq21} \\
		\displaystyle -\mu + \varepsilon \mu > -\delta + \varepsilon(1+\gamma-\mu) \label{systineq22} \\
		\displaystyle -\mu + \varepsilon \mu > -\delta + \varepsilon(\gamma - \mu) \label{systineq23} \\
		\displaystyle -\mu + \varepsilon \mu > - \frac{3}{2} + \varepsilon(3-\mu) \label{systineq24}\\
		\displaystyle -\mu + \varepsilon \mu > -\delta + \varepsilon(1-\mu) \label{systineq25} \\
		\displaystyle -\mu + \varepsilon \mu > -\delta + \varepsilon(2-\mu) \label{systineq26}
	\end{subnumcases}
	\begin{itemize}
		\item \emph{Inequalities \eqref{systineq18}, \eqref{systineq19}, \eqref{systineq21}, \eqref{systineq22}, \eqref{systineq24}.} Equivalent to $\varepsilon < \frac{1}{2}$, which is true by hypothesis. We used once again the relation $\gamma = 2 \delta - 1$ and the fact that $\delta > \mu$ to study \eqref{systineq22}.
		\item \emph{Inequality \eqref{systineq20}.} This inequality is equivalent to $\varepsilon < \frac{1-\mu}{1-2\mu}$. But we can show that $\frac{1}{2} < \frac{1-\mu}{1-2\mu}$ and we recall that $\varepsilon < \frac{1}{2}$, so \eqref{systineq20} is verified.
		\item \emph{Inequality \eqref{systineq23}.} This inequality is equivalent to $\delta - \mu > \varepsilon ( \gamma - 2 \mu)$. Now we distinguish two cases: if $\gamma - 2 \mu \leqslant 0$ then the inequaliy holds since $\delta - \mu$ is  positive. In the other case, \eqref{systineq23} is equivalent to $\varepsilon < \frac{\delta - \mu}{\gamma - 2 \mu}$. Hence the latter is true since $\varepsilon < \frac{1}{2}$ and we can show $\frac{1}{2} < \frac{\delta - \mu}{\gamma - 2 \mu}$.
		\item \emph{Inequality \eqref{systineq25}.} It is equivalent to $\varepsilon < \frac{\delta - \mu}{1 - 2 \mu}$. Since $\frac{1}{2} < \frac{\delta - \mu}{1 - 2 \mu}$ holds and $\varepsilon < \frac{1}{2}$, \eqref{systineq25} is satisfied.
		\item \emph{Inequality \eqref{systineq26}.} It is equivalent to $\varepsilon < \frac{\delta - \mu}{2 - 2 \mu}$. We can show that $\delta - \frac{1}{2} < \frac{\delta - \mu}{2 - 2 \mu}$ and we recall that $\varepsilon < \delta - \frac{1}{2}$, which proves that \eqref{systineq26} is true.
	\end{itemize}
	\end{itemize}
	
\end{proof}

\section{Application to the free Schrödinger equation: propagation of wave packets and anomalous phenomena}

\hspace{2.5ex} In this section, we are interested in the time asymptotic behaviour of the solution of the free Schrödinger equation in one dimension, with initial conditions in a frequency band $[p_1,p_2]$ and we suppose that $p_1$ is a singular frequency. We establish time asymptotic expansions as well as uniform estimates to explore the influence of the compact frequency band and of the singularity on the dispersion.\\

We introduce the free Schrödinger equation on the line
\begin{equation*}
	(S) \qquad \left\{ \hspace{-3mm} \begin{array}{rl}
			& \big[ i \partial_t + \partial_x^2 \big] u(t) = 0 \\
			& \vspace{-0.3cm} \\
			& u(0) = u_0
	\end{array} \right. \; ,
\end{equation*}
for $t > 0$. If $u_0 \in \mathcal{S}'(\R)$ then this initial value problem has a unique solution in the tempered distribution given by
\begin{equation*}
		u(t) = \tf^{-1} \left( e^{-i t p^2} \tf u_0 \right) \; ,
	\end{equation*}	
where $\tf : \mathcal{S}'(\R) \longrightarrow \mathcal{S}'(\R)$ is the Fourier transform.\\
Throughout this section, we shall suppose that the initial data satisfies the following condition,\\

\noindent \textbf{Condition ($\mathbf{C_{[p_1,p_2],\mu}}$).} Let $\mu \in (0,1)$ and $p_1 < p_2$ be two finite real numbers.\\
	A tempered distribution $u_0$ satisfies Condition (C$_{[p_1,p_2],\mu}$) if and only if $\tf u_0 \equiv 0$ on $\R \setminus [p_1,p_2]$ and $\tf u_0$ verifies Assumption (A$_{\mu,1,1}$) on $[p_1,p_2]$, with $\tf u_0 (p_2) = 0$.\\
	
\noindent Under this condition, we note that $\tf u_0$ is a function which has a singularity of order $\mu-1$ at $p_1$ whereas the point $p_2$ is regular. For simplicity, we assume $\tf u_0 (p_2) = 0$ but a similar work can be carried out if $\tf u_0 (p_2) \neq 0$. Moreover in this situation, the solution formula of the free Schrödinger equation defines a smooth function $u : (0,+ \infty) \times \R \longrightarrow \C$ given by the following integral
\begin{equation} \label{solformula}
	u(t,x) = \frac{1}{2 \pi} \int_{p_1}^{p_2} \tf u_0 (p) \, e^{-it p^2 + ixp} \, dp \; .
\end{equation}
We recall that for $v \in L^1(\R)$, the Fourier transform of $v$ is defined by
\begin{equation*}
	\tf v(p) = \int_{\R} v(x) \, e^{-i xp} \, dx \; .
\end{equation*}

In this section, we shall need the following definitions of certain space-time regions.
	
\begin{DEF3}
	Let $a < b$ be two finite real numbers and let $\varepsilon > 0$.
	\begin{enumerate}
		\item We define the space-time cone $\mathfrak{C}(a, b)$ as follows :	\begin{equation*}
		\mathfrak{C}(a, b) := \left\{ (t,x) \in (0,+\infty) \times \R \: \bigg| \: 2 \, a < \frac{x}{t} < 2 \, b \right\} \; .
	\end{equation*}
		\item We define the space-time curve $\mathfrak{G}_{\varepsilon}(a)$ as follows :
	\begin{equation*}
		\mathfrak{G}_{\varepsilon}(a) := \left\{ (t,x) \in (0,+\infty) \times \R \: \bigg| \: \frac{x}{t} = 2 \, a + 2 \, t^{-\varepsilon} \right\} \; .
	\end{equation*}
		\item We define the space-time region $\mathfrak{R}_{\varepsilon}(a,b)$ as follows :
	\begin{equation*}
		\mathfrak{R}_{\varepsilon}(a,b) = \left\{ (t,x) \in (0,+\infty) \times \R \: \bigg| \: 2 \, a + 2 \, t^{-\varepsilon} \leqslant \frac{x}{t} < 2 \, b \; , \; t > T_{\varepsilon}(a,b) \right\} \; ,
	\end{equation*}
	where $T_{\varepsilon}(a,b) := (b-a)^{-\frac{1}{\varepsilon}}$.
	\end{enumerate}
\end{DEF3}

In the first result, we furnish uniform remainder estimates for asymptotic expansions of the solution in the cone $\mathfrak{C}\big( p_1 + \varepsilon, p_2 \big)$. After having rewritten the solution formula as an oscillatory integral, we can apply the results of the preceding section and we use the fact that the distance between the stationary point $\frac{x}{2t}$ and the singularity $p_1$ is bounded from below by $\varepsilon$ to estimate uniformly the remainder. Let us note that the method employed in the proof furnishes asymptotic expansions of the solution in the entire cone $\mathfrak{C}\big( p_1 , p_2 \big)$ with explicit blow-up when $\frac{x}{2t}$ approaches $p_1$. Especially, we see in the proof that restricting the cone $\mathfrak{C}\big( p_1 , p_2 \big)$ to $\mathfrak{C}\big( p_1 + \varepsilon, p_2 \big)$ is sufficient to obtain uniform estimates.\\
It is interesting to note that the singularity diminishes the time decay rate in the cone below the rate of quantum mechanic dispersion $t^{-\frac{1}{2}}$, when leaving the $L^2$-setting.

\begin{THM6} \label{THM6}
	Suppose that $u_0$ satisfies Condition \emph{(C$_{[p_1,p_2],\mu}$)} and choose a real number $\varepsilon > 0$ such that
	\begin{equation*}
		p_1 + \varepsilon < p_2 \; .
	\end{equation*}
	For all $(t,x) \in \mathfrak{C}\big( p_1 + \varepsilon, p_2 \big)$, define $H(t,x,u_0) \in \C$ and $K_{\mu}(t,x,u_0) \in \C$ as follows,
	\begin{align*}
		& \bullet \hspace{5mm} H(t,x,u_0) := \frac{1}{2 \sqrt{\pi}} \, e^{-i \frac{\pi}{4}} \, e^{i \frac{x^2}{4t}} \, \tilde{u} \hspace{-1mm} \left( \frac{x}{2t} \right) \left( \frac{x}{2t} - p_1 \right)^{\mu-1} \; , \\
		& \bullet \hspace{5mm} K_{\mu}(t,x,u_0) := \frac{\Gamma(\mu)}{2^{\mu+1}} \, e^{i \frac{\pi \mu}{2}} \, e^{-i t p_1^{\, 2} + i x p_1} \, \tilde{u}(p_1) \left( \frac{x}{2t} - p_1 \right)^{-\mu} \; .
	\end{align*}
	Then for all $(t,x) \in \mathfrak{C}\big( p_1 + \varepsilon, p_2 \big)$, we have
	\begin{itemize}
		\item \textit{Case $\displaystyle \mu > \frac{1}{2}$ :}
	\begin{equation*}
		\left| u(t,x) - H(t,x,u_0) \, t^{-\frac{1}{2}} \right| \leqslant \sum_{k = 1}^9 \, C_k^{(1)}(u_0, \varepsilon) \, t^{-\tilde{\beta}_k^{(1)}} \; ,
	\end{equation*}
	where $\displaystyle \max_{k \in \{1,...,9\}} \left\{ -\tilde{\beta}_k^{(1)} \right\} < - \frac{1}{2}$.  The exponents $\tilde{\beta}_k^{(1)}$ are defined in Corollary \ref{COR2} and the constants $C_k^{(1)}(u_0, \varepsilon) \geqslant 0$ are given in the proof ;
		\item \textit{Case $\displaystyle \mu = \frac{1}{2}$ :}
	\begin{equation*}
		\left| u(t,x) - \big( H(t,x,u_0) + K_{\mu}(t,x,u_0) \big) \, t^{- \frac{1}{2}} \right| \leqslant \sum_{k = 1}^8 \, C_k^{(2)}(u_0, \varepsilon) \, t^{-\tilde{\beta}_k^{(2)}} \; ,
	\end{equation*}
	where $\displaystyle \max_{k \in \{1,...,8\}} \left\{ -\tilde{\beta}_k^{(2)} \right\} < - \frac{1}{2}$.  The exponents $\tilde{\beta}_k^{(2)}$ are defined in Corollary \ref{COR2} and the constants $C_k^{(2)}(u_0, \varepsilon) \geqslant 0$ are given in the proof ;
		\item \textit{Case $\displaystyle \mu < \frac{1}{2}$ :}
	\begin{equation*}
		\left| u(t,x) - K_{\mu}(t,x,u_0) \, t^{-\mu} \right| \leqslant \sum_{k = 1}^9 \, C_k^{(3)}(u_0, \varepsilon) \, t^{-\tilde{\beta}_k^{(3)}} \; ,
	\end{equation*}
	where $\displaystyle \max_{k \in \{1,...,9\}} \left\{ -\tilde{\beta}_k^{(3)} \right\} < - \mu$.  The exponents $\tilde{\beta}_k^{(3)}$ are defined in Corollary \ref{COR2} and the constants $C_k^{(3)}(u_0, \varepsilon) \geqslant 0$are given in the proof.
	\end{itemize}
\end{THM6}

\begin{proof}
	We shall prove the result in the case $\mu > \frac{1}{2}$; the proofs in the other cases are very similar.\\
	In the solution formula \eqref{solformula}, we factorize the phase function $p \longmapsto -t p^2 + xp$ by $t$, which gives
	\begin{equation*}
		\forall \, (t,x) \in (0,+\infty) \times \R \qquad u(t,x) = \int_{p_1}^{p_2} U(p) \, e^{i t \psi(p)} \, dp \; ,
	\end{equation*}
	where \footnote[2]{See Remark \ref{REM0}}
	\begin{equation*}
		\left\{ \begin{array}{rl}
				& \displaystyle \forall \, p \in (p_1,p_2] \qquad U(p) := \frac{1}{2\pi}  \, \tf u_0(p) = \frac{1}{2\pi}  \, (p-p_1)^{\mu-1} \, \tilde{u}(p) \; , \\
				& \vspace{-0.3cm} \\
				& \displaystyle \forall \, p \in \R \qquad \psi(p) := -p^2 + \frac{x}{t} \, p \; .
		\end{array} \right.
	\end{equation*}
	By hypothesis, $U$ verifies Assumption (A$_{\mu,1,1}$) on $[p_1,p_2]$ and $\psi$ has the form
	\begin{equation*}
		\psi(p) = - \left( p - p_0 \right)^2 + c \; ,
	\end{equation*}
	where $p_0 := \frac{x}{2t}$ and $c := p_0^2 = \frac{x^2}{4 t^2}$. Moreover, we have the following equivalence,
	\begin{equation*}
		(t,x) \in \mathfrak{C}\big( p_1 + \varepsilon, p_2 \big) \qquad \Longleftrightarrow \qquad p_1 + \varepsilon < \frac{x}{2t} = p_0 < p_2 \; ,
	\end{equation*}
	implying the fact that the stationary point $p_0$ belongs to $(p_1,p_2)$. Hence Corollary \ref{COR2} is applicable and we obtain for all $(t,x) \in \mathfrak{C}\big( p_1 + \varepsilon, p_2 \big)$,
	\begin{equation*}
		\left| u(t,x) - H(t,x,u_0) \, t^{-\frac{1}{2}} \right| \leqslant \sum_{k = 1}^9 \, \tilde{R}_k^{(1)} \hspace{-1mm} \left( \frac{1}{2 \pi} \, \tf u_0 \right) \, \left( \frac{x}{2t} - p_1\right)^{-\tilde{\alpha}_k^{(1)}} \, t^{-\tilde{\beta}_k^{(1)}} \; ,
	\end{equation*}
	where the coefficient $H(t,x,u_0)$ is given in the statement of the theorem, the constants $\tilde{R}_k^{(1)}\big( \frac{1}{2 \pi} \, \tf u_0 \big) \geqslant 0$ and the exponents $\tilde{\alpha}_k^{(1)} \geqslant 0$, $\tilde{\beta}_k^{(1)} > \frac{1}{2}$ are provided by Corollary \ref{COR2}. Note that we can choose $\tilde{\alpha}_k^{(1)} \geqslant 0$ according to Remark \ref{REM6} and in this case, if $(t,x) \in \mathfrak{C}\big( p_1 + \varepsilon, p_2 \big)$ then
	\begin{equation*}
		\varepsilon \leqslant \frac{x}{2t} - p_1 \qquad \Longrightarrow \qquad \left( \frac{x}{2t} - p_1\right)^{-\tilde{\alpha}_k^{(1)}} \leqslant \varepsilon^{-\tilde{\alpha}_k^{(1)}} \; .
	\end{equation*}
	By defining
	\begin{equation*}
		\forall \, k \in \{1,...,9\} \qquad C_k^{(1)}(u_0, \varepsilon) := \tilde{R}_k^{(1)} \hspace{-1mm} \left( \frac{1}{2 \pi} \, \tf u_0 \right) \varepsilon^{-\tilde{\alpha}_k^{(1)}} \; ,
	\end{equation*}
	we obtain the result for the case $\mu > \frac{1}{2}$.\\
	We define in a similar way
	\begin{align*}
		& \bullet \hspace{5mm} \forall \, k \in \{1,...,8\} \qquad  C_k^{(2)}(u_0, \varepsilon) := \tilde{R}_k^{(2)} \hspace{-1mm} \left( \frac{1}{2 \pi} \, \tf u_0 \right) \varepsilon^{-\tilde{\alpha}_k^{(2)}} \; ; \\
		& \bullet \hspace{5mm} \forall \, k \in \{1,...,9\} \qquad C_k^{(3)}(u_0, \varepsilon) := \tilde{R}_k^{(3)} \hspace{-1mm} \left( \frac{1}{2 \pi} \, \tf u_0 \right) \varepsilon^{-\tilde{\alpha}_k^{(3)}} \; .
	\end{align*}
\end{proof}

\begin{REM0} \label{REM0}
	\em At this stage, the authors of \cite{fam} introduced the large parameter $\omega := \sqrt{t^2 + x^2}$, and replaced $t$ and $x$ by the bounded parameters $\tau := \frac{t}{\omega}$ and $\chi := \frac{x}{\omega}$. This led to a family of phase functions which was globally bounded in $\mathcal{C}^4$ with respect to $\tau$ and $\chi$. This was necessary for the application of \cite{hormander}. In our context, it is sufficient to control the phase functions in space-time cones. Indeed the explicitness of our remainder estimates shows that their coefficients depend only on the quotient $\frac{x}{t}$, which is bounded in these cones. It is not necessary to have the global boundedness with respect to $t$ and $x$ separately. Therefore we can use $t$ as a large parameter instead of $\omega$, which is conceptually simpler and clearer.
\end{REM0}

In the second result, we study the solution outside the cone $\mathfrak{C} \big( p_1 , p_2 \big)$. In this case, the stationary point is outside the integration interval and so the decay rate is only governed by the singular frequency. As above, we have to restrict slightly the space-time cones in which we expand the solution to bound uniformly the remainder.

\begin{THM7}
	Suppose that $u_0$ satisfies Condition \emph{(C$_{[p_1,p_2],\mu}$)} and choose $\varepsilon > 0$ such that
	\begin{equation*}
		-\varepsilon^{-1} < p_1 - \varepsilon \qquad \text{and} \qquad  p_2 + \varepsilon < \varepsilon^{-1} \; .
	\end{equation*}
	Define 
	\begin{align*}
		& \bullet \hspace{2mm} \forall \, (t,x) \in \mathfrak{C}\left(-\varepsilon^{-1} \, , \, p_1 - \varepsilon \right) \quad K_{\mu}^{(1)}(t,x,u_0) := \frac{\Gamma(\mu)}{2^{\mu+1} \pi} \, e^{- i \frac{\pi \mu}{2}} \, e^{-it p_1^{\, 2} + i x p_1} \, \tilde{u}(p_1) \hspace{-1pt} \left( p_1 - \frac{x}{2t} \right)^{-\mu} \\
		& \bullet \hspace{2mm} \forall \, (t,x) \in \mathfrak{C}\left( p_2 + \varepsilon \, , \, \varepsilon^{-1} \right) \quad K_{\mu}^{(2)}(t,x,u_0) := \frac{\Gamma(\mu)}{2^{\mu+1} \pi} \, e^{i \frac{\pi \mu}{2}} \, e^{-it p_1^{\, 2} + i x p_1} \, \tilde{u}(p_1) \hspace{-1pt} \left( \frac{x}{2t} - p_1 \right)^{-\mu} \, ,
	\end{align*}	
	Then
	\begin{itemize}
		\item for all $(t,x) \in \mathfrak{C}\big(-\varepsilon^{-1} \, , \, p_1 - \varepsilon \big)$, we have
	\begin{equation*}
		\left| u(t,x) - K_{\mu}^{(1)}(t,x,u_0) \, t^{-\mu} \right| \leqslant C^{(1)}(u_0,\varepsilon) \, t^{-1} \; .
	\end{equation*}
	The constant $C^{(1)}(u_0,\varepsilon) \geqslant 0$ is given in the proof ;
	\item for all $(t,x) \in \mathfrak{C}\big( p_2 + \varepsilon \, , \, \varepsilon^{-1} \big)$, we have
	\begin{equation*}
		\left| u(t,x) - K_{\mu}^{(2)}(t,x,u_0) \, t^{-\mu} \right| \leqslant C^{(2)}(u_0,\varepsilon) \, t^{-1} \; .
	\end{equation*}
	The constant $C^{(2)}(u_0,\varepsilon) \geqslant 0$ is given in the proof.
	\end{itemize}
\end{THM7}

\begin{proof}
	The first step of the proof consists in rewritting the solution as an oscillatory integral, namely,
	\begin{equation*}
		u(t,x) = \frac{1}{2\pi} \int_{p_1}^{p_2} \tf u_0(p) \, e^{-i tp^2 + ixp} \, dp = \int_{p_1}^{p_2} U(p) \, e^{i t \psi(p)} \, dp \; ,
	\end{equation*}
	where the functions $U$ and $\psi$ are defined at the beginning of the proof of Theorem \ref{THM6}.\\
	Here the hypothesis on $(t,x)$ implies
	\begin{align*}
		& \bullet \hspace{3mm} (t,x) \in \mathfrak{C}\left(-\varepsilon^{-1} \, , \, p_1 - \varepsilon \right) \qquad \Longleftrightarrow \qquad -\varepsilon^{-1} < p_0 < p_1 - \varepsilon \; , \\
		& \bullet \hspace{3mm} (t,x) \in \mathfrak{C}\left( p_2 + \varepsilon \, , \, \varepsilon^{-1} \right) \qquad \Longleftrightarrow \qquad p_2 + \varepsilon < p_0 < \varepsilon^{-1} \; ,
	\end{align*}
	where $p_0 := \frac{x}{2t}$ is the unique stationary point of the phase. Hence the stationary point does not belong to the integration interval $[p_1,p_2]$ in both situations. Note that in the cone $\mathfrak{C}\big(-\varepsilon^{-1} \, , \, p_1 - \varepsilon \big)$, the phase is decreasing so we can carry out the substitution $p \longmapsto -p$ to make it increasing. It follows that the function $U$ satisfies Assumption (A$_{\mu,1,1}$) and $\psi$ verifies Assumption (P$_{1,1,N}$) (for $N \geqslant 1$) on $[p_1,p_2]$ in both cases and so Theorem \ref{THM2} is applicable. Note that it is not necessary to use a cutting-point to split the integral since only $p_1$ is a singular point. This implies that we are going to use only the expansion of the integral $\tilde{I}^{(1)}(\omega,q)$ (see the proof of Theorem \ref{THM2}) with $\omega = t$ and $q = p_2$ in the present proof. By applying the estimate of the remainder $R_1^{(2)}(\omega, q)$ and by using the inequalities $ -\varepsilon^{-1} < p_0 < p_1 - \varepsilon$ in $\mathfrak{C}\big(-\varepsilon^{-1} \, , \, p_1 - \varepsilon \big)$, and the estimate of $R_1^{(1)}(\omega, q)$ and $p_2 + \varepsilon < p_0 < \varepsilon^{-1}$ in the other cone, we obtain finally the expressions of the constants given in the statement of the present theorem :	\begin{align*}
		& \bullet \hspace{3mm} C^{(1)}(u_0,\varepsilon) := \frac{1}{4 \pi \mu} \left( p_2 + \varepsilon^{-1} \right) (p_2 - p_1)^{\mu} \, \big\| \tilde{u} \big\|_{W^{1,\infty}(p_1,p_2)} \bigg( \frac{1-\mu}{2} \, \left( p_2 + \varepsilon^{-1} \right) \varepsilon^{-4} \\
		& \hspace{4cm} + \varepsilon^{-2} + \varepsilon^{-3} \bigg) \; , \\
		& \bullet \hspace{3mm} C^{(2)}(u_0,\varepsilon) := \frac{1}{4 \pi \mu} \left( \varepsilon^{-1} - p_1 \right) (p_2 - p_1)^{\mu} \big\| \tilde{u} \big\|_{W^{1,\infty}(p_1,p_2)} \bigg( \frac{1-\mu}{2} \, \left( \varepsilon^{-1} - p_1 \right) \varepsilon^{-4} \\
		& \hspace{4cm} + \varepsilon^{-2} + \varepsilon^{-3} \bigg) \; .
	\end{align*}
\end{proof}

The next result is a consequence of Theorem \ref{THM6}. It permits to compute the limit of the $L^2$-norm of the solution on the spatial cross-section of the cone $\mathfrak{C}\big(p_1 + \varepsilon,p_2 \big)$ when the time tends to infinity, assuming $u_0 \in L^2(\R)$.

\begin{COR3} \label{COR3}
	Suppose that the hypotheses of Theorem \ref{THM6} are satisfied with $\mu > \frac{1}{2}$ and define the interval $I_t$ as follows,
	\begin{equation*}
		I_t := \Big( \, 2 \, (p_1 + \varepsilon) \,  t, \, 2 \, p_2 \, t \, \Big) \; .
	\end{equation*}
	Then we have
	\begin{equation*}
		\left| \big\| u(t,.) \big\|_{L^2(I_t)} - \frac{1}{\sqrt{2 \pi}} \, \big\| \tf u_0 \big\|_{L^2(p_1+\varepsilon, p_2)} \right| \leqslant \sum_{k=1}^9 \tilde{C}_k^{(1)}(u_0,\varepsilon) \, t^{- \tilde{\beta}_k^{(1)} + \frac{1}{2}} \; ,
	\end{equation*}
	where $\displaystyle \max_{k \in \{1,...9\}} \left\{ -\tilde{\beta}_k^{(1)} \right\} < - \frac{1}{2}$. The exponents $\tilde{\beta}_k^{(1)}$ are defined in Theorem \ref{THM6} and the constants $\tilde{C}_k^{(1)}(u_0, \varepsilon) \geqslant 0$ are given in the proof.
\end{COR3}

\begin{proof}
	We start by using the triangle inequality as follows,
	\begin{align*}
		\left| \big\| u(t,.) \big\|_{L^2(I_t)} - \left\| H(t,.,u_0) \, t^{-\frac{1}{2}} \right\|_{L^2(I_t)} \right|^2	& \leqslant \left\| u(t,.) - H(t,x,u_0) \, t^{-\frac{1}{2}} \right\|_{L^2(I_t)}^2 \\
						& = \int_{I_t} \left| u(t,x) - H(t,x,u_0) \, t^{-\frac{1}{2}} \right|^2 \, dx \; .
	\end{align*}
	Now we employ the estimate provided by Theorem \ref{THM6} in the case $\mu > \frac{1}{2}$ to control the last term :
	\begin{align*}
		\int_{I_t} \left| u(t,x) - H(t,x,u_0) \, t^{-\frac{1}{2}} \right|^2 \, dx	& \leqslant \int_{I_t} \left| \sum_{k = 1}^9 \, C_k^{(1)}(u_0, \varepsilon) \, t^{-\tilde{\beta}_k^{(1)}} \right|^2 \, dx \\
				& = \sum_{k = 1}^9 \, C_k^{(1)}(u_0, \varepsilon)^2 \, t^{-2 \tilde{\beta}_k^{(1)}} \, |I_t| \; ,	 
	\end{align*}
	where $|I_t|$ is the Lebesgue measure of the interval $I_t$; in our context, $|I_t|$ is equal to $2 ( p_2 - p_1 - \varepsilon) \, t$ . By defining
	\begin{equation*}
		\tilde{C}_k^{(1)}(u_0, \varepsilon) := \sqrt{2 ( p_2 - p_1 - \varepsilon)} \, C_k^{(1)}(u_0, \varepsilon) \; ,
	\end{equation*}
	we obtain
	\begin{equation*}
		\left| \big\| u(t,.) \big\|_{L^2(I_t)} - \left\| H(t,.,u_0) \, t^{-\frac{1}{2}} \right\|_{L^2(I_t)} \right| \leqslant \sum_{k = 1}^9 \, \tilde{C}_k^{(1)}(u_0, \varepsilon) \, t^{-\tilde{\beta}_k^{(1)} + \frac{1}{2}} \; .
	\end{equation*}
	To finish, we compute $\left\| H(t,.,u_0) \, t^{-\frac{1}{2}} \right\|_{L^2(I_t)}$ by using the expression of $H(t,x,u_0)$ given in Theorem \ref{THM6},
	\begin{align*}
		\left\| H(t,.,u_0) \, t^{-\frac{1}{2}} \right\|_{L^2(I_t)}^2	& = \frac{1}{4 \pi} \, \int_{I_t} \left| \tilde{u} \hspace{-1mm} \left(\frac{x}{2t} \right) \left( \frac{x}{2t} - p_1 \right)^{\mu-1} t^{-\frac{1}{2}} \right|^2 \, dx \\
												& = \frac{t}{2 \pi} \int_{p_1 + \varepsilon}^{p_2} \left| \tilde{u}(y) (y-p_1)^{\mu-1} \right|^2 \, dy \, t^{-1} \\
												& = \frac{1}{2 \pi} \big\| \tf u_0 \big\|_{L^2(p_1 + \varepsilon, p_2)}^2 \; .
	\end{align*}
	The proof is now complete.
\end{proof}

The aim of the next result is to show that the time decay rate is $t^{-\frac{\mu}{2}}$ on points moving in space-time with the critical velocity given by the singularity $p_1$. In this case, we do not have to deal with the uniformity of the constants of the remainder since we establish an asymptotic expansion on a line.

\begin{THM8} \label{THM8}
	Suppose that $u_0$ satisfies Condition \emph{(C$_{[p_1,p_2],\mu}$)}. For all $t > 0$, define $L_{\mu}(t,u_0)$ as follows :
	\begin{equation*}
		 \bullet \quad L_{\mu}(t,u_0) := \frac{1}{2} \, \Gamma \hspace{-1mm} \left( \frac{\mu}{2} \right) e^{-i \frac{\pi \mu}{4}} e^{i t p_1^{\, 2}} \, \tilde{u}(p_1) \; .
	\end{equation*}
	Then for all $(t,x) \in (0,+\infty) \times \R$ such that $\frac{x}{t} = 2 \, p_1$,	we have
	\begin{equation*}
		\left| u(t,x) - L_{\mu}(t,u_0) \, t^{-\frac{\mu}{2}} \right| \leqslant C(u_0) \, t^{-\frac{1}{2}} \; .
	\end{equation*}
	The constant $C(u_0) \geqslant 0$ is given in the proof.
\end{THM8}

\begin{proof}
	In this situation, the stationary point $p_0 = \frac{x}{2t}$ and the singularity $p_1$ are equal but the phase is decreasing. So we make the substitution $p \longmapsto -p$ and then we apply Theorem \ref{THM2}. Since $p_1$ is the unique singular point inside the integration interval, we can employ only the expansion of the integral $\tilde{I}^{(2)}(\omega,q)$ given in the proof of Theorem \ref{THM2}. Let us give the constant $C(u_0)$ to conclude the proof,
	\begin{equation*}
		\bullet \quad C(u_0) := \frac{\sqrt{\pi}}{2 \mu} \, (p_2 - p_1)^{\mu} \, \big\| \tilde{u}' \big\|_{L^{\infty}(p_1,p_2)} \; .
	\end{equation*}
\end{proof}

Uniform estimates of the solution in the curved region $\mathfrak{R}_{\varepsilon}(p_1, p_2)$, which is asymptotically larger than any space-time cone contained in $\mathfrak{C}(p_1, p_2)$, are provided in the next theorem. For this purpose, we rewrite the solution as an oscillatory integral that we estimate by using Corollary \ref{COR2}. Here the quantities $\big( \frac{x}{2t} - p_1 \big)^{-\tilde{\alpha}_k^{(j)}}$, which produce the blow-up, are bounded by $t^{\varepsilon \, \tilde{\alpha}_k^{(j)}}$ furnishing a uniform estimate of the solution with modified decay rates for sufficiently small $\varepsilon > 0$. Finally we use Theorem \ref{THM5} to give the preponderant decay rate.

\begin{THM9} \label{THM9}
	Suppose that $u_0$ satisfies Condition \emph{(C$_{[p_1,p_2],\mu}$)} and fix $\varepsilon \in \big(0, \frac{1}{2}\big)$. Then for all $(t,x) \in \mathfrak{R}_{\varepsilon}(p_1, p_2)$, we have
	\begin{itemize}
	\item Case $\displaystyle \mu > \frac{1}{2}$ :
	\begin{equation*}
		\big| u(t,x) \big| \leqslant C_0^{(1)}(u_0) \, t^{-\frac{1}{2} + \varepsilon(1-\mu)} + \sum_{k=1}^9 C_k^{(1)} (u_0) \, t^{- \tilde{\beta}_k^{(1)} + \varepsilon \tilde{\alpha}_k^{(1)}} \; ,
	\end{equation*}
	where $\displaystyle \max_{k \in \{1,..., 9\}} \left\{ -\tilde{\beta}_k^{(1)} + \varepsilon \tilde{\alpha}_k^{(1)} \right\} < -\frac{1}{2} + \varepsilon(1-\mu)$ and the decay rate $t^{-\frac{1}{2} + \varepsilon(1-\mu)}$ is optimal. The exponents $\tilde{\alpha}_k^{(1)}$, $\tilde{\beta}_k^{(1)}$ are defined in Corollary \ref{COR2} and the constants $C_k^{(1)}(u_0) \geqslant 0$ are given in the proof ;
	\item Case $\displaystyle \mu = \frac{1}{2}$ :
	\begin{equation*}
		\big| u(t,x) \big| \leqslant C_0^{(2)}(u_0) \, t^{-\frac{1}{2} + \frac{\varepsilon}{2}} + \sum_{k=1}^8 C_k^{(2)} (u_0) \, t^{- \tilde{\beta}_k^{(2)} + \varepsilon \tilde{\alpha}_k^{(2)}} \; ,
	\end{equation*}
	where $\displaystyle \max_{k \in \{1,..., 8\}} \left\{ -\tilde{\beta}_k^{(2)} + \varepsilon \tilde{\alpha}_k^{(2)} \right\} < -\frac{1}{2} + \frac{\varepsilon}{2}$ and the decay rate $t^{-\frac{1}{2} + \frac{\varepsilon}{2}}$ is optimal. The exponents $\tilde{\alpha}_k^{(2)}$, $\tilde{\beta}_k^{(2)}$ are defined in Corollary \ref{COR2} and the constants $C_k^{(2)}(u_0) \geqslant 0$ are given in the proof ;
	\item Case $\displaystyle \mu < \frac{1}{2}$ :	
	\begin{equation*}
		\big| u(t,x) \big| \leqslant C_0^{(3)}(u_0) \, t^{-\mu + \varepsilon \mu} + \sum_{k=1}^9 C_k^{(3)} (u_0) \, t^{- \tilde{\beta}_k^{(3)} + \varepsilon \tilde{\alpha}_k^{(3)}} \; ,
	\end{equation*}
	where $\displaystyle \max_{k \in \{1,..., 9\}} \left\{ -\tilde{\beta}_k^{(3)} + \varepsilon \tilde{\alpha}_k^{(3)} \right\} < -\mu + \varepsilon \mu$ and the decay rate $t^{-\mu + \varepsilon \mu}$ is optimal. The exponents $\tilde{\alpha}_k^{(3)}$, $\tilde{\beta}_k^{(3)}$ are defined in Corollary \ref{COR2} and the constants $C_k^{(3)}(u_0) \geqslant 0$ are given in the proof.
	\end{itemize}
\end{THM9}

\begin{proof}
	Assuming that the hypotheses of Corollary \ref{COR2} are satisfied, we derive  the following estimate of the oscillatory integral in the case $\mu > \frac{1}{2}$ ,
	\begin{align}
		\left| \int_{p_1}^{p_2} U(p) \, e^{i \omega \psi(p)} \, dp \right|	& \leqslant \left| \tilde{H}(\omega, p_0, U) \right| (p_0 - p_1)^{\mu - 1} \, \omega^{-\frac{1}{2}} + \sum_{k=1}^9 \tilde{R}_k^{(1)}(U) (p_0 - p_1)^{-\tilde{\alpha}_k^{(1)}} \, \omega^{- \tilde{\beta}_k^{(1)}} \nonumber \\
																			& \label{estannexe} \leqslant \sqrt{\pi} \big\| \tilde{u} \big\|_{L^{\infty}(p_1,p_2)} \, (p_0 - p_1)^{\mu - 1} \, \omega^{-\frac{1}{2}} + \sum_{k=1}^9 \tilde{R}_k^{(1)}(U) (p_0 - p_1)^{-\tilde{\alpha}_k^{(1)}} \, \omega^{- \tilde{\beta}_k^{(1)}} \; .
	\end{align}
	We are going to use this estimate to establish the theorem.\\
	Let us establish the result in the case $\mu > \frac{1}{2}$. First we rewrite the solution formula \eqref{solformula} as an oscillatory integral (see the beginning of the proof of Theorem \ref{THM6}). The hypothesis $(t,x) \in \mathfrak{R}_{\varepsilon}(p_1,p_2)$ implies
	\begin{equation} \label{estannexe2}
		\frac{x}{2t} - p_1 \geqslant t^{-\varepsilon} \; ,
	\end{equation}
	and we can easily see that $p_0 := \frac{x}{2t} \in (p_1, p_2)$. Hence we observe that the hypotheses of Corollary \ref{COR2} are satisfied and in particular, estimate \eqref{estannexe} is applicable, furnishing
	\begin{equation*}
		\big| u(t,x) \big| \leqslant \frac{1}{2 \sqrt{\pi}} \, \big\| \tilde{u} \big\|_{L^{\infty}(p_1,p_2)} \left( \frac{x}{2t} - p_1 \right)^{\mu - 1} t^{-\frac{1}{2}} \, + \sum_{k=1}^9 \tilde{R}_k^{(1)} \hspace{-1mm} \left(\frac{1}{2\pi} \, \tf u_0 \right) \left( \frac{x}{2t} - p_1 \right)^{-\tilde{\alpha}_k^{(1)}} t^{- \tilde{\beta}_k^{(1)}} \; .
	\end{equation*}
	According to Remark \ref{REM6}, each $\tilde{\alpha}_k^{(1)}$ is non-negative if the parameter $\delta \in \big( \frac{1}{2} , 1 \big)$ appearing in the proof of Corollary \ref{COR2} is such that $\delta \geqslant \frac{\mu + 1}{2}$; so let us suppose that $\delta \geqslant \frac{\mu + 1}{2}$. Hence we can put estimate \eqref{estannexe2} into the last inequality and we obtain
	\begin{equation*}
		\big| u(t,x) \big| \leqslant \frac{1}{2 \sqrt{\pi}} \, \big\| \tilde{u} \big\|_{L^{\infty}(p_1,p_2)} \, t^{-\frac{1}{2} + \varepsilon ( 1 - \mu)} + \sum_{k=1}^9 \tilde{R}_k^{(1)} \hspace{-1mm} \left(\frac{1}{2\pi} \, \tf u_0 \right) \, t^{- \tilde{\beta}_k^{(1)} + \varepsilon \tilde{\alpha}_k^{(1)}} \; .
	\end{equation*}
	Thanks to Theorem \ref{THM5}, we know that each exponent $- \tilde{\beta}_k^{(1)} + \varepsilon \tilde{\alpha}_k^{(1)}$ is strictly smaller than $-\frac{1}{2} + \varepsilon ( 1 - \mu)$ if $\delta > \mu$ (which is true since $\delta \geqslant \frac{\mu + 1}{2}$) and if $\delta > \frac{1}{2} + \varepsilon$. Defining for all $k \in \{1,...,9\}$,
	\begin{equation*}
		C_0^{(1)}(u_0) := \frac{1}{2 \sqrt{\pi}} \, \big\| \tilde{u} \big\|_{L^{\infty}(p_1,p_2)} \qquad , \qquad C_k^{(1)}(u_0) := \tilde{R}_k^{(1)} \hspace{-1mm} \left(\frac{1}{2\pi} \, \tf u_0 \right) \; ,
	\end{equation*}
	we obtain the result in the case $\mu > \frac{1}{2}$. The optimality is a direct consequence of Theorem \ref{THM10}.\\
	We employ the same arguments to establish the estimates in the two other cases, and we define
	\begin{align*}
		& \bullet \hspace{5mm} C_0^{(2)}(u_0) := \left( \frac{1}{2 \sqrt{\pi}} \, + \frac{\Gamma(\mu)}{2^{\mu + 1}} \right) \big\| \tilde{u} \big\|_{L^{\infty}(p_1,p_2)} \quad , \quad C_k^{(2)}(u_0) := \tilde{R}_k^{(2)} \hspace{-1mm} \left(\frac{1}{2\pi} \, \tf u_0 \right) \; ; \\
		& \bullet \hspace{5mm} C_0^{(3)}(u_0) := \frac{\Gamma(\mu)}{2^{\mu + 1}} \, \big\| \tilde{u} \big\|_{L^{\infty}(p_1,p_2)} \quad , \quad C_k^{(3)}(u_0) := \tilde{R}_k^{(3)} \hspace{-1mm} \left(\frac{1}{2\pi} \, \tf u_0 \right) \; ,
	\end{align*}
	for $k \geqslant 1$.
\end{proof}

The last result is devoted to the optimality of the previous uniform estimates. In the region $\mathfrak{R}_{\varepsilon}(p_1 , p_2)$, we expect that the decay will be slow in parts which are close to the critical direction given by $p_1$, where the influence of the singularity is the strongest. So we use Theorem \ref{THM5} to provide asymptotic expansions of the solution on the space-time curve $\mathfrak{G}_{\varepsilon}(p_1)$, the left boundary of the region $\mathfrak{R}_{\varepsilon}(p_1,p_2)$, and we show that the decay rates obtained in the preceding result are attained on this curve, proving the optimality.

\begin{THM10} \label{THM10}
	Suppose that $u_0$ satisfies Condition \emph{(C$_{[p_1,p_2],\mu}$)} and fix $\varepsilon \in \big(0, \frac{1}{2}\big)$. For all $t > T_{\varepsilon}(p_1,p_2)$, define $H(t,u_0) \in \C$ and $K_{\mu}(t,u_0) \in \C$ as follows :
	\begin{align*}
		& \bullet \hspace{5mm} H(t,u_0) := \frac{1}{2 \sqrt{\pi}} \, e^{-i \frac{\pi}{4}} \, e^{it(p_1 + t^{-\varepsilon})^2} \, \tilde{u}\big( p_1 + t^{-\varepsilon} \big) \; ; \\
		& \bullet \hspace{5mm} K_{\mu}(t,u_0) := \frac{\Gamma(\mu)}{2^{\mu}} \, e^{i \frac{\pi \mu}{2}} \, e^{-i t p_1^{\, 2} + i x p_1} \, \tilde{u}(p_1) \; .
	\end{align*}
	Then for all $(t,x) \in \mathfrak{G}_{\varepsilon}(p_1)$ with $t > T_{\varepsilon}(p_1,p_2)$, we have
	\begin{itemize}
	\item Case $\displaystyle \mu > \frac{1}{2}$ :	
	\begin{equation*}
		\Big| u(t,x) - H(t,u_0) \, t^{-\frac{1}{2} + \varepsilon(1-\mu)} \, \Big| \leqslant \sum_{k=1}^9 C_k^{(1)}(u_0) \, t^{-\tilde{\beta}_k^{(1)} + \varepsilon \tilde{\alpha}_k^{(1)}} \; ,
	\end{equation*}
	where $\displaystyle \max_{k \in \{1,..., 9\}} \left\{ -\tilde{\beta}_k^{(1)} + \varepsilon \tilde{\alpha}_k^{(1)} \right\} < -\frac{1}{2} + \varepsilon(1-\mu)$. The exponents $\tilde{\alpha}_k^{(1)}$, $\tilde{\beta}_k^{(1)}$ are defined in Corollary \ref{COR2} and the constants $C_k^{(1)}(u_0) \geqslant 0$ are given in Theorem \ref{THM9} ;
	\item Case $\displaystyle \mu = \frac{1}{2}$ :
	\begin{equation*}
		\Big| u(t,x) - \big(K_{\mu}(t,u_0) + H(t,u_0) \big) \, t^{-\frac{1}{2} + \frac{\varepsilon}{2}} \, \Big| \leqslant \sum_{k=1}^8 C_k^{(2)}(u_0) \, t^{-\tilde{\beta}_k^{(2)} + \varepsilon \tilde{\alpha}_k^{(2)}} \; ,
	\end{equation*}
	where $\displaystyle \max_{k \in \{1,..., 8\}} \left\{ -\tilde{\beta}_k^{(2)} + \varepsilon \tilde{\alpha}_k^{(2)} \right\} < -\frac{1}{2} + \frac{\varepsilon}{2}$. The exponents $\tilde{\alpha}_k^{(2)}$, $\tilde{\beta}_k^{(2)}$ are defined in Corollary \ref{COR2} and the constants $C_k^{(2)}(u_0) \geqslant 0$ are given in Theorem \ref{THM9} ;
	\item Case $\displaystyle \mu < \frac{1}{2}$ :	
	\begin{equation*}
		\Big| u(t,x) - K_{\mu}(t,u_0) \, t^{-\mu + \varepsilon \mu} \, \Big| \leqslant \sum_{k=1}^9 C_k^{(3)}(u_0) \, t^{-\tilde{\beta}_k^{(3)} + \varepsilon \tilde{\alpha}_k^{(3)}} \; ,
	\end{equation*}
	where $\displaystyle \max_{k \in \{1,..., 9\}} \left\{ -\tilde{\beta}_k^{(3)} + \varepsilon \tilde{\alpha}_k^{(3)} \right\} < -\mu + \varepsilon \mu$. The exponents $\tilde{\alpha}_k^{(3)}$, $\tilde{\beta}_k^{(3)}$ are defined in Corollary \ref{COR2} and the constants $C_k^{(3)}(u_0) \geqslant 0$ are given in Theorem \ref{THM9}.
	\end{itemize}
\end{THM10}

\begin{proof}
	The hypothesis $(t,x) \in \mathfrak{G}_{\varepsilon}(p_1)$ is equivalent to $p_0 - p_1 = t^{-\varepsilon}$ where $p_0 := \frac{x}{2t}$. Using the solution formula \eqref{solformula} and the rewritting given in the proof of Theorem \ref{THM6}, we can apply Theorem \ref{THM5} and we obtain the desired estimates.	
\end{proof}

\section{The core of the method: oscillation control by complex analysis}

\hspace{2.5ex} This final section contains the technical but crucial arguments and calculations to fill the considerable gaps left in the original sketch of the proof of Erdélyi. The results will be presented in the order they appear in the proof of Erdélyi's stationary phase method.\\

Throughout this section, the parameter $\omega > 0$ will be fixed and the integer $j$ will belong to $\{1,2\}$. We shall prove the propositions in the case $j=1$ only; the proofs in the case $j=2$ are very similar.\\

At the beginning of the proof of Erdélyi's theorem, a change of variables is carried out in order to simplify the phase. The aim of the following proposition is to prove that this change of variables is admissible. To this end, we prove that the associated transformation is a diffeomorphism by exploiting substantially the factorization of the zeros of the derivative of the phase.

\begin{LEM1} \label{LEM1}
	Fix $q_1, q_2 \in (p_1,p_2)$. Let $\psi : [p_1,p_2] \longrightarrow \R$ be a function satisfying Assumption \emph{(P$_{\rho_1,\rho_2,N}$)} and consider the functions $\varphi_1 : [p_1,q_1] \longrightarrow \R$ and $\varphi_2 : [q_2,p_2] \longrightarrow \R$  defined by
	\begin{equation*}
		j=1,2 \qquad \varphi_j(p) = \left( (-1)^{j+1} \big( \psi(p) - \psi(p_j) \big) \right)^{\frac{1}{\rho_j}} \; .
	\end{equation*}
	Then $\varphi_1$ (resp. $\varphi_2$) is a $\mathcal{C}^{N+1}$-diffeomorphism between $[p_1,q_1]$ (resp. $[q_2,p_2]$) and $\big[0,\varphi_1(q_1) \big]$ (resp. $\big[0,\varphi_2(q_2) \big]$).
\end{LEM1}

\begin{proof}
	First of all, we check that $\varphi_1 \in \mathcal{C}^{N+1}\big( [p_1,q_1] , \R \big)$. We recall that
	\begin{equation*}
		\psi'(p) = (p-p_1)^{\rho_1-1} \tilde{\psi}_2(p) \; ,
	\end{equation*}
	where we put $\tilde{\psi}_2(p):=(p_2-p)^{\rho_2-1} \tilde{\psi}(p)$, which belongs to $\mathcal{C}^{N+1}\big( [p_1,q_1] , \R \big)$. Applying Taylor's Theorem with the integral form of the remainder to $\psi'$, we obtain the following representation of $\varphi_1$ :
	\begin{align*}
		\forall \, p \in [p_1 , q_1] \qquad \varphi_1(p)	& = (p-p_1) \bigg( \int_0^1 y^{\rho_1 -1} \, \tilde{\psi}_2 \big( y(p-p_1) + p_1 \big) \, dy \bigg)^{\frac{1}{\rho_1}} \\
																		& =: (p-p_1) \, J_1(p)^{\frac{1}{\rho_1}} \; .
	\end{align*}
	We fix $k \in \{1,\ldots,N\}$ and we compute formally the $k$th derivative of the above expression by using the product rule :
	\begin{equation} \label{deriv}
		\frac{d^k}{dp^k} \big[\varphi_1\big] (p) = (p-p_1) \frac{d^k}{dp^k} \hspace{-1mm} \left[ J_1^{\, \frac{1}{\rho_1}} \right]\hspace{-1mm} (p) + k \, \frac{d^{k-1}}{dp^{k-1}} \hspace{-1mm} \left[ J_1^{\, \frac{1}{\rho_1}} \hspace{-1mm} \right](p) \; .
	\end{equation}
	The positivity and the regularity of the function $\tilde{\psi}_2$ allow to differentiate $k$ times under the integral sign the function $J_1$. Hence the $k$ first derivatives of the composite function $\displaystyle J_1^{\, \frac{1}{\rho_1}}$ exist and are continuous; in particular, the expression \eqref{deriv} is well-defined for all $p \in [p_1,q_1]$ and $\frac{d^k}{dp^k} \big[\varphi_1\big]$ is continuous. Concerning the $(N+1)$th derivative, we must be careful because we have not supposed that $\tilde{\psi} \in \mathcal{C}^{N+1}\big( [p_1,p_2] , \R \big)$. However we can formally apply the product rule to $\varphi_1$ once again for $k=N+1$ :
	\begin{equation*}
		\frac{d^{N+1}}{dp^{N+1}} \big[\varphi_1\big] (p) = (p-p_1) \,  \underbrace{\frac{d^{N+1}}{dp^{N+1}} \hspace{-1mm} \left[ J_1^{\, \frac{1}{\rho_1}} \right] \hspace{-1mm} (p)}_{(i)} + (N+1) \, \underbrace{\frac{d^{N}}{dp^{N}} \hspace{-1mm} \left[ J_1^{\,\frac{1}{\rho_1}} \right] \hspace{-1mm}(p)}_{(ii)} \; .
	\end{equation*}
	Note that the term $(ii)$ is well-defined by the previous work. So it remains to study the term $(i)$. Firstly, let us define the function $h_1 : s \longmapsto s^{\frac{1}{\rho_1}}$. Then we obtain by applying Fa\`a di Bruno's Formula to $J_1^{\, \frac{1}{\rho_1}} = h_1 \circ J_1$, 
	\begin{equation*} \label{deriv2}
		\frac{d^{N+1}}{dp^{N+1}} \hspace{-1mm} \left[ J_1^{\, \frac{1}{\rho_1}} \right](p) = \sum C_N \underbrace{\left( \frac{d^{m_1+\ldots+m_{N+1}}}{dp^{m_1+\ldots+m_{N+1}}} \big[h_1\big]  \circ J_1 \right) \hspace{-1mm} (p)}_{(iii)} \, \prod_{l=1}^{N+1} \underbrace{\left( \frac{d^l}{dp^l} \big[J_1\big](p) \right)^{m_l}}_{(iv)}
	\end{equation*}
	where the sum is over all the $(N+1)$-tuples $(m_1, ..., m_{N+1})$ satisfying: $1m_1+2m_2+3m_3+\ldots+(N+1)m_{N+1}=N+1$. We note that the term $(iii)$ is well-defined by the positivity of $J_1$; moreover by the previous study, the term $(iv)$ is well-defined and continuous for any $l \neq N+1$. So we have to study $\displaystyle \left( \frac{d^{N+1}}{dp^{N+1}} \big[J_1\big](p) \right)^{m_{N+1}}$, where $m_{N+1} \leqslant 1$ by the above constraint. Since the case $m_{N+1} = 0$ is clear, we suppose that $m_{N+1} = 1$. We have
	\begin{equation*}
		\begin{aligned}
			\frac{d}{dp} \hspace{-1mm} \left[ \frac{d^N}{dp^N} \big[ J_1 \big] \right] \hspace{-1mm} (p)	& = \frac{d}{dp} \hspace{-1mm} \left[ \int_0^1 y^{N+\rho_1-1} \frac{d^N}{dp^N} \big[ \tilde{\psi}_2 \big] \big((p-p_1)y+p_1\big) \, dy \right] \\
							& = \frac{d}{dp} \hspace{-1mm} \left[ \frac{1}{(p-p_1)^{\rho_1+N}} \int_{p_1}^p (s-p_1)^{\rho_1+N-1} \frac{d^N}{dp^N} \big[ \tilde{\psi}_2 \big](s) \, ds \right] \\
							& = \frac{-(\rho_1 +N)}{(p-p_1)^{\rho_1+N+1}} \int_{p_1}^p (s-p_1)^{\rho_1+N-1} \frac{d^N}{dp^N} \big[ \tilde{\psi}_2 \big](s) \, ds \\
							& \qquad + \frac{1}{(p-p_1)^{\rho_1+N}} \, (p-p_1)^{N+\rho_1 -1} \, \frac{d^N}{dp^N} \big[ \tilde{\psi}_2 \big](p) \\
							& = \frac{-(\rho_1 +N)}{(p-p_1)} \int_0^1 y^{\rho_1+N-1} \frac{d^N}{dp^N} \big[ \tilde{\psi}_2 \big] \big(y(p-p_1)+p_1 \big) \, dy \\
							& \qquad  + \frac{1}{(p-p_1)} \, \frac{d^N}{dp^N} \big[ \tilde{\psi}_2 \big](p) \; .
		\end{aligned}
	\end{equation*}
	Multiplying this equality by $(p-p_1)$, we observe that the function $\displaystyle p \in [p_1,q_1] \longmapsto (p-p_1) \frac{d^{N+1}}{dp^{N+1}} \left[J_1^{\frac{1}{\rho_1}} \right] \hspace{-1mm}(p)$ is well-defined and continuous. Then $\displaystyle \frac{d^{N+1}}{dp^{N+1}} \big[\varphi_1\big]$ is continuous on $[p_1,q_1]$, proving that $\varphi_1 \in \mathcal{C}^{N+1}\big( [p_1,q_1] , \R \big)$.\\
	Furthermore, one remarks that
	\begin{equation*}
	\left\{ \hspace{-0.3cm} \begin{array}{rl}
			& \displaystyle \varphi_1'(p) = \frac{1}{\rho_1} \, \psi'(p) \, \big(\psi(p) - \psi(p_1) \big)^{\frac{1}{\rho_1}-1} > 0 \qquad \forall \, p \in (p_1,q_1] \; , \\
			& \displaystyle \varphi_1'(p_1) = \frac{1}{\rho_1^{\rho_1}} \, \tilde{\psi}_2 (p_1)^{\frac{1}{\rho_1}} > 0 \; ,
	\end{array} \right. \; ,
\end{equation*}
	so by the inverse function theorem, $\varphi_1 \, : \, [p_1,q_1] \, \longrightarrow \, \big[ 0, \varphi_1(q_1) \big]$ is a $\mathcal{C}^{N+1}$-diffeomorphism.
\end{proof}

As a result of this change of variables in the proof of the result of Erdélyi, the integrand is factorized into a holomorphic function and a function on a real interval. The aim of the next result is to prove that this second function is regular.

\begin{LEM2} \label{LEM2}
	Let $U : (p_1,p_2) \longrightarrow \C$ be a function satisfying Assumption \emph{(A$_{\mu_1,\mu_2,N}$)} and consider the functions $k_j : \big(0, \varphi_j(q_j) \big] \longrightarrow \C$ defined by
	\begin{equation*}
		k_j(s) = U\big( \varphi_j^{-1}(s) \big) \, s^{1-\mu_j} \, \big( \varphi_j^{-1}\big)'(s) \; ,
	\end{equation*}
	where the functions $\varphi_j$ and the points $q_j$ are defined in Proposition \ref{LEM1}.\\
	Then $k_j$ can be extended to $\big[0,\varphi_j(q_j)\big]$ and belongs to $\displaystyle \mathcal{C}^{N}\left(\big[0,\varphi_j(q_j)\big],\C \right)$.
\end{LEM2}

\begin{proof}
	We define $\tilde{u}_2(p) := (p_2-p)^{\mu_2-1} \tilde{u}(p)$ for all $p \in \, [p_1,p_2)$. Then we have by the definition of $k_1$,
	\begin{align}
		k_1(s)	& = \left( \varphi_1^{-1}(s) - \varphi_1^{-1}(0) \right)^{\mu_1 -1} \, \tilde{u}_2 \big(\varphi_1^{-1}(s) \big) \, s^{1-\mu_1} \, \big(\varphi_1^{-1}\big)'(s) \nonumber \\
				& = \left( \frac{\varphi_1^{-1}(s) - \varphi_1^{-1}(0)}{s} \right)^{\mu_1-1} \tilde{u}_2 \big(\varphi_1^{-1}(s) \big) \big(\varphi_1^{-1}\big)'(s) \nonumber \\
				& \label{defk1} = \left( \int_0^1 \big( \varphi_1^{-1} \big)'(sy) \, dy \right)^{\mu_1-1} \tilde{u}_2 \big(\varphi_1^{-1}(s) \big) \big(\varphi_1^{-1}\big)'(s)
	\end{align}
	for all $s \in \, \big( 0, \varphi_1(q_1) \big]$, and $\displaystyle k_1(0) := \tilde{u}_2(p_1) \big(\varphi_1^{-1}\big)'(0)^{\mu_1}$ by taking the limit in \eqref{defk1}. The conclusion comes from the regularity of $\tilde{u}$ and Proposition \ref{LEM1}.
\end{proof}

The next step in the proof of Erdélyi consists in creating the expansion of the integral by the classical procedure of integrations by parts. Thanks to the above factorization of the integrand, we differentiate the regular function and we calculate successive primitives under integral forms of the holomorphic factor. The next task will be to exploit the holomorphy property and Cauchy's theorem to shift the integration path of the primitives in a region where we control the oscillations of the complex exponential in view of precise estimates of the remainder

In the following lemma, we establish an estimate of the complex exponential on a well-chosen half-line which will be the integration path of the primitives. This result will be essential in the proof of Theorem \ref{LEM4}.

\begin{LEM3} \label{LEM3}
	Let $\rho_j \geqslant 1$ and $s \geqslant 0$. Then we have
	\begin{equation*}
		\forall \, t \geqslant 0 \qquad \left| \, e^{(-1)^{j+1} i \omega \left( s + t e^{(-1)^{j+1} i \frac{\pi}{2\rho_j}} \right)^{\rho_j}} \right| \leqslant e^{- \omega t^{\rho_j}} \; .
	\end{equation*}
\end{LEM3}

\begin{proof}
	Let us fix $s,t \geqslant 0$. By a simple calculation, we furnish the following equality,
	\begin{equation} \label{inteq}
		i \rho_1 \omega \int_0^s \left( \xi + t e^{i \frac{\pi}{2 \rho_1}} \right)^{\rho_1 -1} d\xi = i \omega \left( s + t e^{i \frac{\pi}{2\rho_1}} \right)^{\rho_1} + \omega t^{\rho_1} \; .
	\end{equation}
	Moreover one can see that
	\begin{equation*}
		\forall \, \xi \in [0,s] \qquad 0 \leqslant \text{Arg}\left(\xi + te^{i \frac{\pi}{2 \rho_1}}\right) \leqslant \frac{\pi}{2 \rho_1} \; ;	\end{equation*}
	and since $\rho_1 \geqslant 1$, it follows
	\begin{equation*}
		0 \leqslant \text{Arg}\left( \left(\xi + te^{i \frac{\pi}{2 \rho_1}} \right)^{\rho_1-1} \right) \leqslant \frac{\pi(\rho_1-1)}{2 \rho_1} \leqslant \frac{\pi}{2} \; .
	\end{equation*}
	Hence the imaginary part of the complex number $\big(\xi + te^{i \frac{\pi}{2 \rho_1}}\big)^{\rho_1-1}$ is positive and so the real part of the right-hand side in \eqref{inteq} is negative. Therefore we have
	\begin{equation*}
		\Re{\big( i \omega z^{\rho_1} + \omega t^{\rho_1} \big)} \leqslant 0 \quad \Longrightarrow \quad \big| \, e^{i \omega z^{\rho_1}} \big| \, e^{\omega t^{\rho_1}} = \big| \, e^{i \omega z^{\rho_1} + \omega t^{\rho_1}} \big| = e^{\Re\left({i \omega z^{\rho_1} + \omega t^{\rho_1}} \right)} \leqslant 1 \; ,
	\end{equation*}
	which yields the result in the case $j=1$. To treat the case $j=2$, we use the following equality
	\begin{equation*}
		-i \rho_2 \, \omega \int_0^s \left( \xi + t e^{-i \frac{\pi}{2 \rho_2}} \right)^{\rho_2 -1} \, d\xi = -i \omega z^{\rho_2} + \omega t^{\rho_2} \; ,
	\end{equation*}
	and we carry out a similar work. This ends the proof.
\end{proof}

Now we compute the limit of a sequence of primitives of a certain function related to the following one,
\begin{equation*}
	s \in (0, s_j] \longmapsto s^{\mu_j - 1} e^{(-1)^{j+1} i \omega s^{\rho_j}} \in \C \; ,
\end{equation*}
which is the holomorphic factor appearing after Erdélyi's substitution. The sequence is constructed in such a way that each primitive is given by an integral on a finite path, and the sequence of these paths tends to the half-line considered in Lemma \ref{LEM3}, called $\Lambda^{(j)}(s)$. Exploiting the completeness of the space of holomorphic functions $\mathcal{H}(\Omega)$, where $\Omega \subset \C$ is an non-empty open subset of $\C$, and the continuity of the derivative in this space, we show that the resulting limit is also a primitive and its integration path is $\Lambda^{(j)}(s)$.\\
This result will permit to derive the first primitive of the above holomorphic factor, given in Corollary \ref{COR1}.

\begin{LEM4} \label{LEM4}
	Let $s_j > 0$. Define the parallelogram $D_j \subset \C$ and the domain $U \subset \C$ as follows :
	\begin{equation*}
		\begin{aligned}
			& \bullet \quad D_j := \left\{ v^* + t_ve^{(-1)^{j+1} i \frac{\pi}{2\rho_j}} \in \C \, \bigg| \, v^* \in (0, s_j + 1) \; , \; |t_v| < 1 \right\} \\
			& \bullet \quad U := \C \setminus \Big\{z \in \C \, \Big| \, \Re(z) \leqslant 0 \; , \; \Im(z) =0 \Big\}
		\end{aligned}
	\end{equation*}
	Fix $\mu_j \in (0,1]$, $\rho_j \geqslant 1$ and $n \in \N$. Let $F_{n,\omega}^{(j)}(.,.) : U \times \C \longrightarrow \C$ be the function defined by
	\begin{equation*}
		F_{n,\omega}^{(j)}(v,w) := \frac{(-1)^{n}}{n!} \, (v-w)^n \, v^{\mu_j - 1} e^{(-1)^{j+1} i \omega v^{\rho_j}} \; .
	\end{equation*}
	Then for every $w \in D_j$, $F_{n,\omega}^{(j)}(.,w)$ has a primitive $H_{n,\omega}^{(j)}(.,w)$ on $D_j$ given by
	\begin{equation*}
		H_{n,\omega}^{(j)}(v,w) := - \int_{\Lambda^{(j)}(v)} F_{n,\omega}^{(j)}(z,w) \, dz = \frac{(-1)^{n+1}}{n!} \int_{\Lambda^{(j)}(v)} (z-w)^n \, z^{\mu_j - 1} \, e^{(-1)^{j+1} i\omega z^{\rho_j}} dz \; ,
	\end{equation*}
	where
	\begin{equation} \label{lambda}
		\Lambda^{(j)}(v) := \left\{ v + t e^{(-1)^{j+1} i \frac{\pi}{2\rho_j}} \, \Big| \, t \geqslant 0 \right\} 
	\end{equation}
\end{LEM4}

\begin{proof}
	Let us fix $w \in D_1$ and $n \in \N$. Firstly, we show that the integral defining $H_{n,\omega}^{(1)}(v,w)$ is well-defined for every $v \in D_1$. Since $v \in D_1$, we can write $v= v^* + t_v e^{i\frac{\pi}{2\rho_1}}$ where $0 < v^* < s_1 + 1$ and $-1 < t_v < 1$, and we observe that
	\begin{equation*}
		-H_{n,\omega}^{(1)}(v,w) = \int_{\Lambda^{(1)}(v)} F_{n,\omega}^{(1)}(z,w) \, dz = \int_{\Lambda^{(1)}(v,v^*)} \dots \quad + \; \int_{\Lambda^{(1)}(v^*)} \dots \; ,
	\end{equation*}	
	where $\Lambda^{(1)}(v,v^*)$ is the segment which starts from the point $v$ and goes to $v^*$, and $\Lambda^{(1)}(v^*)$ is given by \eqref{lambda}. Since $F_{n,\omega}^{(1)}(.,w)$ is continuous on the segment $\Lambda^{(1)}(v,v^*)$, then the integral on $\Lambda^{(1)}(v,v^*)$ is well-defined. Concerning the second integral, we give a parametrization of the integration path $\Lambda^{(1)}(v^*)$,
	\begin{equation*}
		\forall \, t \in [0,+\infty) \qquad \lambda_{v^*}^{(1)}(t) := v^* + te^{i\frac{\pi}{2\rho_1}} \in \Lambda^{(1)}(v^*) \; .
	\end{equation*}
	We obtain
	\begin{align}
		\left| F_{n,\omega}^{(1)} \Big( \lambda_{v^*}^{(1)}(t), w \Big) \right|	&	\leqslant \frac{1}{n!} \left| v^* + t e^{i \frac{\pi}{2\rho_1}} - w \right|^n \left| v^* + t e^{i \frac{\pi}{2\rho_1}} \right|^{\mu_1-1} \left| e^{i \omega \left(v^* + t e^{i \frac{\pi}{2\rho_1}}\right)^{\rho_1}} \right| \nonumber \\
															&	\label{ineqrefv} \leqslant \frac{1}{n!} \sum_{k=0}^n \binom{n}{k} |v^* - w|^{n-k} \, (v^*)^{\mu_1-1} \, t^k \, e^{-\omega t^{\rho_1}} \; ,
	\end{align}
	where \eqref{ineqrefv} comes from the binomial Theorem, Lemma \ref{LEM3} and the geometric observation:
	\begin{equation*}
		\left| v^* + t e^{i \frac{\pi}{2\rho_1}} \right| \geqslant v^* \; .	\end{equation*}
	Since the right-hand side of \eqref{ineqrefv} defines an integrable function on $[0,+\infty)$, and since $\big| (\lambda_{v^*}^{(1)})'(t) \big| = 1$, the function $F_{n,\omega}^{(1)}(., w)$ is integrable on $\Lambda^{(1)}(v^*)$ and hence, $H_{n,\omega}^{(1)}(v,w)$ is well-defined for any $v \in D_1$.\\
	Secondly, we prove that $H_{n,\omega}^{(1)}(.,w) \, : \, D_1 \, \longrightarrow \, \C$ is a primitive of $F_{n,\omega}^{(1)}(.,w)$ on $D_1$. To this end, we show that $F_{n,\omega}^{(1)}(.,w)$ is a uniform limit on all compact subsets of $D_1$ of a sequence of functions $\big(H_{m,n,\omega}^{(1)}(.,w)\big)_{m\geqslant 1}$ which are primitives of $F_{n,\omega}^{(1)}(.,w)$ on $D_1$. Here we build this sequence of functions as follows: first of all, fix an arbitrary point $v_0 \geqslant s_1 + 1$, for instance $v_0 := s_1 + 1$, and define the following sequence of complex numbers :
	\begin{equation*}
		\forall \, m \in \N \backslash \{0\} \qquad v_m := v_0 + m \, e^{i\frac{\pi}{2 \rho_1}} \; .
	\end{equation*}	
	Let $m \in \N \backslash \{0\}$, let $v = v^* + t_v e^{i\frac{\pi}{2 \rho_1}} \in D_1$ and let $\Lambda_m(v)$ be the path which is the juxtaposition of the segment that starts from the point $v$ and goes to the point $v^* + me^{i\frac{\pi}{2\rho_1}}$ and of the horizontal segment that joins the points $v^* + me^{i\frac{\pi}{2\rho_1}}$ and $v_m$. We can now define the sequence of functions $\big(H_{m,n,\omega}^{(1)}(.,w) : D_1 \longrightarrow \C\big)_{m\geqslant1}$ as follows :
	\begin{equation*}
		H_{m,n,\omega}^{(1)}(v,w) := -\int_{\Lambda_m(v)} F_{n,\omega}^{(1)}(z,w) \,  dz \; .
	\end{equation*}
	It is clear that $F_{n,\omega}^{(1)}(.,w)$ is holomorphic on $U$, which is simply connected, and for any $v \in D_1$, $\Lambda_m(v)$ is included in $U$. The Cauchy integral Theorem affirms that each $H_{m,n,\omega}^{(1)}(.,w) : D_1 \longrightarrow \C$ is a primitive of the function $F_{n,\omega}^{(1)}(.,w)$.\\ Now we prove that this sequence converges to $H_{n,\omega}^{(1)}(.,w)$ uniformly on any compact subset $K$ of $D_1$. Let $K \subset D_1$ be a compact and for every $v \in K$, we have
	\begin{equation*}
		H_{m,n,\omega}^{(1)}(v,w) - H_{n,\omega}^{(1)}(v,w) = \int_{\Lambda_m^{c,1}(v)} F_{n,\omega}^{(1)}(z,w) \, dz + \int_{\Lambda_m^{c,2}(v)} F_{n,\omega}^{(1)}(z,w) \, dz \; ,
	\end{equation*}
	where $\Lambda_m^{c,1}(v)$ is the horizontal segment which starts from $v_m$ and goes to $v^* + m e^{i\frac{\pi}{2 \rho_1}}$, and $\Lambda_m^{c,2}(v)$ is the half-line with angle $\frac{\pi}{2 \rho_1}$ that starts from $v^* + m e^{i\frac{\pi}{2 \rho_1}}$ and goes to infinity. Let $\lambda_m^{c,1} : [0, v_0 - v^*] \longrightarrow \C$ and $\lambda_m^{c,2} : [0,+\infty) \longrightarrow \C$ be parametrizations of $\Lambda_m^{c,1}(v)$ and $\Lambda_m^{c,2}(v)$ respectively, and defined by
	\begin{equation*}
		\begin{aligned}
			& \bullet \quad \forall \, t \in \big[0, v_0-v^*\big] \qquad \lambda_m^{c,1}(t) := -t + v_0 + m e^{i\frac{\pi}{2\rho_1}} \in \Lambda_m^{c,1}(v) \; ,\\
			& \bullet \quad \forall \, t \in [0,+\infty) \qquad \lambda_m^{c,2}(t) := v^* + (t+m) e^{i\frac{\pi}{2\rho_1}} \in \Lambda_m^{c,2}(v) \; .
		\end{aligned}
	\end{equation*}
	Then we have the following estimates:
	\begin{align}
		\Big| F_{n,\omega}^{(1)}\big(\lambda_m^{c,1}(t),w\big) \Big|	& \label{inequa5} \leqslant \frac{1}{n!} \sum_{k=0}^n \binom{n}{k} \, | v_0 -w |^{n-k} \, \left| - t + m e^{i\frac{\pi}{2 \rho_1}}\right|^k \, m^{\mu_1 -1} \, e^{-\omega m^{\rho_1}} \\
																& \label{inequa6} \leqslant \frac{1}{n!} \sum_{k=0}^n \binom{n}{k} \, | v_0 -w |^{n-k} \, \big(C(K) + m \big)^k \, m^{\mu_1 -1} \, e^{-\omega m^{\rho_1}}
	\end{align}
	\begin{itemize}
		\item \eqref{inequa5}: use the binomial Theorem, Lemma \ref{LEM3} and $\big| \lambda_m^{c,1}(t) \big| \geqslant m$ ;
		\item \eqref{inequa6}: employ the compactness of $K$ which leads to $0 \leqslant t \leqslant v_0-v^* \leqslant C(K)$, for a certain constant $C(K) > 0$.
	\end{itemize}
	Inequality \eqref{inequa6} permits to estimate uniformly the integral on $\Lambda_m^{c,1}(v)$,
	\begin{equation*}
		\begin{aligned}
			\left| \int_{\Lambda_m^{c,1}(v)} F_{n,\omega}^{(1)}(z,w) \, dz \right|	& \leqslant \int_0^{|v_0-v^*|} \frac{1}{n!} \sum_{k=0}^n \binom{n}{k} \, | v_0 -w |^{n-k} \, \big(C(K) + m \big)^k \, m^{\mu_1 -1} \, e^{-\omega m^{\rho_1}} \, dt \\
		&\leqslant \frac{1}{n!} \sum_{k=0}^n \binom{n}{k} \, | v_0 -w |^{n-k} \, \big(C(K) + m \big)^k \, m^{\mu_1 -1} \, e^{-\omega m^{\rho_1}} C(K) \\
		& \longrightarrow \; 0 \quad , \quad m \longrightarrow +\infty \; ,
		\end{aligned}
	\end{equation*}
	where we used $0 \leqslant v_0-v^* \leqslant C(K)$ one more time. Here, the convergence is uniform with respect to $v$. Furthermore,
	\begin{align}
		\Big| F_{n,\omega}^{(1)}\big(\lambda_m^{c,2}(t),w\big) \Big|	& \label{inequa1} \leqslant  \frac{1}{n!} \sum_{k=0}^n \binom{n}{k} \, |v^* - w|^{n-k} \, |v^*|^{\mu_1-1} \, (t+m)^k \, e^{-\omega (t+m)^{\rho_1}} \\
																& \label{inequa2} \leqslant \frac{C_{n,w}(K)}{n!} \, \sum_{k=0}^n \binom{n}{k} \, (t+m)^k \, e^{-\omega (t+m)^{\rho_1}} \\
																& \label{inequa3} \leqslant \frac{C_{n,w}(K)}{n!} \, \sum_{k=0}^n \binom{n}{k} \, m^k \, e^{-\omega m^{\rho_1}} \, (1+t)^k \, e^{-\omega t^{\rho_1}} \\
																& \label{inequa4} \leqslant \frac{C_{n,w}(K) \, M_{k,\omega}}{n!} \, \sum_{k=0}^n \binom{n}{k} \, (1+t)^k \, e^{-\omega t^{\rho_1}}
	\end{align}
	\begin{itemize}
		\item \eqref{inequa1}: use the binomial Theorem, Lemma \ref{LEM3} and $\big| \lambda_m^{c,2}(t) \big| \geqslant v^*$ ;
		\item \eqref{inequa2}: use the compactness of $K$ and the fact that $v \in K$ to bound uniformly $|v^* - w|$ and $|v^*|$ ;
		\item \eqref{inequa3}: use the inequalities $(m+t)^k \leqslant m^k (1+t)^k$ and $e^{-\omega(t+m)^{\rho_1}} \leqslant e^{-\omega m^{\rho_1}} e^{-\omega t^{\rho_1}}$ ;
		\item \eqref{inequa4}: use the boundedness of the sequences $\big(m^k e^{-\omega m^{\rho_1}}\big)_{m\geqslant1}$ for $k=0,\ldots,n$.
	\end{itemize}
	We remark that \eqref{inequa3} tends to $0$ as $m$ tends to infinity for all $t \geqslant 0$ and \eqref{inequa4} gives an integrable function independent on $m$. So by the dominated convergence Theorem,
	\begin{align}
			\bigg| \int_{\Lambda_m^{c,2}(v)} F_{n,\omega}^{(1)}(z,w) \, dz \bigg|	& \label{inequa7} \leqslant \int_0^{+\infty} \frac{C_{n,w}(K)}{n!} \, \sum_{k=0}^n \binom{n}{k} \, m^k \, e^{-\omega m^{\rho_1}} \, (1+t)^k \, e^{-\omega t^{\rho_1}} \, dt \\
																					& \longrightarrow \; 0 \quad , \quad m \longrightarrow + \infty \; , \nonumber
	\end{align}
	and the convergence is uniform with respect to $v$ since the right-hand side term in \eqref{inequa7} is independent from $v$. Hence the hypotheses of a theorem of Weierstrass are satisfied and therefore the function $H_{n,\omega}^{(1)}(.,w): D_1 \longrightarrow \C$ is holomorphic and its derivative is given by
	\begin{equation*}
		\forall \, v \in D_1 \qquad \frac{\partial}{\partial v} \big[ H_{n,\omega}^{(1)} \big](v,w) = \lim_{n \rightarrow + \infty} \frac{\partial}{\partial v} \big[ H_{m,n,\omega}^{(1)}(v,w) \big] = F_{n,\omega}^{(1)}(v, w) \; ,
	\end{equation*}
	and the convergence is uniform on every compact subset.
\end{proof}

Since we integrate by parts many times in Erdélyi's proof, we need successive primitives of the holomorphic part of the integrand. For this purpose, we establish this second intermediate but essential result by employing the preceding theorem as well as complex analysis in several variables.\\
The desired primitives will be deduced from the next result in Corollary \ref{COR1}.

\begin{LEM5} \label{LEM5}
	Let $n \in \N \backslash \{0\}$, define the function $h : \C \longrightarrow \C \times \C$ by
	\begin{equation*}
		h(u) := (u,u) \; ,
	\end{equation*}
	and let $H_{n,\omega}^{(j)}(.,.) : D_j \times D_j \longrightarrow \C$ be the function defined in Theorem \ref{LEM4}. Then the composite function $H_{n,\omega}^{(j)}(.,.) \circ h$ is holomorphic on $D_j$ and its derivative is given by
	\begin{equation*}
		\forall \, u \in D_j \qquad \frac{d}{d u} \hspace{-1mm} \left[H_{n,\omega}^{(j)}\circ h \right] \hspace{-1mm} (u) = \left(H_{n-1,\omega}^{(j)} \circ h \right) \hspace{-1mm} (u) \; .
	\end{equation*}
\end{LEM5}

\begin{proof}
	The aim of the proof is to differentiate the composite function. For this purpose, we must ensure that this function is holomorphic with respect to each variable.\\
	Fix $n \in \N \backslash \{0\}$. We remark that each component of $h$ is holomorphic on $\C$, so is $h$ on $\C \times \C$. Moreover for any fixed $w \in  D_1$, $H_{n,\omega}^{(1)}(.,w) : D_1 \longrightarrow \C$ is a primitive of $F_{n,\omega}^{(1)}(.,w)$ on $D_1$ by Theorem \ref{LEM4}, so it is holomorphic. Now let us show that $H_{n,\omega}^{(1)}(v,.) : D_1 \longrightarrow \C$ belongs to $\mathcal{C}^1(D_1)$ and satisfies the Cauchy-Riemann equations for fixed $v \in  D_1$. To do so, we employ the holomorphy of $F_{n,\omega}^{(1)}(v,.) : \C \longrightarrow \C$ which provides the following relations :
	\begin{equation} \label{holomorphy}
		\forall \, w = x + iy \in \C \qquad  \frac{\partial}{\partial w} \hspace{-1mm} \left[ F_{n,\omega}^{(1)} \right] \hspace{-1mm} (v,w) = \frac{\partial}{\partial x} \hspace{-1mm} \left[ F_{n,\omega}^{(1)} \right] \hspace{-1mm} (v,w) = -i \frac{\partial}{\partial y} \hspace{-1mm} \left[ F_{n,\omega}^{(1)} \right] \hspace{-1mm} (v,w) \; .
	\end{equation}
	And by a quick calculation, we obtain
	\begin{equation} \label{form1}
		\frac{\partial}{\partial w} \hspace{-1mm} \left[ F_{n,\omega}^{(1)} \right] \hspace{-1mm} (v,w) = \frac{(-1)^{n-1}}{(n-1)!} \, (v-w)^{n-1} v^{\mu_1-1} e^{i \omega v^{\rho_1}} = F_{n-1,\omega}^{(1)}(v,w) \; .
	\end{equation}
	Furthermore, one can bound $F_{n-1,\omega}^{(1)}(.,w)$ on each path $\Lambda^{(1)}(v,v^*)$ and $\Lambda^{(1)}(v^*)$ by integrable functions independent on $w$. To do so, one can parametrize each path $\Lambda^{(1)}(v,v^*)$ and $\Lambda^{(1)}(v^*)$ and employ similar arguments to the previous ones as well as the boundedness of $D_1$. So we obtain the ability to differentiate under the integral sign which yields the following equalities :
	\begin{align}
		-\frac{\partial}{\partial x} \hspace{-1mm} \left[ H_{n,\omega}^{(1)} \right] (v,w)	& = \frac{\partial}{\partial x} \hspace{-1mm} \left[ \int_{\Lambda^{(1)}(v,v^*)} F_{n,\omega}^{(1)}(z,w) \, dz \right] + \frac{\partial}{\partial x} \hspace{-1mm} \left[ \int_{\Lambda^{(1)}(v^*)} F_{n,\omega}^{(1)}(z,w) \, dz \right] \nonumber \\
												& = \label{equa1} \int_{\Lambda^{(1)}(v,v^*)} \frac{\partial}{\partial x} \hspace{-1mm} \left[ F_{n,\omega}^{(1)} \right] \hspace{-1mm} (z,w) \, dz \: + \int_{\Lambda^{(1)}(v^*)} \frac{\partial}{\partial x} \hspace{-1mm} \left[ F_{n,\omega}^{(1)} \right] \hspace{-1mm} (z,w) \, dz \\
												& = \label{equa2} \int_{\Lambda^{(1)}(v,v^*)} \frac{\partial}{\partial w} \hspace{-1mm} \left[ F_{n,\omega}^{(1)} \right] \hspace{-1mm} (z,w) \, dz \: + \int_{\Lambda^{(1)}(v^*)} \frac{\partial}{\partial w} \hspace{-1mm} \left[ F_{n,\omega}^{(1)} \right] \hspace{-1mm} (z,w) \, dz \\
												& = \label{equa3} \int_{\Lambda^{(1)}(v,v^*)} F_{n-1,\omega}^{(1)}(z,w) \, dz \: + \int_{\Lambda^{(1)}(v^*)} F_{n-1,\omega}^{(1)} (z,w) \, dz \\
												& = \int_{\Lambda^{(1)}(v)} F_{n-1,\omega}^{(1)}(z,w) \, dz \nonumber \\
												& = -H_{n-1,\omega}^{(1)}(v,w) \nonumber												
	\end{align}
	\begin{itemize}
		\item \eqref{equa1}: apply of the theorem of differentiation under the integral sign ;
		\item \eqref{equa2}: use equalities \eqref{holomorphy} coming from the holomorphy of the function $F_{n,\omega}^{(1)}(v,.)$ ;
		\item \eqref{equa3}: employ relation \eqref{form1} .
	\end{itemize}
	In a similar way, we obtain
	\begin{equation*}
		-i \frac{\partial}{\partial y} \hspace{-1mm} \left[ H_{n,\omega}^{(1)} \right] \hspace{-1mm} (v,w) = H_{n-1,\omega}^{(1)}(v,w) \; .
	\end{equation*}
	Then the Cauchy-Riemann quations are satisified and $\displaystyle \frac{\partial}{\partial x} \big[ H_{n,\omega}^{(1)} \big] (v,.)$ and $\displaystyle \frac{\partial}{\partial y} \big[ H_{n,\omega}^{(1)} \big](v,.)$ are continuous on $D_1$ by the continuity of $F_{n-1,\omega}^{(1)}(z,.) : \C \longrightarrow \C$. So $H_{n,\omega}^{(1)}(v,.) : D_1 \longrightarrow \C$ is holomorphic, with
	\begin{equation*}
		\frac{\partial}{\partial w} \hspace{-1mm} \left[ H_{n,\omega}^{(1)} \right] \hspace{-1mm} (v,w) = H_{n-1,\omega}^{(1)}(v,w) \; .
	\end{equation*}
	Finally the composite function $H_{n,\omega}^{(1)} \circ h$ is holomorphic on $D_1 \times D_1$ and we have the formula
	\begin{align*}
		\frac{d}{d u} \big[H_{n,\omega}^{(1)} \circ h \big](u)	& = \left( \frac{\partial}{\partial v } \big[ H_{n,\omega}^{(1)} \big] \big( h(u) \big) \quad \frac{\partial}{\partial w} \big[ H_{n,\omega}^{(1)} \big] \big(h(u) \big) \right) \left( \begin{array}{c}
						1 \\
						1 \\
						\end{array}
						\right) \\
																& =  \frac{\partial}{\partial v} \big[ H_{n,\omega}^{(1)} \big] (u,u) +  \frac{\partial}{\partial w} \big[ H_{n,\omega}^{(1)} \big] (u,u) \; ;
	\end{align*}
	And a short computation shows that $\displaystyle \frac{\partial}{\partial v} \big[ H_{n,\omega}^{(1)} \big] (u,u) = F_{n,w}^{(1)}(u,u) = 0$, so for all $u \in D_1$,
	\begin{equation*}
		\frac{d}{d u} \big[H_{n,\omega}^{(1)}\circ h \big](u) = \left(H_{n-1,\omega}^{(1)} \circ h \right)(u) = \frac{(-1)^n}{(n-1)!} \int_{\Lambda^{(1)}(u)} (z-u)^{n-1} z^{\mu_1-1} e^{i \omega z^{\rho_1}} dz \; .
	\end{equation*}
\end{proof}

Finally by restricting the domain of definition of the functions introduced in the two preceding theorems to the interval $(0,s_j]$, we derive formulas for the primitives of the function $s \in (0,s_j] \longmapsto s^{\mu_j-1} \, e^{(-1)^{j+1} i \omega s^{\rho_j}}$, the holomorphic part of the integrand.

\begin{COR1} \label{COR1}
	Fix $s_j > 0$, $\rho_j \geqslant 1$ and $\mu_j \in (0,1]$. For any $\omega >0$, the sequence of functions $\big( \phi_{n}^{(j)}(.,\omega,\rho_j,\mu_j) : \, (0,s_j] \longrightarrow \C \big)_{n \geqslant 1}$ defined in Theorem \ref{THM1} or in Theorem \ref{THM2} satisfies the recursive relation :
	\begin{equation*}
		\forall \, s \in \; (0,s_j] \qquad \left\{	\hspace{-0.3cm} \begin{array}{rl}
														& \displaystyle \frac{\partial}{\partial s} \hspace{-1mm} \left[ \phi_{n+1}^{(j)} \right] \hspace{-1mm} (s,\omega,\rho_j,\mu_j) = \phi_{n}^{(j)}(s, \omega,\rho_j,\mu_j) \qquad \forall \, n \geqslant 1 \; , \\[3mm]
														& \displaystyle \frac{\partial}{\partial s} \hspace{-1mm} \left[ \phi_{1}^{(j)} \right] \hspace{-1mm} (s,\omega,\rho_j,\mu_j) = s^{\mu_j-1} \, e^{(-1)^{j+1} i \omega s^{\rho_j}} \; .
		\end{array} \right.
	\end{equation*}
\end{COR1}

\begin{proof}
	It suffices to note that $\phi_{n+1}^{(j)}(.,\omega,\rho_j,\mu_j)$ is the restriction to $(0,s_j] \subset D_j$ of the function $H_{n,\omega}^{(j)}\circ h$. Hence Theorem \ref{LEM4} affirms that $\phi_{1}^{(j)}(.,\omega,\rho_j,\mu_j) : (0,s_j] \longrightarrow \C$ is a primitive of $s \in (0,s_j] \longmapsto s^{\mu_j-1} e^{(-1)^{j+1} i \omega s^{\rho_j}}$, and use Theorem \ref{LEM5} to show that a primitive of $\phi_{n}^{(j)}(.,\omega,\rho_j,\mu_j) : (0,s_j] \longrightarrow \C$ is given by $\phi_{n+1}^{(j)}(.,\omega,\rho_j,\mu_j) : (0,s_j] \longrightarrow \C$, for $n \geqslant 1$.
\end{proof}

\begin{REM2} \label{REM2}
	 \emph{The function $\phi_{n+1}^{(j)}(.,\omega,\rho_j,\mu_j) : \, (0,s_j] \longrightarrow \C$ can be extended to $[0,+\infty)$. Indeed, we recall a parametrization of the curve $\Lambda^{(j)}(s)$ given by		\begin{equation*}
			\lambda_{s}^{(j)} : t \in (0,+\infty) \longmapsto s+t e^{(-1)^{j+1} i\frac{\pi}{2 \rho_j}} \in \Lambda^{(j)}(s) \; ,
		\end{equation*}
		and we consider the following estimate
		\begin{equation} \label{eqrefe}
			\forall \,  t > 0 \qquad \Big|F_{n,\omega}^{(j)} \big( \lambda_{s}^{(j)}(t) , s \big) \Big| \leqslant \frac{1}{n!} \, t^{n+\mu_j-1} \, e^{-\omega t^{\rho_j}} \; ,
		\end{equation}
		which was obtained by noting that
		\begin{equation*}
			t \leqslant \left| s+t e^{(-1)^{j+1} i\frac{\pi}{2 \rho_j}} \right| = \left| \lambda_s^{(j)}(t) \right| \qquad \Longrightarrow \qquad t^{\mu_j-1} \geqslant \left| \lambda_s^{(j)}(t) \right|^{\mu_j-1} \; .
		\end{equation*}
		We notice that the right-hand side of \eqref{eqrefe} is an integrable function with respect to $t$ on $[0,+\infty)$ and independent from $s > 0$. So $\phi_{n+1}^{(j)}(s,\omega,\rho_j,\mu_j)$ is well-defined for all $s \geqslant 0$. In particular, $\phi_{n+1}^{(j)}(0,\omega,\rho_j,\mu_j)$ is defined as follows,
		\begin{align*}
				\phi_{n+1}^{(j)}(0,\omega,\rho_j,\mu_j)	& := \lim_{s \rightarrow 0^+} \phi_{n+1}^{(j)}(s,\omega,\rho_j,\mu_j) \\
															& = \frac{(-1)^{n+1}}{n!} \int_{\Lambda^{(j)}(0)} z^{n+\mu_j -1} e^{(-1)^{j+1} i\omega z^{\rho_j}} dz \; .
			\end{align*}
			}
\end{REM2}

\end{document}